\documentclass[a4paper,article,10pt]{memoir}
\usepackage{CCLAuthor}
\usepackage[reqno]{amsmath}
\usepackage{amssymb}
\usepackage{amsthm}  
\theoremstyle{plain}

\theoremstyle{definition}

\numberwithin{theorem}{chapter}
\usepackage[round,authoryear]{natbib} 
\usepackage[]{graphicx}               
\usepackage{bm}
\usepackage{esvect}
\usepackage{ulem}
\usepackage{adjustbox}
\usepackage{accents}
\usepackage{stmaryrd}
\usepackage{xcolor}
\usepackage{enumitem}
\usepackage{hyperref}
\usepackage{cleveref}
\usepackage{booktabs}
\usepackage{longtable}
\usepackage[section]{placeins}
\usepackage{tabu}

\definecolor{chcol}{rgb}{0.4,0.,0.9}

\renewcommand{\Pr}{\text{Pr}}     
\newcommand{\halfComma}{\kern 0.083334em}
\newcommand{\pderivative}[2]{\frac{\partial #1}{\partial #2}}
\newcommand{\avg}[1]{\left\{\hspace*{-3pt}\left\{#1\right\}\hspace*{-3pt}\right\}}
\newcommand{\jump}[1]{\ensuremath{\left\llbracket #1 \right\rrbracket}}
\newcommand{\dS}{{\,\operatorname{dS}}}         
\newcommand{\PBT}{{\,\operatorname{PBT}}}    
\newcommand{\ec}{{\mathrm{EC}}}                     
\newcommand{\ent}{{S}}                                     
\newcommand{\dx}{\,\text{d}x}                           
\newcommand\iprod[1]{\left\langle #1\right\rangle}                                             
\newcommand\inorm[1]{\left |\left| #1\right|\right|}                                               
\newcommand\iprodN[1]{\left\langle #1\right\rangle_{\!N}}                                 
\newcommand\inormN[1]{\left |\left| #1\right|\right|_{N}}                                     
\newcommand\irefInt{\int\limits_{-1}^{1}}                                                            
\newcommand\ivolN{\int\limits_{ N }\! }                                                              
\newcommand\volEN{\mkern-11mu\int\limits_{E , N}\mkern-5mu }                    
\newcommand\isurfEN{\mkern-11mu\int\limits_{\partial E , N}\mkern-11mu }    
\renewcommand\vec[1]{\accentset{\,\rightarrow}{#1}}                              
\newcommand\spacevec[1]{\accentset{\,\rightarrow}{#1}}                        
\newcommand\contravec[1]{\tilde{ #1}}                                                     
\newcommand\contraspacevec[1]{\spacevec{\tilde{#1}}}                          
\newcommand\statevec[1]{\mathbf #1}                                                     
\newcommand\statevecGreek[1]{\boldsymbol #1}                                     
\newcommand\contrastatevec[1]{\tilde{\mathbf #1}}                                 
\newcommand\acclrvec[1]{\accentset{\,\leftrightarrow}{#1}}                      
\newcommand\bigstatevec[1]{\acclrvec{{\mathbf #1}}}                              
\newcommand\biggreekstatevec[1]{\acclrvec{\boldsymbol #1}}                
\newcommand\bigcontravec[1]{\acclrvec{\tilde{\mathbf{#1}}}}                   
\newcommand\vecNabla{\accentset{\,\rightarrow}\nabla}                         
\newcommand\vecNablaXi{\accentset{\,\rightarrow}\nabla_{\!\xi}}            
\newcommand\vecNablaX{\accentset{\,\rightarrow}\nabla_{\!x}}              
\newcommand\mmatrix[1]{\underbar{#1}}				
\newcommand\matrixvec[1]{\mathcal #1}                           
\newcommand\bigmatrix[1]{\mathfrak #1}                          
\newcommand{\dmat}{\matrixvec{D}}     
\newcommand{\qmat}{\matrixvec{Q}}    
\newcommand{\mmat}{\matrixvec{M}}   
\newcommand{\bmat}{\matrixvec{B}}    
\newcommand\interiorfaces{\genfrac{}{}{0pt}{}{\mathrm{interior}}{\mathrm{faces}}}
\newcommand\boundaryfaces{\genfrac{}{}{0pt}{}{\mathrm{boundary}}{\mathrm{faces}}}
\newcommand\allfaces{\genfrac{}{}{0pt}{}{\mathrm{all}}{\mathrm{faces}}}
\newcommand{\testfuncOne}{\statevecGreek{\varphi}}  
\newcommand{\testfuncTwo}{\biggreekstatevec{\vartheta}}

\newcommand{\DD}{\spacevec{\mathbb{D}}\cdot}
\newcommand{\DDs}{\spacevec{\mathbb{D}}^{s}\cdot}
\newcommand{\IN}[1]{\mathbb I^{N}\!\!\left(#1\right)} 
\newcommand{\PN}[1]{\mathbb P^{#1}}
\newcommand{\LTwo}[1]{\mathbb L^{2}\!\left(#1\right)}
\newcommand\overRe{\frac{1}{{\operatorname{Re} }}}
\newcommand\twooverRe{\frac{2}{{\operatorname{Re} }}}
\newcommand\oneHalf{\frac{1}{2}}

\usepackage{xargs}    
\usepackage[colorinlistoftodos,prependcaption,textsize=tiny]{todonotes}
\newcommandx{\unsure}[2][1=]{\todo[linecolor=blue,backgroundcolor=blue!25,bordercolor=blue,#1]{#2}}
\newcommandx{\changeThis}[2][1=]{\todo[linecolor=red,backgroundcolor=red!25,bordercolor=red,#1]{#2}}
\newcommandx{\info}[2][1=]{\todo[linecolor=OliveGreen,backgroundcolor=OliveGreen!25,bordercolor=OliveGreen,#1]{#2}}
\newcommandx{\improvement}[2][1=]{\todo[linecolor=gray,backgroundcolor=gray!25,bordercolor=gray,#1]{#2}}
\newcommandx{\thiswillnotshow}[2][1=]{\todo[disable,#1]{#2}}
\begin{document}
%
%
%
\title{Construction of Modern Robust Nodal Discontinuous Galerkin Spectral Element Methods for the Compressible Navier-Stokes Equations}
%
%
\author{%
    Andrew R. Winters\textsuperscript{\ddag}
    David A. Kopriva\textsuperscript{\dag},
    Gregor J. Gassner\textsuperscript{*},
    Florian Hindenlang\textsuperscript{**},
    \\ \smallskip\small
    \textsuperscript{\ddag}
    Link\"{o}ping University, Link\"{o}ping, Sweden
    \\
    \textsuperscript{\dag} 
    The Florida State University, FL, USA and \\
    San Diego State University, San Diego, CA, USA
    \\
    \textsuperscript{*} 
    University of Cologne, Cologne, Germany                     
    \\
    \textsuperscript{**}
    Max Planck Institute for Plasma Physics, Garching, Germany
    }
    \maketitle
%
%
%
%
    \begin{abstract}
Discontinuous Galerkin (DG) methods have a long history in computational physics and engineering to approximate solutions of partial differential equations due to their high-order accuracy and geometric flexibility. However, DG is not perfect and there remain some issues. Concerning robustness, DG has undergone an extensive transformation over the past seven years into its modern form that provides statements on solution boundedness for linear and nonlinear problems. 

This chapter takes a constructive approach to introduce a modern incarnation of the DG spectral element method for the compressible Navier-Stokes equations in a three-dimensional curvilinear context. The groundwork of the numerical scheme comes from classic principles of spectral methods including polynomial approximations and Gauss-type quadratures. We identify aliasing as one underlying cause of the robustness issues for classical DG spectral methods. Removing said aliasing errors requires a particular differentiation matrix and careful discretization of the advective flux terms in the governing equations.
%
%
%
%
%
%
    \end{abstract}
    \newpage
%
%
\begin{KeepFromToc}
   \tableofcontents
\end{KeepFromToc}
%
%
\CCLsection{Prologue}

The discontinuous Galerkin (DG) method dates back to the work of \citet{Nitsche} for the solution of elliptic problems and to the work of \citet{ReedHill} for the solution of linear hyperbolic advection problems. However, it was the work by Cockburn, Shu and others starting in 1989, e.g., \citep{CockburnHou&Shu,Cockburn1998199,cockburn1991,CBS-convection-diffusion}, that initiated a substantial interest in DG methods for the approximation of nonlinear hyperbolic and mixed hyperbolic-parabolic conservation laws. \citet{Bassi&Rebay:1997:B&F97} were among the first who applied the DG method to the compressible Navier-Stokes equations. 

As time went on, the DG methodology gained more and more traction in many different applications fields, such as e.g. compressible flows \cite{Black1999,black2000,Rasetarinera:2001zr}, electromagnetics and optics \cite{Koprivaetal1999,Koprivaetal2000,ISI:000246595900002,ISI:000226090600009}, acoustics \cite{chan2017,Rasetarineraetal2001,Stanescu:2002tg,Stanescu:2002kl,wilcox2010}, meteorology \cite{giraldo2002,Giraldo:2008eu,Restelli:2009nx,bonov2018}, and geophysics \cite{GroundwaterDSEM,DamBreak2004}. Nowadays, DG is applied in almost all sciences where high fidelity computational approximations of differential equations is necessary. 

The first available book on DG methods was published in 1999 \citet{Cockburn1999}. This book was a collection of proceedings articles and hence still left many practical issues related to an actual implementation unanswered. The situation changed however in 2005, when \citet{Karniadakis:2005fj} released their book, which not only includes the theory but also provided detailed explanations of the scheme and the algorithms. 
Focus of that early work was mostly on the modal variant of the DG scheme on hybrid mixed meshes. 

Whereas the mathematical formulation of the discontinuous Galerkin scheme is agnostic to the particular choice of element types and basis function type, the actual scheme, i.e. the algorithms and the numerical properties such as efficiency and accuracy depend strongly on the choice of basis functions (and many other choices such as the type of discrete integration and type of element shapes). Later, \citet{hesthaven2008} published a DG book with focus on nodal basis function on simplex shaped elements, while \citet{Kopriva:2009nx} published a book on nodal DG methods on quadrilateral (and hexahedral) elements. These three books cover the vast majority of commonly used DG variants and somewhat represented the state of the art, at least at their publication times. 

The textbooks together with the promising theoretical properties of DG methods, such as e.g. high dispersion accuracy and low dissipation errors, e.g., \citep{Hu:1999,Gassner:2010pd}, high potential for parallel computing, e.g., \citep{altmann2013,Stanescu:2002tg,baggag2000},  and natural stabilisation for advection dominated problems via (approximate) Riemann solvers, certainly helped to attract more and more interest in the application of DG in research and industry.

However, it became clear that additional development and advancement in the state of the art was needed to make DG methods really competitive, as for instance documented in large European research collaborations ADIGMA\footnote{https://cordis.europa.eu/docs/results/30/30719/122807181-6\_en.pdf} and IDIHOM\footnote{https://www.dlr.de/as/desktopdefault.aspx/tabid-7027/11654\_read-27492/}. The main issues that hold back high-order DG are (i) efficiency of time integration, (ii)  high-order grid generation and (iii) robustness. See e.g. \cite{Wang2013} for more discussion and details. The authors also identify efficient $hp$-adaptivity as a major issue, however control of adaptivity and error estimates is a general problem for all numerical schemes, not specific to high-order DG.   

Discontinuous Galerkin methods are in general very well suited for explicit time stepping, e.g., low storage Runge-Kutta time integrators. As the mass matrix is local and most often diagonal, no inversion is needed so that the computational complexity of a single time step is very low. Typical for explicit time integration, a CFL-like time step restriction is necessary to keep the simulations stable. For DG, the maximum CFL number not only depends on the mesh size and the fastest wave speed (or equivalent viscous speed), but also on the polynomial degree of the approximation. In industrial applications, such an explicit time step restriction could turn out to be prohibitive, resulting in an inefficient overall framework. Naturally, a remedy to this issue are implicit time integration methods, in particular implicit high-order Runge-Kutta methods. However, as it turns out, high-order DG methods in three spatial dimensions result in huge block dense algebraic equation systems that are notoriously difficult to solve efficiently. And without a proper pre-conditioner, no benefit from implicit time integration is left over. (In many cases, one can observe a negative speed up compared to explicit Runge-Kutta time integration). 

High-order grid generation is a subtle problem not easy to realise at first. Often high-order DG methods are praised for their capabilities of using unstructured meshes. A major issue, however, is that gaining the full benefit of the high-order DG approach requires meshes with high-order approximations at curved boundaries. The automatic construction of curvilinear high quality meshes, with e.g. boundary layers, is still an open research problem. There are strategies and work arounds available, e.g. an open source software solution that post processes a straight-sided mesh\footnote{https://www.hopr-project.org}. But for a fully operable process chain, much more research and development is needed. 

The third major issue, that of robustness, is more subtle than the time integration efficiency and high-order mesh generation, but at least as important, if not more important -- without a stable approximation, there is no need to talk about efficiency at all. Also no high-order mesh is a remedy when the scheme is not robust and blows up during the simulation. DG methods have natural inbuilt upwind-like dissipation from (approximate) Riemann solvers and are often characterized as more stable than, for instance, their continuous counterpart, the standard Galerkin finite element methods. 

However, with the high-order of DG methods there comes a whole bag of robustness issues. Most prominently, when simulating problems with discontinuities such as shock waves, Gibbs-type oscillations occur along with possible violations of physical solution bounds, such as positivity of density and pressure, that lead to failure of the algorithm. The inbuilt upwind-like dissipation in the DG methodology is not enough for stability. 

Even without discontinuities, other nonlinear solution features may bring high-order DG methods to failure. In realistic applications of turbulent flows, the resolution is orders of magnitude smaller than necessary for grid convergence. While most DG variants work perfectly fine for well resolved problems, e.g., in a grid convergence study, it turns out to be non-trivial to construct robust DG schemes for underresolved problems. 

Not all of the issues presented above have been solved to a satisfactory level for the DG method; in fact most of these issues are active areas of research. A comprehensive discussion on all aspects and possible solution strategies would clearly go beyond the scope of a single book chapter, especially since we aim to present  the details of the theory in the spirit of the three DG books mentioned above, i.e. with maximum detail. 

In comparison to the three DG books that essentially cover the prior state of the art, we focus on the issue of robustness for underresolved flows, e.g., when simulating turbulence, and aim to present the advances made in the last decade. We first introduce the mathematical and algorithmic building blocks of the DG method in Section \ref{sec:spectral_calc}. In Section \ref{sec:section_NSE_GG}, we present the underlying partial differential equations for compressible fluid dynamics, with a focus on their mathematical stability properties. In Section \ref{sec:section_SEM_GG}, we provide a detailed description of the spectral element framework on curvilinear unstructured hexahedral grids and its corresponding DG variant in Section \ref{sec:curvi_DG}. The main strategy to get a provably stable nodal DG scheme is presented in Section \ref{sec:section_splitformDG_GG} with an outlook and possible extensions discussed in the last Section \ref{sec:epilogue_GG}.

\CCLsubsubsection{Nomenclature}

Notation used throughout this chapter is adapted from \citet{Gassner_BR1}:

{\centering
\begin{longtable}[c]{ll}
$\PN{N}$ & Space of polynomials of degree $\leqslant N$\\[0.075cm]
$\mathbb{I}^{N}$& Polynomial interpolation operator\\[0.075cm]
$(x,y,z)$ &Physical space coordinates\\[0.075cm]
$\left(\xi,\eta,\zeta\right)$ &Reference space coordinates\\[0.075cm]
$\spacevec v $ & Vector in three-dimensional space\\[0.075cm]
$\spacevec n= n_{1}\hat x + n_{2}\hat y + n_{3}\hat z$ & Physical space normal vector\\[0.075cm]
$\hat n= \hat n^{1}\hat \xi + \hat n^2\hat \eta + \hat n^3\hat \zeta$ & Reference space normal vector\\[0.075cm]
$\statevec u$ &Continuous quantity\\[0.075cm]
$\statevec {U}$ & Polynomial approximation\\[0.075cm]
$\bigstatevec f,\, \bigcontravec{\!f}$& Block vector of Cartesian flux and contravariant flux\\[0.075cm]
$\dmat$& $(N+1)\times (N+1)$ Matrix\\[0.075cm]
$\mmatrix B$& $5\times 5$ Matrix\\[0.075cm]
$\bigmatrix B$& $15\times 15$ Block matrix\\[0.075cm]
\end{longtable}
\setcounter{table}{0} 
}

\CCLsection{Spectral Calculus Toolbox}\label{sec:spectral_calc}

Algorithms developed to numerically model physical problems are typically designed to solve discrete approximations of partial differential equations (PDEs). Often, such PDEs are derived in the framework of differential and integral calculus, and can be formulated in terms of first-order differential operators, such as the divergence. These PDEs express fundamental physical laws like the conservation of mass, momentum, and total energy. 

The continuum operators are equipped with important differential and integral identities, e.g., the derivative of a constant function is zero, or integration-by-parts.  To demonstrate conservation and stability of numerical approximations as well as to accurately capture the physics of a solution it is beneficial for the discretization to mimic as many of these important differential and integral identities as possible.

This section provides the groundwork and discussion of a discrete spectral calculus for nodal spectral element methods. The basic principles of polynomial interpolation and high-order Gauss-type quadrature provide the tools needed to develop high-order, conservative, and stable approximations for PDEs written as conservation laws or balance laws. 

Spectral methods owe their roots to the solutions of PDEs by orthogonal polynomial expansions \citep{gottlieborszag77,Kreiss&Oliger:1973}. Spectral element methods today
approximate solutions of PDEs with piecewise polynomials equivalent to finite series of Legendre polynomials, e.g.
\begin{equation}
U(x) = \sum\limits_{k = 0}^N {{{\hat C}_k}{L_k}(x)},
\label{eq:LegendrePolynomialInterpolant_DAK}
\end{equation}
which is a polynomial of degree $N$. As a shorthand, we let $\PN{N}$ denote the space of polynomials of degree less than or equal to $N$ and write $U\in\PN{N}$.
Whereas finite difference methods approximate the solutions of PDEs only at a finite set of discrete points, $U_{j}$, spectral element methods
are akin to finite element methods in that the approximation is well-defined at all points.

Four features characterize Discontinuous Galerkin (DG) spectral element methods:
\begin{itemize}
\item Approximation of the solutions and fluxes by piecewise high-order polynomials that represent polynomial expansions.
\item Approximation of integrals and inner products by high-order Gauss-type quadratures.
\item A weak formulation of the original differential equations.
\item Coupling of elements through the use of a numerical flux (aka ``Riemann Solver'').
\end{itemize}

In this section, we review the background for the first two features, namely the approximation of functions by Legendre polynomial expansions and Gauss-Lobatto quadratures used by discontinuous Galerkin
spectral element methods. Along the way, we develop a discrete calculus that mirrors the continuous one, which will allow us to use familiar operations on discretely defined functions.

\CCLsubsection{Legendre Polynomials and Series}

That the Legendre polynomial $L_{k}(x)$ in \eqref{eq:LegendrePolynomialInterpolant_DAK} is a polynomial
of degree $k$ is seen through the three term recurrence relation it satisfies,
 \begin{equation}
 \begin{gathered}
  {L_{k + 1}}(x) = \frac{{2k + 1}}{{k + 1}}x{L_k}(x) - \frac{k}{{k + 1}}{L_{k - 1}}(x), \quad x\in[-1,1],\hfill \\
  {L_0}(x) = 1, \hfill \\
  {L_1}(x) = x. \hfill \\ 
\end{gathered}
  \end{equation}
Two Legendre polynomials $L_{k}$ and $L_{n}$ are orthogonal with respect to the $\LTwo{-1,1}$ inner product
\begin{equation}
\iprod {{L_k},{L_n}}  = \irefInt {{L_k}{L_n}\dx}  = \frac{2}{{2k + 1}}{\delta _{kn}},
\label{eq:L2InnerProductOfL_DAK}
\end{equation}
where $\LTwo{-1,1}$ is the space of square integrable functions on the interval $[-1,1]$. That is, all functions $u$ for which ${\left\| u \right\|^2} = \iprod{u,u} < \infty$. The Kronecker delta is nonzero only when the subscripts match
\begin{equation}{\delta _{kn}} = \left\{ \begin{gathered}
  1\quad k = n \hfill \\
  0\quad \text{otherwise} \hfill \\ 
\end{gathered}  \right. ,\end{equation}
so it follows that the norm of $L_{k}$ is given by
\begin{equation}
{\left\| {{L_k}} \right\|^2}  = \frac{2}{{2k + 1}}.
\end{equation}

The Legendre polynomials form a basis for the space $\LTwo{-1,1}$, which means that any
square integrable function $u$ on $[-1,1]$ can be represented as an \textit{infinite} series in Legendre polynomials
\begin{equation}
	u = \sum\limits_{k = 0}^\infty  {{{\hat u}_k}{L_k}(x)}\quad {\text{for all }}\;u\in\LTwo{-1,1}. 
	\label{eq:InfiniteLegendreSeries_DAK}
\end{equation} 
The spectral coefficients $\hat u_{k}$ are found as usual through orthogonal projection
\begin{equation}\iprod{u,{L_n}} = \sum\limits_{k = 0}^\infty  {{{\hat u}_k}\iprod{{L_k},{L_n}} }  = \sum\limits_{k = 0}^\infty  {{{\hat u}_k}{{\left\| {{L_k}} \right\|}^2}{\delta _{kn}}}  = {{\hat u}_n}{\left\| {{L_n}} \right\|^2},\end{equation}
so that
\begin{equation}
{{\hat u}_n} = \frac{{\iprod{u,{L_n}}}}{{{{\left\| {{L_n}} \right\|}^2}}},\qquad n = 0,1, \ldots ,\infty.
 \end{equation}

The best approximation of $u$ by the polynomial $U\in\PN{N}$ defined in \eqref{eq:LegendrePolynomialInterpolant_DAK} is to choose $\hat C = \hat u$ for $k=0,1,2,\ldots, N$, because then the error is orthogonal to the approximation space,
\begin{equation}
\inorm{u-U} = \inorm{\sum\limits_{k = N+1}^\infty  {{{\hat u}_k}{L_k}(x)}}.
\label{eq:TruncationError}
\end{equation}
Series truncation is the orthogonal projection of $\LTwo{-1,1}$ onto $\PN{N}(-1,1)$ with respect to the continuous inner product, and we call that
approximation $U = \mathbb{P}^{N}\!(u)$, where $\mathbb{P}^{N}$ is called the truncation operator.

\CCLsubsection{Legendre Polynomial Interpolation}

Alternatively, a function $u$ can be approximated by a polynomial interpolant of degree $N$ that passes through $N+1$ nodal points. Polynomial spectral element methods approximate a function $u(x)$ as a high-order Legendre polynomial interpolant $U(x)$ represented in 
\eqref{eq:LegendrePolynomialInterpolant_DAK} where the spectral interpolation coefficients $\hat C_{k}$ are determined so that 
\begin{equation}
\IN{u}( x_{j}) = U(x_{j}) = u(x_{j})\quad j = 0, 1, \ldots,N.
\label{eq:InterpolationConditions_DAK}
\end{equation}
We will use the notation $\IN{u}(x)$ or just $\IN{u}$ without the argument to denote the polynomial interpolant of order $N$, with $\mathbb{I}^{N}$ being the interpolation operator. Upper case functions like $U(x)$ will denote polynomial interpolants, whereas lower case functions can be anything, unless otherwise noted or convention dictates otherwise.

The preferable approach to find the coefficients $\hat C_{k}$ to satisfy \eqref{eq:InterpolationConditions_DAK} is to mimic the orthogonal projection process used to find the $\hat u_{k}$. Suppose that we have a
\textit{discrete inner product} $\iprodN{\cdot,\cdot}$ with the property
\begin{equation}
\iprodN{L_{k},L_{n}}=\inormN{L_{n}}^{2}\delta_{kn},
\end{equation}
where $\inormN{L_{n}}^{2}=\iprodN{L_{n},L_{n}}$. Then the interpolation coefficients $\hat C_{n}$, for $n=0,1,2,\ldots,N$ could be computed without solving a full Vandermonde matrix system as
\begin{equation}
{{\hat C}_n} = \frac{{{{\iprodN{u,{L_n}}}}}}{\inormN{L_{n}}^{2}}\quad n = 0,1, \ldots ,N,
\label{eq:ModalPolynomialRep_DAK}
\end{equation}
since
\begin{equation}
\iprodN{u,L_{n}} = \sum\limits_{k = 0}^\infty  {{{\hat C}_k}{\iprod{{L_k}(x),{L_n}}_N}}  = \sum\limits_{k = 0}^\infty  {{{\hat C}_k}\inormN{L_n}^2{\delta _{kn}}}  = {{\hat C}_n}{\inorm{L_n}^2}.
\end{equation}

An alternative -- and equivalent -- representation of the interpolant $U$ is to use the \textit{Lagrange} or \textit{nodal} form
\begin{equation}
U(x) = \sum\limits_{j = 0}^N {{U_j}{\ell _j}(x)},
\label{eq:LagrangeFormInterpolant_DAK}
 \end{equation}
 where $U_{j}= U\left(x_{j}\right)$ and the $\ell_{j}(x)$ are the Lagrange interpolating polynomials with nodes at the same points $x_j
$\begin{equation}
{\ell _j}(x) = \prod\limits_{i = 0;i \ne j}^N {\frac{{\left( {x - {x_i}} \right)}}{{\left( {{x_j} - {x_i}} \right)}}}  \in {\PN{N}},
\end{equation}
 that clearly possess the Kronecker delta property
 \begin{equation}
 \ell_{j}(x_{n}) = \delta_{jn}.
 \label{eq:kron_delta_Lagrange}
 \end{equation}

It remains, then to find an appropriate quadrature for the discrete inner product and interpolation nodes.

\CCLsubsection{Legendre Gauss Quadrature and the Discrete Inner Product}\label{sec:LGLQuadrature}

We can construct a discrete inner product with the desired properties by approximating the true inner product with a form of Gauss-Legendre quadrature known as the \textit{Gauss-Lobatto-Legendre} (or just Gauss-Lobatto) quadrature rule. The Gauss-Lobatto quadrature approximation of a function $f$ is a weighted sum of nodal values,
\begin{equation}
\ivolN {f(x)\dx}  \equiv \sum\limits_{j = 0}^N {f( {{x_j}} ){w_j}} ,	\label{eq:quadratureFormula_DAK}
\end{equation}
where the nodes $x_{j}$ are
\begin{equation}
	{x_j} =  + 1, - 1,{\text{ and the zeros of }}{L_N'}\left( x \right),
	\label{eq:GaussLobattoNodesDef_DAK}
\end{equation}
and the quadrature weights are 
\begin{equation}
{w_j} = \frac{2}{{N\left( {N + 1} \right)}}\frac{1}{\left[{L_N( {{x_j}} )}\right]^{2}}.
\label{eq:LobattoQuadratureWeights_DAK}
\end{equation}
The Gauss-Lobatto rule is exact when the integrand is a polynomial of degree $2N-1$ so we can write
\begin{equation}
\ivolN {f(x)\dx} = \irefInt {f(x)\dx} \quad {\text{for all }}f \in {\PN{2N - 1}}.
\label{eq:GaussLobQuadExactness_DAK}
\end{equation}

The quadrature rule \eqref{eq:quadratureFormula_DAK} allows us to define the discrete inner product as
\begin{equation}
\iprodN{u,v} \equiv \sum\limits_{j = 0}^N {u( {{x_j}}) v({{x_j}}){w_j}},
\label{eq:DiscreteInnerProductQuadrature_DAK}
 \end{equation} 
which has the desired orthogonality properties. If we replace $u$ and $v$ by the Legendre polynomials $L_{k}$ and $L_{n}$, then provided $k+n\leqslant 2N-1$,
\begin{equation}
\iprodN{{L_k},{L_n}} = \irefInt {{L_k}{L_n}\dx}  = \inormN{{L_k}}^{2}{\delta _{kn}}.
\end{equation}
The quadrature is therefore exact for all the inner products except $L_{N}$ with itself. That means for all $k<N$, $\inormN{L_{k}}=\inorm{L_{k}}$. The last case needs to be computed directly and separately leading to 
\begin{equation}
\inormN{L_{k}}^{2} = \left\{ \begin{gathered}
  \frac{2}{{2k + 1}}\quad k < N \hfill \\
  \frac{2}{N}\quad k = N \hfill \\ 
\end{gathered}  \right. .
\end{equation}

We can expose several useful properties from the definition of the discrete inner product. For any function $g(x)$, \eqref{eq:DiscreteInnerProductQuadrature_DAK} and the interpolation conditions \eqref{eq:InterpolationConditions_DAK} imply
that
\begin{equation}
\iprodN{\IN{g},V} = \sum\limits_{j = 0}^N {{g_j}{V_j}{w_j}} =\iprodN{g,V}\quad {\text{for all }}V\in \PN{N}.
\label{eq:DiscreteIPProjection_DAK}
\end{equation}
So any time the argument on the left is seen, it can be viewed as the interpolant of that argument. For example, for $U,V,W\in\PN{N}$,
\begin{equation}
	\iprodN{UV,W} =\iprodN{\IN{UV},W}.
	\label{eq:TripleInnerProduct_DAK}
\end{equation}

Also, the exactness of the quadrature \eqref{eq:GaussLobQuadExactness_DAK} implies that if the product $UV\in\PN{2N-1}$,
\begin{equation}
\iprodN{U,V} = \iprod{U,V}.
\end{equation}
For instance, the discrete inner product of a polynomial approximation with its first derivative is exact
\begin{equation}
\iprodN{U,U_{x}} = \iprod{U,U_{x}}.
\end{equation}

However, it also follows that for $U\in\mathbb{P}^{N}$
\begin{equation}
	\inormN{U}\ne\inorm{U}.
\end{equation}
Nevertheless, the discrete norm is \textit{equivalent} to the continuous norm in that for $U\in\mathbb{P}^{N}$ the discrete norm is bounded from above
and below by the continuous norm \citep{ISI:A1982NE30900005}
\begin{equation}
\inorm{U}\leqslant\inormN{U}\le\sqrt{2+\frac{1}{N}}\inorm{U}\leqslant\sqrt{3}\inorm{U}.
\label{eq:NormEquivalence_DAK}
\end{equation}
The discrete norm is never smaller than the true $\mathbb{L}^{2}$ norm, and is never more than about 73\% larger than the true norm. Equivalence says that if
an approximation converges in the discrete norm, it converges in the continuous norm, too.

\CCLsubsection{Aliasing Error}

Since the Gauss-Lobatto quadrature \eqref{eq:quadratureFormula_DAK} is exact only for polynomials of degree $2N-1$, the question is raised as to what errors are created by using the discrete inner product \eqref{eq:DiscreteInnerProductQuadrature_DAK} compared to the exact integral \eqref{eq:L2InnerProductOfL_DAK}. This is of particular interest when the function being interpolated is not a polynomial or not a polynomial of sufficiently low-order. The answer comes by comparing the exact coefficients $\hat u_{k}$ and the interpolation coefficients $\hat C_{k}$. 

To relate the exact and interpolation coefficients, we substitute 
\eqref{eq:InfiniteLegendreSeries_DAK} into \eqref{eq:ModalPolynomialRep_DAK} and use the discrete orthogonality to see that
\begin{equation}
\resizebox{\textwidth}{!}{$\displaystyle{
{{\hat C}_k} = \frac{{{\iprodN{\sum\limits_{n = 0}^\infty  {{{\hat u}_n}{L_n}} ,{L_k}}}}}{\inormN{L_k}^{2}} = {{\hat u}_k} + \frac{{\sum\limits_{n = N + 1}^\infty  {{{\hat u}_n}{\iprodN{{L_n},{L_k}}}} }}{\inorm{L_k}_{N}^{2}}\equiv \hat u_{k} + \hat A_{k}\quad k = 0,1, \ldots ,N.
}$}
\label{eq:AliasingExpression_DAK}
\end{equation}
So \eqref{eq:AliasingExpression_DAK} shows that the interpolation coefficients computed with quadrature are the exact coefficients plus an \textit{aliasing} error that depends on the discrete projection of $L_{n}$ onto $L_{k}$ for $n > N$. The interpolant can therefore be written as
\begin{equation}
\resizebox{\textwidth}{!}{$\displaystyle{
\begin{aligned}
\IN{u} &= \sum\limits_{k=0}^N \hat C_k L_k(x) = \sum\limits_{k = 0}^N {{{\hat u}_k}L_k(x)}  
+ \sum\limits_{k = 0}^N {\left\{ {\frac{1}{{\inormN{L_k}^{2}}}\sum\limits_{n = N + 1}^\infty  {{{\hat u}_n}\iprodN{{L_n},{L_k}}} } \right\}{L_k}(x)}
\\& = u - \sum\limits_{k = N + 1}^\infty  {{{\hat u}_k}L_k( x )} + \sum\limits_{k = 0}^N {\left\{ {\frac{1}{{\inormN{L_k}^{2}}}\sum\limits_{n = N + 1}^\infty  {{{\hat u}_n}\iprodN{{L_n},{L_k}}} } \right\}{L_k}(x)}
\\& = u +\left\{- \sum\limits_{k = N + 1}^\infty  {{{\hat u}_k}L_k( x )}\right\} + \sum\limits_{k = 0}^N {\hat A_{k}{L_k}(x)}.
\label{eq:AliasingPolynomial_DAK}
 \end{aligned}
 }$}
\end{equation}

So the interpolant is the actual function plus two errors. The first error is the truncation error due to the finite number of modes available. The second error is the aliasing error, which arises because discretely the higher order modes have a non-zero contribution to the low-order modes, as represented in the fact that the discrete inner products do not vanish. 

The indices on the sums over $k$ in \eqref{eq:AliasingPolynomial_DAK} show that the truncation and aliasing errors are orthogonal to each other, since
\[
\iprod{\sum\limits_{k = N + 1}^\infty  {{{\hat u}_k}{L_k}} ,\sum\limits_{n = 0}^N {{A_n}{L_n}}} = \sum\limits_{k = N + 1}^\infty {\sum\limits_{n = 0}^N {{{\hat u}_k}{A_n}\iprod{{L_k},{L_n}}} }  = 0.
\]
 Note also that \eqref{eq:AliasingPolynomial_DAK} says that if $u\in\PN{N}$ for which $\hat u_{k}=0,\; k>N$, then $\IN{u}=u$, as expected.

The projection result \eqref{eq:TripleInnerProduct_DAK} shows that the discrete inner product of a compound argument (function of polynomials) with a polynomial introduces aliasing errors. Compound arguments appear when projecting a flux, like $f(U) = \frac{1}{2}U^{2}$ for the Burgers equation, onto the basis functions. 

For \textit{polynomial} compound functions, where the result is the product of polynomials, the aliasing error created by discrete inner products can be eliminated by \textit{consistent integration}, more commonly referred to as ``overintegration,'' at the cost of extra evaluations of the function. The idea is to evaluate the inner product at $M>N$ points so that the discrete inner product is exact. For the product $UV\in\PN{2N}$, for example, the problem is to find $M$ so that for $W\in\PN{N}$
\begin{equation}
\iprod{UV,W}_{M} = \irefInt {UVW\dx} .
\end{equation}
The product $UVW\in\mathbb{P}^{3N}$ and the Gauss-Lobatto quadrature is exact for arguments in $\mathbb{P}^{2M-1}$. Therefore, there is no aliasing error if
\begin{equation}
3N=2M-1\quad \text{or }\quad M=\frac{3N+1}{2}.
\end{equation}
In other words, aliasing due to the discrete inner product can be eliminated by evaluating the inner product at 
\begin{equation}
M > \frac{3}{2}N
\end{equation}
points, 50\% more than used in the interpolation.  More generally, for a polynomial function $F\in\mathbb{P}^{p}$, aliasing can be avoided when taking
\begin{equation}
M>\frac{p+1}{2}.
\end{equation}
It should be clear, however, that if $F$ is not a polynomial function, then \eqref{eq:TripleInnerProduct_DAK} implies that aliasing will be present except in the limit as $M\rightarrow\infty$, i.e. taking an infinite number of quadrature points and converging the discrete inner product to the continuous one.

\CCLsubsection{Spectral Differentiation}

Derivatives of functions are approximated by the derivatives of their polynomial interpolants
\begin{equation}
u' \approx {\left( \IN{u} \right)^\prime }.
\end{equation}
The interpolant can be represented in either the nodal or modal form. The choice can be made solely on efficiency considerations. 

In Legendre spectral element methods, derivative approximations are computed by matrix-vector multiplication where the vector holds the nodal values of the approximation and the matrix holds derivatives of the Lagrange interpolating polynomials. The derivative of the Lagrange form interpolant is
\begin{equation}
U' = \sum\limits_{n = 0}^N {{U_n}{\ell_n'}( x)}. 
\end{equation}
When evaluated at the Gauss-Lobatto points $x_{j},\; j = 0, 1, \ldots,N$,
\begin{equation}
{{U}'_j} = \sum\limits_{n = 0}^N {{U_n}{{\ell}'_n}( {{x_j}} )}  = \sum\limits_{n = 0}^N {{U_n}{\dmat_{jn}}}
\label{eq:DerivativeOperation_DAK}
 \end{equation}
 where $\dmat_{jn}= {{\ell}'_n}( {{x_j}} )$  are the elements of the \textit{derivative matrix}, $\dmat$. Thus, derivatives are computed by matrix-vector multiplication
 \begin{equation}
 \matrixvec{U}'=\dmat\matrixvec{U},
 \end{equation}
 where $\matrixvec{U} = [U_{0}\; U_{1}\;\ldots\; U_{N}]^{T}$. 
 
 One often noted feature of the approximation \eqref{eq:DerivativeOperation_DAK} is that the approximation to the derivative includes only points in the domain, even up to the boundary points, independent of the approximation order. This is a feature not held, for instance, with high-order finite difference methods, where external ``ghost points'' appear if the stencil is used near the boundaries.

Derivatives of functions of a polynomial, like a flux $f(U)$, a product $UV$, or other compound quantity are computed nodally. For instance, the product $Q=UV\in\mathbb{P}^{2N}$ is computed by $\dmat\matrixvec{Q}$ where $Q_{j} = U_{j}V_{j}$. Differentiation of the product is therefore equivalent to
\begin{equation}
{\left( {UV} \right)^\prime } \approx {\left( \IN{UV} \right)^\prime }.
\end{equation}
As a result, there is an aliasing error associated with representing the product as a polynomial of degree $N$. 

One unfortunate aspect of polynomial differentiation is that \textit{differentiation and interpolation do not commute}, i.e.,
\begin{equation}
{\left( \IN{u} \right)^\prime } \ne {\mathbb{I}^{N-1}}\!\left( {u'} \right).
\end{equation}
 As a consequence of interpolation and differentiation not commuting, common rules like the product and chain rules do not hold except in special cases. For example, the product rule 
 \begin{equation}
 {\left( \IN{UV} \right)^\prime } \ne {\mathbb{I}^{N - 1}}\!\left( {U'V} \right) + {\mathbb{I}^{N - 1}}\!\left( {UV'} \right)
 \label{eq:ProductRuleError}
 \end{equation}
 unless $UV\in\mathbb{P}^{N}$. 
\CCLsubsection{Spectral Accuracy}\label{sec:SpectralAccuracy}
The truncation, $\mathbb{P}^{N}(u)$, interpolation, $\IN{u}$, derivative, $\left(\IN{u} \right)'$, and Gauss-Lobatto quadrature, $\int_N {u\dx}$, approximations all possess what is known as \emph{spectral accuracy}: The rate of convergence depends only on the smoothness of $u$. For very smooth functions they converge extremely fast, and with enough smoothness, they converge exponentially fast. In this section, we review these facts and their meaning, and leave the technical derivations and more precise forms to  \cite{canuto2007}.

Spectral accuracy follows from the fact that the truncation error (c.f. \eqref{eq:TruncationError})
\begin{equation}
\inorm{u-\mathbb{P}^{N}(u)} = \inorm{\sum\limits_{k = N+1}^\infty  {{{\hat u}_k}{L_k}(x)}} = \sum\limits_{k = N+1}^\infty  {\frac{2}{2k+1}\left|{{\hat u}_k}\right|^{2}}
\end{equation}
depends on the size of the modal coefficients, $\hat u_{k}$. From Fourier analysis, we know that the smoother $u$ is, the faster those coefficients decay. With $\hat u_{N+1}$ being the dominant mode, the truncation error decays $\sim \left|\hat u_{N+1}\right|$. So we get the relationship between smoothness and truncation error through the modal coefficients. 

To write the error convergence more precisely, we define the Sobolev norm,
\begin{equation}
\inorm{u}_{H^{m}}^{2} = \sum_{n=0}^{m}\inorm{\frac{d^{n}u}{dx^{n}}}^{2},
\end{equation}
so the smoother the function $u$ is, the larger the index $m$ for which the norm is finite, i.e. $\inorm{u}_{H^{m}}<\infty$. Note that when $m=0$, the
Sobolev norm reduces to the $\LTwo{-1,1}$ energy norm. Also, if $\inorm{u}_{H^{m}}^{2}<\infty$ for some $m$, then the energy norm of $u$ and each derivative
individually up to order $m$ is also bounded.

In terms of the Sobolev norm, the truncation error satisfies
\begin{equation}
\inorm{u - \mathbb{P}^{N}(u)}\leqslant CN^{-m}\inorm{u}_{H^{m}},
\label{eq:TruncationErrorSobolevNorm}
\end{equation}
where $C$ is a generic constant.
Eq. \eqref{eq:TruncationErrorSobolevNorm} is the statement of spectral accuracy. We see directly that for a fixed smoothness implied by $m$, the error converges like $N^{-m}$, which is rapid for large $m$. For fixed $N$, the convergence rate increases as $m$ increases. If all derivatives exist, so that we can take  $m\rightarrow\infty$, the approximation is said to have \emph{infinite order convergence}. 

The interpolation error seen in \eqref{eq:AliasingPolynomial_DAK} is the sum of the truncation error, which decays spectrally fast, plus the aliasing error \eqref{eq:AliasingExpression_DAK}, which
also depends on the rate of decay of the modal coefficients. As a result, the interpolation error is also spectral, though larger than the truncation error and requires slightly more smoothness \citep{canuto2007}. Ultimately, the interpolation error is also bounded like in \eqref{eq:TruncationErrorSobolevNorm} but only if $m>1/2$.

Since the interpolation error converges spectrally fast, it is not surprising that the error in the derivative of the interpolant is spectrally accurate too, though at a lower rate.
For the derivative,
\begin{equation}
\inorm{u' - \left(\IN{u}\right)'}\leqslant CN^{1-m}\inorm{u}_{H^{m}}.
\label{eq:DerivativeErrorSobolevNorm}
\end{equation}
Note that since the differentiation error is spectrally accurate, it follows that the product rule error in \eqref{eq:ProductRuleError}
also converges spectrally fast.

Finally, the discrete inner product, and by extension the quadrature is spectrally accurate. For a function $u$ and polynomial $\phi\in\PN{N}$,
\begin{equation}
\left| \iprod{u,\phi} -\iprodN{u,\phi}\right|\leqslant CN^{-m}\inorm{u}_{H^{m}}\inorm{\phi}.
\end{equation} 

Spectral convergence becomes exponential convergence if $u$ is so smooth that it is analytic (in the complex variables sense) in some ellipse in the complex plane around the foci $-1,1$. Exponential convergence is sometimes confused with spectral accuracy, whereas it is instead a special case.
Recently \citet{XIE:2013hb} have shown that
\begin{equation}
\mathop {\max }\limits_{|x| \leqslant 1} \left| {u - \IN{u}} \right| \leqslant C\left( \rho  \right){N^{3/2}}{e^{ - N\ln (\rho )}},
\end{equation} 
and
\begin{equation}
\mathop {\max }\limits_{0 \leqslant j \leqslant N} \left| {u' - {{\left( {\IN{u}} \right)}^\prime }} \right| \leqslant C\left( \rho  \right){N^{7/2}}{e^{ - N\ln (\rho )}},
\end{equation} 
where $\rho$ increases with the size of the ellipse, and hence the region of analyticity. For large enough $N$, the exponential decay dominates the polynomial growth. Gauss quadrature is also exponentially convergent for analytic functions by virtue of the interpolation convergence.

\CCLsubsection{The Discrete Inner Product and Summation-by-Parts}\label{sec:disc_IP_and_SBP}

One of the most important properties of the discrete inner product \eqref{eq:DiscreteInnerProductQuadrature_DAK}  with regards to the methods we develop here is the summation-by-parts (SBP) property. The summation-by-parts property is the discrete equivalent of the integration-by-parts property 
\begin{equation}
\iprod{u,{v'}} = \irefInt {u{v'}\dx}  = \left. {uv} \right|_{ - 1}^1 - \iprod{{u'},v},
\label{eq:IntegrationByParts1D_DAK}
\end{equation}
or equivalently
\begin{equation}
\iprod{u,{v'}} + \iprod{{u'},v} = \left. {uv} \right|_{ - 1}^1.
\label{eq:IntegrationByParts1D2_DAK}
\end{equation}
If $U, V\in\PN{N}$ then $UV'\in\PN{2N-1}$ and $U'V\in\PN{2N-1}$. Using the exactness between the integral and the Gauss-Lobatto quadrature,
\begin{equation}
\iprodN{U,V'}=\iprod{U,V'}=\left. {UV} \right|_{ - 1}^1-\iprod{U',V}= \left. {UV} \right|_{ - 1}^1-\iprodN{U',V}.
\end{equation}
Therefore, the \textit{summation-by-parts} formula is
\begin{equation}
	\iprodN{U,V'}= \left. {UV} \right|_{ - 1}^1-\iprodN{U',V},
	\label{eq:SBP1D1_DAK}
\end{equation}
or equivalently,
\begin{equation}
	\iprodN{U,V'}+\iprodN{U',V}= \left. {UV} \right|_{ - 1}^1,
	\label{eq:SBP1D2_DAK}
\end{equation}
which is the discrete equivalent of the integration-by-parts property \eqref{eq:IntegrationByParts1D2_DAK} held by the continuous integral.

It is interesting to note that if $U\in\mathbb{P}^{N}$ and $V\in\mathbb{P}^{N}$, the summation-by-parts
 formula \eqref{eq:SBP1D2_DAK} gives
\begin{equation}
 \left. {UV} \right|_{ - 1}^1 = \ivolN {{{\left( \IN{UV} \right)}^\prime }\dx}  = \ivolN {U'V\dx}  + \ivolN {UV'\dx} ,
 \end{equation}
 which says that the mean value of the error due to the lack of commutativity of interpolation and differentiation is zero.
 
\CCLsubsubsection{Integral Quantities in Matrix-Vector Form}

Since the nodal degrees of freedom can be represented as a vector, integral quantities like the inner product and summation-by-parts can be written in matrix-vector form. For instance, let us define the diagonal mass matrix, whose entries are the Gauss-Lobatto quadrature weights,
\begin{equation}\label{eq:mass_matrix}
\mmat = 
\left[\begin{array}{ccc}
w_0 &   & \text{\huge0} \\ 
  & \ddots  &   \\
\text{\huge0} &   & w_N
\end{array}\right],
\end{equation}
with which we can write the discrete inner product as
\begin{equation}
\iprodN{U,V} = \sum_{j=0}^{N}U_{j}w_{j}V_{j}=\matrixvec U^T \mmat \matrixvec  V.
\end{equation}
Similarly, the quadrature of $F$ can be expressed as
\begin{equation}
\ivolN {F(x)\dx} = \iprodN{1,F}= 1^{T}\mmat \matrixvec F. 
\end{equation}

Written in matrix form, we can show that that the summation-by-parts formula \eqref{eq:SBP1D2_DAK} is solely a property of the derivative and mass matrices and a boundary matrix 
\begin{equation}
\bmat = \left[\begin{array}{ccccc}-1 &  &  &  & \text{\huge0} \\ & 0 &  &  &  \\ &  & \ddots &  &  \\ &  &  & 0 &  \\\text{\huge0} &  &  &  & 1\end{array}\right].
\end{equation}
If we write \eqref{eq:SBP1D2_DAK} in summation form, 
\begin{equation}
\sum_{i=0}^{N}U_{i}w_{i}\left(\dmat \matrixvec V\right)_{i}+ \sum_{i=0}^{N}w_{i}\left(\dmat \matrixvec U\right)_{i}V_{i}=\left\{U_{N}V_{N}-U_{0}V_{0}\right\},
\end{equation}
we see the equivalent matrix-vector equation
\begin{equation}
\matrixvec U^{T}\mmat\dmat\matrixvec V + \left(\mmat\dmat\matrixvec U\right)^{T}\matrixvec V=\matrixvec U^{T}\bmat\matrixvec V,
\end{equation}
which can be factored as
\begin{equation}
\matrixvec U^{T}\left\{\mmat\dmat + \left(\mmat\dmat\right)^{T}-\bmat\right\}\matrixvec V=0.
\end{equation}
Since the polynomials from which the nodal values in the vectors $\matrixvec U,\matrixvec V$ are arbitrary, the components are linearly independent and it follows that
\begin{equation}
\mmat\dmat + \left(\mmat\dmat\right)^{T}=\bmat.
\label{eq:SBP_with_D}
\end{equation}
Commonly, the matrix $\qmat = \mmat\dmat$ is defined leaving
\begin{equation}
\qmat + \qmat^{T}=\bmat.
\label{eq:SBPWithQ_{DAK}}
\end{equation}

The relation \eqref{eq:SBPWithQ_{DAK}} is a matrix expression of the summation-by-parts property, \eqref{eq:SBP1D2_DAK}. Summation-by-parts with matrix operators was introduced in the finite difference community, e.g. \citep{kreiss1,kreiss1977,kreiss1974finite,strand1994} and, e.g., the review article \citet{svard2014}, with the goal to mimic finite element type energy estimates with local stencil based differentiation operators, i.e. finite differences. The matrix expression \eqref{eq:SBPWithQ_{DAK}} shows that collocation type spectral elements with Gauss-Lobatto nodes may structurally be interpreted as diagonal norm summation-by-parts finite difference methods.  

From \eqref{eq:SBPWithQ_{DAK}} it is possible to assess the structure of $\qmat$ and determine many of its entries. With the aid of the Lagrange nodal polynomial basis functions and the collocated Gauss-Lobatto quadrature, it is easy to see that the entries of the $\qmat$ matrix can be directly computed as
\begin{equation}
\qmat_{ij} = \iprodN{\ell_j',\ell_i},\quad i,j=0,...,N.
\label{eq:sbp_property0_GG}
\end{equation} 
Due to the consistency of polynomial interpolation, it follows that it is possible to exactly represent a constant function. From this consistency it follows that the derivative of a constant can be computed exactly, which translates into the matrix condition 
\begin{equation}
\dmat\, 1 = 0 \;\Rightarrow\; \qmat \, 1 = 0 \;\Rightarrow\; \sum\limits_{j=0}^N \qmat_{ij} = 0,\;\; i=0,...,N.
\label{eq:sbp_property1_GG}
\end{equation} 
That is, the sum of rows of the matrix $\qmat$ (or $\dmat$) are equal to zero. Directly multiplying \eqref{eq:SBPWithQ_{DAK}} with a vector containing only ones, $1$, and using \eqref{eq:sbp_property1_GG}, we get for the sum of the columns
\begin{equation}
\sum\limits_{i=0}^N \qmat_{ij} = 
\begin{cases}
-&1,\quad i=0,\\
&0,\quad i=1,...,N-1,\\
+&1,\quad i=N.
\end{cases}
\label{eq:sbp_property2_GG}
\end{equation} 
If we assess the diagonal entries of \eqref{eq:SBPWithQ_{DAK}}, we immediately get 
\begin{equation}
\qmat_{00} = -\frac{1}{2},\quad \qmat_{NN} = \frac{1}{2},\quad \qmat_{ii} = 0,\,i=1,...,N-1.
\label{eq:sbp_property3_GG}
\end{equation} 
Lastly, the $\qmat$ matrix is almost skew-symmetric, i.e., 
\begin{equation}
\qmat_{ij} = -\qmat_{ji},\quad\forall\, i,j\text{ (except for $\qmat_{00}$ and $\qmat_{NN}$)}.
\label{eq:sbp_property4_GG}
\end{equation}

\CCLsubsection{Extension to Multiple Space Dimensions}\label{sec:InterpMultipleSpace_DAK}

In multiple space dimensions, functions are approximated by tensor products of the one dimensional interpolants on the quadrilateral
 $E_{2} = [-1,1]\times[-1,1]$ or hexahedron  $E_{3} = [-1,1]^3$. We let 
 \begin{equation}
 \spacevec{x}=(x,y,z)=\left(x_{1},x_{2},x_{3}\right) =x\hat x + y\hat y+z\hat z=\sum\limits_{i=0}^{3}x_{i}\hat x_{i},
 \label{eq:spaceVariables_DAK}
 \end{equation}
where $\hat x,\hat y,\hat z$ are the unit vectors in the three coordinate directions
, with the similar definition in two.
For functions $u(x,y)$ and $u(x,y,z)$ defined in $E_{2}$ and $E_{3}$, the Lagrange forms of the interpolant are
\begin{equation}
\IN{u}(x,y) = \sum\limits_{i,j = 0}^N {{u_{ij}}{\ell _i}( x ){\ell _j}( y )},
\label{eq:Interpolant2D}
 \end{equation}
 and
 \begin{equation}
 \IN{u}(x,y,z) = \sum\limits_{i,j,k = 0}^N {{u_{ijk}}{\ell _i}( x ){\ell _j}( y ){\ell _k}( z )},
 \label{eq:Interpolant3D}
\end{equation}
where $u_{ijk}=u\left(x_{i},y_{j},z_{k}\right)$, etc. The $x_{i}$, $y_{j}$ and $z_{k}$ are located at the Gauss-Lobatto nodes.  
For notational simplicity, we have assumed the same polynomial order in each space direction, though this is not necessary in practice. With this assumption, we will denote the space of polynomials of degree $N$ in each space direction also as $\mathbb{P}^{N}$.

\CCLsubsubsection{Tensor Product Spectral Differentiation}

 The tensor product makes the computation of partial derivatives simple and efficient. For example
 \begin{equation}\label{eq:disc_der_3D}
 {\left. {\frac{{\partial U}}{{\partial x}}} \right|_{ijk}} = \sum\limits_{n,m,l = 0}^N {{U_{nml}}{\ell'_n}( {{x_i}} ){\ell _m}( {{y_j}} ){\ell _l}( {{y_k}} )}  = \sum\limits_{n = 0}^N {{U_{njk}}{\ell'_n}( {{x_i}} )}  = \sum\limits_{n = 0}^N {{U_{njk}}{\dmat_{in}}},
  \end{equation}
  where the $j$ and $k$ sums drop out due to the Kronecker delta property of the Lagrange basis \eqref{eq:kron_delta_Lagrange}.
  
Let us assume that the nodal values are stored in an array format, and let us represent an array slice that defines a vector by a colon, ``:". Then the $x-derivative$ can be computed plane--by--plane
  \begin{equation}
  \matrixvec{U}_{jk}=\dmat{U}_{:jk},\quad j,k=0,1,\ldots,N.
  \end{equation}
Similar relations hold for the other partial derivatives, which allows us to write the spectral gradient as
\begin{equation}
\begin{split}
\vecNablaX {U_{ijk}} &= \sum\limits_{n = 0}^N {{U_{njk}}{\dmat_{in}}}\hat x  + \sum\limits_{n = 0}^N {{U_{ink}}{\dmat_{jn}}}\hat y  + \sum\limits_{n = 0}^N {{U_{ijn}}{\dmat_{kn}}}\hat z \\ &
= \dmat U_{:jk}\hat x + \dmat U_{i:j}\hat y + \dmat U_{ij:}\hat z\,,
\end{split}
\label{eq:spectralGradient_DAK}
\end{equation}
and the divergence as
\begin{equation}
\begin{split}
\vecNablaX \cdot {\spacevec F_{ijk}} &= \sum\limits_{n = 0}^N {{(F_{1})_{njk}}{\dmat_{in}}}  + \sum\limits_{n = 0}^N {{(F_{2})_{ink}}{\dmat_{jn}}}  + \sum\limits_{n = 0}^N {{(F_{3})_{ijn}}{\dmat_{kn}}}\\ &
= \dmat \left( F_1\right)_{:jk} + \dmat \left( F_2\right)_{i:k} + \dmat \left( F_3\right)_{ij:}\,.
\end{split}
\label{eq:spectralDivergence_DAK}
\end{equation}


\CCLsubsubsection{Discrete Inner Products and summation-by-parts}\label{sec:DIPAndSBP_DAK}

The discrete inner product is defined using the Gauss-Lobatto rule in each space direction. In 3D,
\begin{equation}\label{eq:disc_IP_3D}
\iprodN{U,V} = \sum\limits_{i,j,k = 0}^N {{U_{ijk}}{V_{ijk}}{w_i}{w_j}{w_k}}  \equiv \sum\limits_{i,j,k = 0}^N {{U_{ijk}}{V_{ijk}}{w_{ijk}}} .
\end{equation}
As in one space dimension, the discrete inner product is exact when the degree of $UV$ is $2N-1$ or less in each direction, i.e.
\begin{equation}
\iprodN{U,V} = \iprod{U,V}\quad \text{for all }UV\in\mathbb{P}^{2N-1}.
\end{equation}

 Because of the tensor product, it turns out that summation-by-parts still holds in multiple space dimensions. We show how in three space dimensions and 
for the partial derivatives in $x$. Let $U,V\in\mathbb{P}^{N}$. Then
\begin{equation}
\begin{gathered}
  U = \sum\limits_{n,m,l = 0}^N {{U_{nml}}{\ell _n}( x ){\ell _m}( y ){\ell _l}( z )}  \hfill \\
  {V_x} = \sum\limits_{\mu ,\nu ,\lambda  = 0}^N {{V_{\mu \nu \lambda }}{\ell'_\mu }( x ){\ell _\nu }( y ){\ell _\lambda }( z )}  \hfill \\ 
\end{gathered}
 \end{equation}
 and so
\begin{equation}
\iprodN{U,{V_x}} = \sum\limits_{n,m,l = 0}^N {\sum\limits_{\mu ,\nu ,\lambda  = 0}^N {{U_{nml}}{V_{\mu \nu \lambda }}\iprodN{{\ell _n}{\ell _m}{\ell _l},{\ell'_\mu }{\ell _\nu }{\ell _\lambda }}} } .
\label{eq:LHSofSBP_DAK}
\end{equation} 
 The discrete inner products in the sums factorize,
\begin{equation}
\iprodN{{\ell _n}{\ell _m}{\ell _l},{\ell'_\mu }{\ell _\nu }{\ell _\lambda }} = \left( {\ivolN {{\ell _n}{\ell'_\mu }\dx} } \right)\left( {\ivolN {{\ell _m}{\ell _\nu }\,\text{d}y} } \right)\left( {\ivolN {{\ell _l}{\ell _\lambda }\,\text{d}z} } \right).
\end{equation}
 We then use summation-by-parts on the first factor
 \begin{equation}
 \resizebox{\textwidth}{!}{$\displaystyle{
\iprodN{{\ell _n}{\ell _m}{\ell _l},{\ell'_\mu }{\ell _\nu }{\ell _\lambda }} = \left( {\left. {{\ell _n}{\ell _\mu }} \right|_{x =  - 1}^1 - \ivolN {{\ell'_n}{\ell _\mu }\dx} } \right)\left( {\ivolN {{\ell _m}{\ell _\nu }\,\text{d}y} } \right)\left( {\ivolN {{\ell _l}{\ell _\lambda }\,\text{d}z} } \right)
  }$}
\end{equation}
 and recombine the discrete inner product
 \begin{equation}
\iprodN{{\ell _n}{\ell _m}{\ell _l},{\ell'_\mu }{\ell _\nu }{\ell _\lambda }} = \ivolN {\left. {{\ell _n}{\ell _\mu }} \right|_{x =  - 1}^1{\ell _m}{\ell _\nu }{\ell _l}{\ell _\lambda }\,\text{d}y\text{d}z}  -\iprodN{{\ell'_n}{\ell _m}{\ell _l},{\ell _\mu }{\ell _\nu }{\ell _\lambda }}.
\label{eq:rewrite_3d_sbp}
 \end{equation}
 Substituting the discrete inner product \eqref{eq:rewrite_3d_sbp} into \eqref{eq:LHSofSBP_DAK} gives the summation-by-parts formula
 \begin{equation}
 \iprodN{U,V_{x}} = \ivolN {\left. {UV} \right|_{x =  - 1}^1\,\text{d}y\text{d}z}  - \iprodN{U_{x},V}.
 \label{eq:SBPProperty-x_{DAK}}
 \end{equation}
 The surface quadrature is precisely
 \begin{equation}
 \resizebox{\textwidth}{!}{$\displaystyle{
 \begin{aligned}
 \ivolN {\left. {UV} \right|_{x =  - 1}^1\,\text{d}y\text{d}z}  &\equiv \sum\limits_{j,k = 0}^N {U( {1,{y_j},{z_k}} )V( {1,{y_j},{z_k}} ){w_j}{w_k}}  - \sum\limits_{j,k = 0}^N {U( { - 1,{y_j},{z_k}} )V( { - 1,{y_j},{z_k}} ){w_j}{w_k}}  \\&= \sum\limits_{j,k = 0}^N {{U_{Njk}}{V_{Njk}}{w_j}{w_k}}  - \sum\limits_{j,k = 0}^N {{U_{0jk}}{V_{0jk}}{w_j}{w_k}}.
 \end{aligned}
 }$}
\end{equation}
 Equivalent results hold for the $y$ and $z$ derivatives and in two space dimensions.
\CCLsubsubsection{Multidimensional Summation-by-Parts and Divergence Theorem}
 The summation-by-parts property extends to three dimensions and to the divergence theorem. Let $\spacevec{F}\in \mathbb{P}^{N}$ be a vector
 \begin{equation}
 \spacevec{F} = \sum\limits_{i = 1}^d {{F_i}{{\hat x}_i}} ,
 \end{equation}
 where $d =2,3$ is the number of space dimensions. Then by adding the summation-by-parts property \eqref{eq:SBPProperty-x_{DAK}} for any $V\in\PN{N}$ and for each component of the vector $\spacevec F$ and its corresponding derivative, we get the multidimensional summation-by-parts theorem
 \begin{equation}
 \iprodN{\vecNablaX  \cdot \spacevec F,V} = \isurfEN {\spacevec F \cdot \hat nV\dS}  - \iprodN{\spacevec F,\vecNablaX V}.
 \label{eq:DiscreteGreens_DAK}
 \end{equation}
 If we then set $V=1$, then we get the discrete divergence theorem
  \begin{equation}
\volEN {\vecNablaX\cdot\spacevec F\,\text{d}\spacevec x}= \iprodN{\vecNablaX  \cdot \spacevec F,1} = \isurfEN {\spacevec F \cdot \hat n \dS}.
 \label{eq:DiscreteDivergence_DAK}
 \end{equation}
The divergence theorem is used to show conservation. In fact, we can say even more. With $V=1$ the quadrature is exact and so the discrete conservation is actually conservative in the integral sense.

\CCLsubsection{Summary}
We summarize the results of this section by showing the continuous and discrete equivalents for the spectral calculus in Table \ref{tab:CalcComparison_DAK}.
{\begin{table}[!ht]
		\centering
		\caption{Summary of calculus computations and rules on $E_{3}=[-1,1]^3$.}
		\label{tab:CalcComparison_DAK}
\begingroup
\setlength{\tabcolsep}{6pt} 
\renewcommand{\arraystretch}{2.6} 
\aboverulesep=0ex
\belowrulesep=0ex
\begin{adjustbox}{max width=\textwidth}
			\begin{tabular}{@{}l|r@{}}
				\toprule
				Continuous setting & Discrete setting\\
				\midrule
				\midrule
				$u,v,\spacevec f \in\LTwo{E}$ & $U,V,\spacevec F\in\PN{N}$ \\
				$\iprod{u,v}=\displaystyle\int\limits_E uv\,\text{d}\spacevec{x}$ & $ \iprodN{U,V} = \displaystyle\sum\limits_{i,j,k=0}^N U_{ijk}V_{ijk}w_iw_jw_k$\\
				$\inorm{u}^2=\iprod{u,u}$ & $\inormN{U}^2=\iprodN{U,U}$ \\
				$\vecNablaX u = u_x\hat x + u_y\hat y+u_z\hat z$ & 
				$\vecNablaX {U_{ijk}} = \left(\displaystyle\sum\limits_{n = 0}^N {{U_{njk}}{\dmat_{in}}}\right)\hat x  + \left(\displaystyle\sum\limits_{n = 0}^N {{U_{ink}}{\dmat_{jn}}}\right)\hat y  + \left(\displaystyle\sum\limits_{n = 0}^N {{U_{ijn}}{\dmat_{kn}}}\right)\hat z$ \\
				$\vecNablaX\cdot\vec f = \left( f_1\right)_x + \left( f_2\right)_y +  \left( f_3\right)_z $& $\vecNablaX \cdot {\spacevec F_{ijk}} = \displaystyle\sum\limits_{n = 0}^N {{(F_{1})_{njk}}{\dmat_{in}}}  + \displaystyle\sum\limits_{n = 0}^N {{(F_{2})_{ink}}{\dmat_{jn}}}  + \displaystyle\sum\limits_{n = 0}^N {{(F_{3})_{ijn}}{\dmat_{kn}}} $\\
				$ \iprod{\vecNablaX\cdot\spacevec f,v} = \mkern-4mu\displaystyle\int\limits_{\partial E}\mkern-3mu \spacevec f\cdot\hat nv\dS - \iprod{\spacevec f,\vecNablaX v}$ & $ \iprodN{\vecNablaX\cdot\spacevec F,V} = \mkern-9mu\displaystyle\int\limits_{\partial E,N}\mkern-9mu \spacevec F\cdot\hat n V\dS - \iprodN{\spacevec F,\vecNablaX V}$\\
				$  \displaystyle\int\limits_{E}\vecNablaX\cdot\spacevec f \,\text{d}\spacevec{x} =  \mkern-4mu\displaystyle\int\limits_{\partial E}\mkern-3mu\spacevec f\cdot\hat n \dS$ & $\displaystyle\int\limits_{\partial E,N}\mkern-5mu\vecNablaX\cdot\spacevec F \,\text{d}\spacevec{x} =\mkern-10mu\displaystyle\int\limits_{\partial E,N}\mkern-8mu\spacevec F\cdot\hat n \dS$ \\
				$ \left(uv\right)' = u'v + uv'$ & $ {\left( \IN{UV} \right)^\prime } \ne {\mathbb{I}^{N - 1}}\!\left( {U'V} \right) + {\mathbb{I}^{N - 1}}\!\left( {UV'} \right)$\\
				\bottomrule
		\end{tabular}
		\end{adjustbox}
		\endgroup
	\end{table}}

\CCLsection{The Compressible Navier-Stokes Equations}\label{sec:section_NSE_GG}

Compressible viscous flows are modelled by the Navier-Stokes equations,
\begin{equation}
\label{eq:nse}
{\statevec u_t} + \sum\limits_{i = 1}^3 {\frac{{\partial {\statevec f_{i}}}}{{\partial {x_i}}} = \overRe\sum\limits_{i = 1}^3 {\frac{{\partial \statevec f^v_{i}\left( {\statevec u,\vecNablaX\statevec u} \right)}}{{\partial {x_i}}}} }.
\end{equation}
The state vector contains the conservative variables of the density, $\rho$, the momenta, $\rho\spacevec v = (\rho v_1 \; \rho v_2 \; \rho v_3)^{T}$ and total energy $\rho E$ per unit volume,
\begin{equation}\statevec u = \left[ {\begin{array}{*{20}{c}}
  \rho  \\ 
  {\rho \spacevec v} \\ 
  {\rho E} 
\end{array}} \right] = \left[ {\begin{array}{*{20}{c}}
  \rho  \\ 
  {\rho v_1} \\ 
  {\rho v_2} \\ 
  {\rho v_3} \\ 
  {\rho E} 
\end{array}} \right].\end{equation}
In standard form, the components of the advective flux are
\begin{equation}
\statevec f_{1}  = \left[ {\begin{array}{*{20}c}
   {\rho v_1}  \\
   {\rho v_1^2  + p}  \\
   {\rho v_1 v_2}  \\
   {\rho v_1 v_3}  \\
   {\rho v_1 H}  \\

 \end{array} } \right]\quad \statevec f_{2}  = \left[ {\begin{array}{*{20}c}
   {\rho v_2}  \\
   {\rho v_2 v_1}  \\
   {\rho v_2 v_2  + p}  \\
   {\rho v_2 v_3}  \\
   {\rho v_2 H}  \\

 \end{array} } \right]\quad \statevec f_{3}  = \left[ {\begin{array}{*{20}c}
   {\rho v_3}  \\
   {\rho v_3 v_1}  \\
   {\rho v_3 v_2}  \\
   {\rho v_3 v_3  + p}  \\
   {\rho v_3 H}  \\

 \end{array} } \right],
 \label{eq:NSEquationFluxes_DAK}
\end{equation}
where $p$ is the pressure and
\begin{equation}
 H = E + \frac{p}
{\rho },\quad E = e + \frac{1}
{2}\left| {\spacevec v} \right|^2, \quad e = \frac{1}
{{\gamma  - 1}}\frac{p}
{\rho }.
\end{equation}
The equations have been scaled with respect to free-stream reference values so that the Reynolds number is
\begin{equation}
\mathrm{Re} = \frac{\rho_{\infty}V_{\infty}L}{\mu_{\infty}},
\end{equation}
where $L$ is the length scale and $V_{\infty}$ is the free-stream velocity. Additionally, the Mach number and Prandtl numbers are 
\begin{equation}{\mathrm{M}_\infty } = \frac{{{V_\infty }}}{{\sqrt {\gamma \mathrm{R}{T_\infty }} }},\quad \Pr  = \frac{{{\mu _\infty }{C_p}}}{{{\lambda _\infty }}}\,.\end{equation}
Written in terms of the primitive variables, the viscous fluxes are
\begin{equation}
  \label{eq:navierstokes_difffluxes}
  \begin{aligned}
  \statevec f^v_{1} &= \left[   0 \quad 
  {{\tau _{11}}} \quad 
  {{\tau _{12}}} \quad 
  {{\tau _{13}}} \quad 
  \left(\left( {\sum}_{j=1}^3 v_j\tau _{1j}\right) + \lambda\frac{\partial T}{\partial x}\right)\;
  \right]^T\,, \\ 
\statevec   f^v_{2} &= \left[
  0 \quad 
  {{\tau _{21}}} \quad 
  {{\tau _{22}}} \quad 
  {{\tau _{23}}} \quad 
  \left( \left( {\sum}_{j=1}^3 v_j\tau _{2j}\right) + \lambda\frac{\partial T}{\partial y}\right)\;
  \right]^T \,,\\ 
\statevec   f^v_{3} &= \left[ 
  0 \quad 
  {{\tau _{31}}} \quad 
  {{\tau _{32}}} \quad 
  {{\tau _{33}}} \quad 
  \left( \left({\sum}_{j=1}^3 v_j\tau _{3j}\right) + \lambda\frac{\partial T}{\partial z}\right)\;
  \right]^T  \,,
  \end{aligned} 
\end{equation}
where
\begin{equation}
\tau _{ij}  = \mu \left( {\frac{{\partial v_j }} {{\partial x^i }} + \frac{{\partial v_i }} {{\partial x^j }}} \right) - \frac{2}{3}\mu \left( {\vecNablaX \cdot\spacevec v} \right)\delta _{ij} \,,\quad \lambda=\frac{\mu}{{(\gamma  - 1)\Pr \mathrm{M}_\infty ^2}}\,,
\end{equation}
and the temperature is
\begin{equation}
T = \gamma \mathrm{M}_{\infty}^{2}\frac{p}{\rho}\,.
\end{equation}

To simplify the presentation, we define \textit{block vectors} (with the double arrow), for instance the block vector flux,
\begin{equation}
\bigstatevec{f} =
 \left[ {\begin{array}{*{20}{c}}
  {{\statevec f_1}} \\ 
  {{\statevec f_2}} \\ 
  {{\statevec f_3}} 
\end{array}} 
\right]\,.
\end{equation}
The spatial gradient of a state vector is a block vector,
\begin{equation}
\vecNablaX \statevec u = \left[ {\begin{array}{*{20}{c}}
  {{\statevec u_x}} \\ 
  {{\statevec u_y}} \\ 
  {{\statevec u_z}} 
\end{array}} 
\right]\,.
\end{equation}
The dot product of two block vectors is defined by
\begin{equation}
\bigstatevec f \cdot \bigstatevec g = \sum\limits_{i = 1}^3 {{{\statevec f}_i}^T{{\statevec g}_i}},\, 
\end{equation}
and the dot product of a block vector with a vector is a state vector,
\begin{equation} 
\spacevec g\cdot\bigstatevec f  = \sum\limits_{i = 1}^3 {{{ g}_i}{{\statevec f}_i}}\,. 
\end{equation}

With this notation the divergence of a flux is defined as
\begin{equation}
\vecNablaX  \cdot \bigstatevec f = \sum\limits_{i = 1}^3 {\frac{{\partial {\statevec f_i}}}{{\partial {x_i}}}}\, ,
\end{equation}
which allows the Navier-Stokes equations to be written compactly as an advection-diffusion like equation
\begin{equation}
  {{\statevec u}_t} + {\vecNablaX} \cdot \bigstatevec f = \overRe{\vecNablaX} \cdot {\bigstatevec f^v}\left( {\statevec u,\vecNablaX\statevec u} \right)  .
  \label{eq:NSEqnsInBlockForm_{DAK}}
\end{equation}

As part of the approximation procedure, it is customary to represent the solution gradients as a new variable to get a first order system of equations
\begin{equation}
 \begin{split}
  \label{eq:navierstokes_mixed}
  {{\statevec u}_t} + {\vecNablaX} \cdot \bigstatevec f &= \overRe{\vecNablaX} \cdot {\bigstatevec{f}^v}( {\statevec u,\bigstatevec q} )   \\
  \bigstatevec q &= {\vecNablaX}\statevec u\,.  
\end{split} 
\end{equation}

To understand the growth of small perturbations in the fluid state, one also studies linearized forms of the Navier-Stokes equations, \eqref{eq:NSEqnsInBlockForm_{DAK}}.
When linearized about a constant state, the Navier-Stokes equations can be written in the form
\begin{equation}
{\statevec u_t} + \sum\limits_{j = 1}^3 {\frac{{\partial {\mmatrix A_j}\statevec u}}{{\partial {x_j}}}}  = \overRe\sum\limits_{i = 1}^3 {\frac{\partial }{{\partial {x_i}}}\left( {\sum\limits_{j = 1}^3 {{\mmatrix B_{ij}}\frac{{\partial \statevec u}}{{\partial {x_j}}}} } \right)},
\label{eq:CNSE}
 \end{equation}
 where $\statevec u = [\delta\rho \;\delta v_{1}\; \delta v_{2}\; \delta v_{3}\; \delta p]^{T}$ is the perturbation from the constant-state reference values. 
The coefficient matrices $\mmatrix A_{j}$ and $\mmatrix B_{ij}$ are constant in the linear approximation of the equations.  

To again simplify the notation for use in the analysis, we define 
a block vector of matrices, e.g.
\begin{equation}\bigstatevec A = \left[ {\begin{array}{*{20}{c}}\
  {{\mmatrix A_1}} \\ 
  {{\mmatrix A_2}} \\ 
  {{\mmatrix A_3}} 
\end{array}} \right],
\end{equation} 
and full block matrix
\begin{equation}
\bigmatrix B = \left[ {\begin{array}{*{20}{c}}
  {{\mmatrix B_{11}}}&{{\mmatrix B_{12}}}&{{\mmatrix B_{13}}} \\ 
  {{\mmatrix B_{21}}}&{{\mmatrix B_{22}}}&{{\mmatrix B_{23}}} \\ 
  {{\mmatrix B_{31}}}&{{\mmatrix B_{32}}}&{{\mmatrix B_{33}}} \\ 
\end{array}} \right].
\label{eq:BlockMatrixB_{DAK}}
\end{equation}
Then the product rule applied to the divergence of the flux in \eqref{eq:CNSE} can be written as
\begin{equation}
{\vecNablaX} \cdot \bigstatevec f = \left({\vecNablaX} \cdot \bigstatevec A\right)\statevec u +  \left(\bigstatevec A\right)^{\!T}{\vecNablaX}\statevec u,
\end{equation}
where 
\begin{equation}
\bigstatevec f = \left[ {\begin{array}{*{20}{c}}
  {{\mmatrix A_1}\statevec u} \\ 
  {{\mmatrix A_2}\statevec u} \\ 
  {{\mmatrix A_3}\statevec u} 
\end{array}} \right],\quad \left(\bigstatevec A\right)^{\!T}=\left[\mmatrix A_1 \;\mmatrix A_2 \;\mmatrix A_3\right ].
\end{equation}
The nonconservative advective form of the linearized Navier-Stokes equations can therefore be written as an advection-diffusion equation
\begin{equation}
\statevec u_{t} + \left({\vecNablaX} \cdot \bigstatevec A\right)\statevec u + \left(\bigstatevec A\right)^{\!T}{\vecNablaX}\statevec u = \overRe\vecNablaX\cdot \left(\bigmatrix B\vecNablaX\statevec u\right).
\end{equation}
Averaging the conservative and nonconservative forms gives a \emph{split form} of the PDE
\begin{equation}
\statevec u_{t} + \frac{1}{2}\left\{{\vecNablaX} \cdot \bigstatevec f+\left({\vecNablaX} \cdot \bigstatevec A\right)\statevec u + \left(\bigstatevec A\right)^{\!T}{\vecNablaX}\statevec u\right\} = \overRe\vecNablaX\cdot \left(\bigmatrix B\vecNablaX\statevec u\right).
\label{eqSplitFormLinear_{DAK}}
\end{equation}
Note that with constant advection matrices, the divergence ${\vecNablaX} \cdot \bigstatevec A$ is zero. 

\CCLsubsection{Boundedness of Energy and Entropy}\label{sec:energy_entropy_continuous}

With suitable boundary conditions, small perturbations,
\begin{equation}
\statevec u = [\delta\rho \;\delta v_{1}\; \delta v_{2}\; \delta v_{3}\; \delta p]^{T},
\end{equation}
are bounded by the initial data in that
\begin{equation}
\inorm{\statevec u(T)}\leqslant C\inorm{\statevec u_{0}},
\end{equation}
where $\inorm{\statevec u}^{2} = \int_{\Omega}| \statevec u |^{2}\,\mathrm{dV}$ is the ``energy norm'' on a domain $\Omega$. Bounding the energy norm guarantees that the individual components of the perturbed state are bounded at any fixed time. They won't blow up.

We show boundedness of the perturbations by multiplying the split form \eqref{eqSplitFormLinear_{DAK}} by an arbitrary $\LTwo{\Omega}$ test function, $\testfuncOne$, and integrating over the domain to get a weak form of the equation. In inner product notation, that weak form is
\begin{equation}
\iprod{\statevec u_{t},\testfuncOne} + \frac{1}{2}\left\{\iprod{\vecNablaX\cdot \bigstatevec f,\testfuncOne}  + \iprod{\left(\bigstatevec A\right)^{\!T}{\vecNablaX}\statevec u,\testfuncOne}\right\} = \overRe\iprod{\vecNablaX\cdot \left(\bigmatrix B\vecNablaX\statevec u\right),\testfuncOne}.
\end{equation}
As with the nonlinear equations in \eqref{eq:navierstokes_mixed}, we introduce the intermediate block vector $\bigstatevec q = \vecNablaX \statevec u$ to produce a first order system
\begin{equation}
 \begin{split}
  \iprod{\statevec u_{t},\testfuncOne} + \frac{1}{2}\left\{\iprod{\vecNablaX\cdot \bigstatevec f,\testfuncOne}  + \iprod{\left(\bigstatevec A\right)^{\!T}{\vecNablaX}\statevec u,\testfuncOne}\right\} &= \overRe\iprod{\vecNablaX\cdot \left(\bigmatrix B\bigstatevec q\right),\testfuncOne}
   \\
  \iprod{\bigstatevec q,\testfuncTwo}  &= \iprod{\vecNablaX \statevec u,\testfuncTwo}   \,,
 \end{split} 
 \label{eq:SplitLinearSystem0}
\end{equation}
where the auxiliary equation for $\bigstatevec{q}$ is multiplied by the test function $\testfuncTwo$ and integrated over the domain.

We then apply the multidimensional integration-by-parts law
to the second and fourth terms, which contain flux divergence, to separate surface (physical boundary) and volume contributions. If, further, we define 
\begin{equation}
 \bigstatevec f^{\,(T)}\left(\testfuncOne \right) = \left[ {\begin{array}{*{20}{c}}
  {{\mmatrix A^{T}_1}\testfuncOne } \\ 
  {{\mmatrix A^{T}_2}\testfuncOne } \\ 
  {{\mmatrix A^{T}_3}\testfuncOne } 
\end{array}} \right]\,,
\label{eq:TestFunctionFlux_{DAK}}
\end{equation}
then we can rewrite the first equation of \eqref{eq:SplitLinearSystem0} as
\begin{equation}
 \begin{split}
\iprod {{{\statevec u}_t},\testfuncOne }
& + \frac{1}{2}\left\{ {\iprod{{\vecNablaX}\statevec u,\bigstatevec{f}^{\,(T)}\left( \testfuncOne  \right)} - \iprod {\bigstatevec f,{\vecNablaX}\testfuncOne } } \right\} 
\\&+ \int\limits_{\partial \Omega } {\left( \frac{1}{2}\left(\bigstatevec f\cdot\spacevec n \right)
- \overRe\left(\left(\bigmatrix B{\vecNablaX}\statevec u\right)\cdot\spacevec n\right) \right)^{\!T}\testfuncOne \dS} 
\\&=  - \overRe\iprod {\bigmatrix B\bigstatevec q,{\vecNablaX}\testfuncOne } ,
 \end{split}
 \label{eq:SplitLinearSystem}
\end{equation}
where
 $\spacevec n$ is the physical space outward normal to the surface.
 
 Looking ahead a few steps, we see that the first term in \eqref{eq:SplitLinearSystem} becomes the time derivative of the energy if we replace the test function $\testfuncOne$ by the solution $\statevec u$. Furthermore, the second term would vanish with the same substitution if the advection matrices $\mmatrix A_i$ were symmetric, leaving only boundary terms, which can be controlled with boundary conditions. This roadmap suggests the need for symmetry in the equations.

The system \eqref{eq:SplitLinearSystem} for the linearized compressible Navier-Stokes equations, although not symmetric, is known to be symmetrizable by a single constant symmetrization matrix, $\mmatrix S$, and there are multiple symmetrizers \citep{Isi:A1981Lw20700001} from which to choose. We denote the symmetrized matrices as
\begin{equation}
\mmatrix A^{s}_{j}=\mmatrix S^{-1}\mmatrix A_{j}\mmatrix S = \left(\mmatrix A^{s}_{j}\right)^{T}\quad\mathrm{and}\quad
\mmatrix B^{s}_{ij}=\mmatrix S^{-1}\mmatrix B_{ij}\mmatrix S = \left(\mmatrix B^{s}_{ij}\right)^{T}\geqslant  0.
\end{equation}
 Explicit representations of the symmetrizer and coefficient matrices are presented in \citet{Nordstrom:2005qy}. Furthermore, the symmetrized block matrix $\bigmatrix B^{s} = \bigmatrix S^{-1}\bigmatrix B\bigmatrix S$,
where
\begin{equation}
\label{eq:symmetrizer_matrix}
\bigmatrix S = \left[ {\begin{array}{*{20}{c}}
  {{\mmatrix S}}\;&0&0 \\ 
  0&{{\mmatrix S}}\;&0 \\ 
  0&0&{{\mmatrix S}}\; 
\end{array}} \right] ,
\end{equation}
is the diagonal block matrix of the symmetrizer, is symmetric and non-negative.


To symmetrize the system and obtain an energy bound at the same time, we let $\testfuncOne = \left(\mmatrix S^{-1}\right)^{T}\mmatrix S^{-1}\statevec u$ in \eqref{eq:SplitLinearSystem}, which includes symmetrization as part of the test function. Then
\begin{equation}
 \begin{split}
\iprod{\mmatrix S^{-1}{{\statevec u}_t},\mmatrix S^{-1}\statevec u } 
&+ \frac{1}{2}\left\{ {\iprod{\vecNablaX\statevec u,\bigstatevec f^{\,(T)}\left( \left(\mmatrix S^{-1}\right)^{T}\mmatrix S^{-1}\statevec u  \right)}}
- \iprod {\bigmatrix S^{-1}\bigstatevec f,{\vecNablaX}\left(\mmatrix S^{-1}\statevec u\right) }  \right\} 
\\&
 +\int\limits_{\partial \Omega } {\left(\frac{1}{2}{\mmatrix S^{-1}\left(\bigstatevec{f}\cdot\spacevec n\right) - \overRe\mmatrix S^{-1}\left(\left(\bigmatrix B{\vecNablaX}\statevec u\right)\cdot\spacevec n\right)} \right)^{T}\mmatrix S^{-1}\statevec u \dS}  
\\&
=  - \overRe\iprod {\bigmatrix S^{-1}\bigmatrix B\bigstatevec q,{\vecNablaX}\left(\mmatrix S^{-1}\statevec u\right )} .
 \end{split}
 \label{eq:SplitLinearSystemStab2}
\end{equation}

To simplify the notation, let us define the symmetric state vector as $\statevec u^{s}=\mmatrix S^{-1}\statevec u$ and examine the terms in \eqref{eq:SplitLinearSystemStab2} separately.
First,
\begin{equation}
\iprod {\mmatrix S^{-1}{{\statevec u}_t},\mmatrix S^{-1}\statevec u }=\frac{1}{2}\frac{d}{dt}\inorm{\statevec u^{s}}^{2}
\label{eq:continuousTimeDepNorm}
\end{equation}
provides the time derivative of the energy of the symmetrized state.
Next, the volume term for the diffusion can be written as
\begin{equation}
\iprod{\bigmatrix S^{-1}\bigmatrix B\bigstatevec q,{\vecNablaX}\left(\mmatrix S^{-1}\statevec u \right)}
= 
\iprod{\bigmatrix B^{s}\bigstatevec q^{s},{\vecNablaX}\statevec u^{s} }.
\end{equation}
Making the changes on the boundary terms,
\begin{equation} 
\begin{split}
&\int\limits_{\partial \Omega } {\left( \frac{1}{2}{\mmatrix S^{-1}\left(\bigstatevec f\cdot\spacevec n\right) - \overRe\mmatrix S^{-1}\left(\left(\bigmatrix B{\vecNablaX}\statevec u\right)\cdot\spacevec n\right )} \right)^{\!T}\mmatrix S^{-1}\statevec u \dS}  \\&=
 \int\limits_{\partial \Omega } 
 {\left( \frac{1}{2}\left(\bigstatevec f^{s}\cdot\spacevec n\right) 
  - \overRe\left(\left(\bigmatrix B^{s}{\vecNablaX}\statevec u^{s}\right)\cdot\spacevec n\right) \right)^{\!T}\statevec u^{s} \dS}  ,
  \end{split}
\end{equation}
where
\begin{equation}
\bigstatevec f^{s} = \left[ {\begin{array}{*{20}{c}}
  {\mmatrix A_1^s\statevec u^{s}} \\ 
  {\mmatrix A_2^s\statevec u^{s}} \\ 
  {\mmatrix A_3^s\statevec u^{s}} 
\end{array}} \right].
\end{equation}

The most interesting terms in \eqref{eq:SplitLinearSystemStab2} are the volume flux terms. The solution flux term is
\begin{equation}
 \iprod{\bigmatrix S^{-1}\bigstatevec f,{\vecNablaX}\left(\mmatrix S^{-1}\statevec u\right) } = \iprod {\bigstatevec f^{s},{\vecNablaX}\statevec u^{s} } ,
 \label{eq:VolumeFluxEnergy_{DAK}}
\end{equation}
and the test function flux term are now the same, for
\begin{equation}
\begin{split}
\iprod{\vecNablaX\statevec u,\bigstatevec f^{\,(T)}\left( \left(\mmatrix S^{-1}\right)^{T}\mmatrix S^{-1}\statevec u  \right)}  &=
\iprod{\bigmatrix S{\vecNablaX}\mmatrix S^{-1}\statevec u,{\left(\bigmatrix S^{-1}\bigstatevec f\left( \statevec u^{s}  \right)\right)^{\!T}}} 
\\&=\iprod{\vecNablaX\mmatrix S^{-1}\statevec u,{\left(\bigmatrix S^{-1}\bigstatevec f\left( \statevec u^{s}  \right)\bigmatrix S\right)^{\!T}}}
\\&= \iprod{\vecNablaX\statevec u^{s},\bigstatevec f^{s}\left( \statevec u^{s}  \right)}.
\end{split}
\label{eq:ExactFluxCancellation}
 \end{equation}

Finally, when we set   $\testfuncTwo = \left(\bigmatrix S^{-1}\right)^{T}\bigmatrix S^{-1}\bigmatrix B\bigstatevec q$ in the second equation of \eqref{eq:SplitLinearSystem0}
\begin{equation}
\begin{split}
  \iprod{\bigstatevec q,\left(\bigmatrix S^{-1}\right)^{T}\bigmatrix S^{-1}\bigmatrix B\bigstatevec q}
  &=  \iprod{\vecNablaX \statevec u,\left(\bigmatrix S^{-1}\right)^{T}\bigmatrix S^{-1}\bigmatrix B\bigstatevec q}
  \end{split}
\end{equation}
we see that
\begin{equation}
\iprod{\vecNablaX \statevec u^{s},\bigmatrix B^{s}\bigstatevec q^{s}}=   \iprod{\bigstatevec q^{s},\bigmatrix B^{s}\bigstatevec q^{s}}\geqslant 0.
\end{equation}

Gathering all the terms, the flux volume terms cancel due to the equivalence of \eqref{eq:VolumeFluxEnergy_{DAK}} and \eqref{eq:ExactFluxCancellation}, leaving
\begin{equation}
\begin{split}
\frac{1}{2}\frac{d}{dt}\inorm{\statevec u^{s}}^{2}&+ \int\limits_{\partial \Omega } {\left( \frac{1}{2}\left(\bigstatevec f^{s}\cdot\spacevec n\right)
   - \overRe\left(\left(\bigmatrix B^{s}{\vecNablaX}\statevec u^{s}\right)\cdot\spacevec n\right) \right)^{\!T}\statevec u^{s} \dS}   \\&= - \overRe\iprod{\bigstatevec q^{s},\bigmatrix B^{s}\bigstatevec q^{s}} \leqslant 0.
\end{split}
\end{equation}
We see, then, that any growth in the energy, defined as the $\mathbb{L}^{2}$ norm, is determined by the boundary integral,
\begin{equation}
\frac{1}{2}\frac{d}{dt}\inorm{\statevec u^{s}}^{2}\leqslant -\int\limits_{\partial \Omega } {\left( \frac{1}{2}\left(\bigstatevec f^{s}\cdot\spacevec n\right)
   - \overRe\left(\left(\bigmatrix B^{s}{\vecNablaX}\statevec u^{s}\right)\cdot\spacevec n\right) \right)^{\!T}\statevec u^{s} \dS} .
   \label{eq:EnergyNormTimeDerivativeContinuous_{DAK}}
\end{equation}
Integrating in time over the interval $[0,T]$,
\begin{equation}\label{eq:linear_analysis_statment_withBCs}
\inorm{\statevec u^{s}(T)}^{2}\leqslant\inorm{\statevec u^{s}(0)}^{2}-\int\limits_{0}^{T}{\int\limits_{\partial \Omega } {\left( \left(\bigstatevec f^{s}\cdot\spacevec n\right)
   - \twooverRe\left(\left(\bigmatrix B^{s}{\vecNablaX}\statevec u^{s}\right)\cdot\spacevec n\right) \right)^{\!T}\statevec u^{s} \dS}}.
\end{equation}

To properly pose the problem, initial and boundary data must be specified. The value at $t=0$ is replaced by initial data $\statevec u^{s}_{0}$. As for the physical boundary terms, \citet{Nordstrom:2005qy} show that the matrices can be split in characteristic fashion into incoming and outgoing information with boundary data specified along the incoming characteristics 
\begin{equation}
\begin{split}
\PBT&= \int\limits_{\partial \Omega } {\left( \left(\bigstatevec f^{s}\cdot\spacevec n\right)
   - \twooverRe\left(\left(\bigmatrix B^{s}{\vecNablaX}\statevec u^{s}\right)\cdot\spacevec n\right) \right)^{\!T}\statevec u^{s}\dS}  \\&=   \int\limits_{\partial \Omega } {{\statevec w^{ + T}}{\Lambda ^ + }{\statevec w^ + }\dS}  - \int\limits_{\partial \Omega } {{\statevec g^T}\left| {{\Lambda ^ - }} \right|\statevec g\dS} ,
   \end{split}
\end{equation}
where $\Lambda^{+}>0$ and $\Lambda^{-}<0$. We will assume here that no energy is introduced by the boundary data, and so set $\statevec g = 0$. As a result,
\begin{equation}
\resizebox{\textwidth}{!}{$\displaystyle{
\PBT = \int\limits_{\partial \Omega } {\left( \left(\bigstatevec f^{s}\cdot\spacevec n\right)
   - \twooverRe\left(\left(\bigmatrix B^{s}{\vecNablaX}\statevec u^{s}\right)\cdot\spacevec n\right) \right)^{\!T}\statevec u^{s}\dS} =  \int\limits_{\partial \Omega } {{\statevec w^{ + T}}{\Lambda ^ + }{\statevec w^ + }\dS}  \geqslant 0 ,
   }$}
\end{equation}
so that
\begin{equation}
\inorm{\statevec u^{s}(T)}\leqslant\inorm{\statevec u^{s}_{0}}.
\end{equation}
Finally, since $\statevec u^{s} = \mmatrix S^{-1}\statevec u$, $\statevec u = \mmatrix S\statevec u^{s}$, 
\begin{equation}
\frac{1}{{{{\left\| \mmatrix S \right\|}_2}}}\left\| \statevec u \right\| \leqslant \left\| {{\statevec u^s}} \right\| \leqslant {\left\| {{\mmatrix S^{ - 1}}} \right\|_2}\left\| \statevec u \right\|,
\end{equation}
where $\|\cdot\|_2$ is the matrix $2$-norm. Therefore, we have the desired result,
\begin{equation}
\inorm{\statevec u(T)}\leqslant C\inorm{\statevec u_{0}}.
\label{eq:ContinuousBound}
\end{equation}

The analysis of the linearized compressible Navier-Stokes equations \eqref{eq:CNSE} provides an $\LTwo{\Omega}$ bound, \eqref{eq:ContinuousBound}, on the solution ``energy''. One might think that an analogous statement of the solution ``energy'' should hold for a nonlinear systems of PDEs. Generally, though, linear stability estimates in the $\LTwo{\Omega}$ sense are insufficient to exclude unphysical solution behavior like expansion shocks \citep{merriam1987}. To eliminate the possibility of such phenomena the notion of the ``energy'' estimate must be generalized for nonlinear systems.

To motivate this generalized solution estimate strategy, we take a slight detour to examine important concepts from thermodynamics. Thermodynamic laws provide rules to decide how physical systems \textbf{cannot} behave and act as guidelines for what solution behavior is physically meaningful and what is not. 

The first law of thermodynamics concerns the conservation of the total energy in a closed system (already present as the fifth equation in the Navier-Stokes equations, \eqref{eq:nse} -- \eqref{eq:NSEquationFluxes_DAK}). The second law states that the entropy of a closed physical system tends to increase over time and, importantly, that it cannot decrease. 

Though somewhat esoteric, the second law of thermodynamics regulates how energies are allowed to transfer within a system. For reversible processes the entropy remains constant over time (isentropic) and the time derivative of the total system entropy is zero. For irreversible processes the entropy increases and that time derivative is positive. Solution dynamics where the total system entropy shrinks in time are never observed and are deemed unphysical. A smooth solution that satisfies the system of nonlinear PDEs, like \eqref{eq:NSEqnsInBlockForm_{DAK}}, corresponds to a reversible process. One of the difficulties, either analytically or numerically, of nonlinear PDEs with a dominant hyperbolic character is that the solution may develop discontinuities (e.g. shocks) regardless of the continuity of the initial conditions \citep{evans2010}. Such a discontinuous solution corresponds to an irreversible process and increases entropy. 

So, the laws of thermodynamics play a pivotal role because they intrinsically provide admissibility criteria and select physically relevant solutions \citep{lax1954,lax1967,Tadmor1987_2}. As given in the compressible Navier-Stokes equations \eqref{eq:NSEqnsInBlockForm_{DAK}}, the total entropy is not part of the state vector of conservative variables $\statevec{u}$. However, we know that the total entropy \textit{is} a conserved quantity for reversible (isentropic) processes. So where is this conservation law ``hiding''? It turns out that there are additional conserved quantities, including the entropy, that are not explicitly built into the nonlinear system, but are still a consequence of the PDE. 
 
To reveal an auxiliary conservation law for the second law of thermodynamics, we define a convex (\emph{mathematical}) entropy function $s = s(\statevec{u})$ that is a scalar function and depends nonlinearly on the conserved variables. From this it will be possible to generalize the previous $\LTwo{\Omega}$ bound for the solution energy $\|\statevec{u}\|$ and instead develop a stability bound on the mathematical entropy of the form
\begin{equation}
\int\limits_\Omega s(\statevec u(T))\,\mathrm{dV}\leqslant \int\limits_\Omega s(\statevec u(0))\,\mathrm{dV} + \PBT,
\label{eq:Entropy_stab_general_statment}
\end{equation}
where $\PBT$ are the physical boundary terms.
This statement of \emph{entropy stability} \eqref{eq:Entropy_stab_general_statment} provides a bound on the entropy function in terms of the initial condition and appropriate boundary conditions, analogous to the linear bound \eqref{eq:linear_analysis_statment_withBCs}. Further, as the mathematical entropy is a convex function of the solution $\statevec u$, the entropy bound also leads to a bound on an associated norm of $\statevec u$ \citep{merriam1987,Dutt:1988}.

Entropy stability for a single entropy does not give nonlinear stability, but it does give a stronger estimate than linear stability \citep{merriam1989,tadmor2003}, which is formally only appropriate for small perturbations to the equations. For the compressible Navier-Stokes equations an appropriate entropy pair $(s,\spacevec{f}^\ent)$, consists of the scalar entropy function 
\begin{equation}
\label{eq:entropy_function}
s = s(\statevec{u}) = -\frac{\rho \varsigma}{\gamma-1} =-\frac{\rho\left(\ln{p}-\gamma\ln{\rho}\right)}{\gamma-1}\,,
\end{equation}
where $\varsigma = \ln{p}-\gamma\ln{\rho}$ is the physical entropy, and the associated entropy flux
\begin{equation}
\spacevec{f}^{\halfComma\ent} = s\,\spacevec{v}\,.
\end{equation}

Note that the mathematical entropy $s$ is taken as the negative of the physical entropy so that the mathematical entropy is bounded in time, just like the energy measured by the $\mathbb L^2$ norm. This allows us to write the more mathematically common type of bound \eqref{eq:Entropy_stab_general_statment}. From here on we will use the term entropy to refer to the mathematical entropy, not the physical, thermodynamic entropy.

We also introduce entropy variables, the vector $\statevec w$ being the derivative of the entropy with respect to the conservative state variables,  
\begin{equation}\label{eq:entVariables}
\statevec{w} = \frac{\partial s}{\partial\statevec{u}} = \left[
      \frac{\gamma - \varsigma}{\gamma -1} - \frac{\rho|\spacevec{v}|^2}{2p},\quad
      \frac{\rho v_1}{p},\quad 
      \frac{\rho v_2}{p},\quad
      \frac{\rho v_3}{p},\quad
      \frac{-\rho}{p}\right]^T,
\end{equation}
with the convexity property
\begin{equation}
\statevec{k}^T\frac{\partial^2 s}{\partial\statevec{u}^2}\statevec{k} >0,\quad \forall \statevec{k}\neq 0,
\end{equation}
if $\rho>0$ and $p>0$ \citep{carpenter_esdg,tadmor2003,Dutt:1988}. 

The positivity requirement on the density and the temperature, $T\propto p/\rho$, ensures a one-to-one mapping between conservative and entropy variables. This constraint is unfortunately not a by-product of the entropy stability estimate for the thermodynamic entropy. Hence, entropy stability is not a true nonlinear stability statement and further criteria (up to this point unknown for the three-dimensional compressible Navier-Stokes equations) are necessary. Consequently, entropy stable discretizations can (and do) produce invalid solutions with negative density or temperature and need additional strategies to guarantee positivity. 

The entropy variables are introduced because they contract the entropy pair, meaning that they satisfy the relations
\begin{equation}
\statevec{w}^T\,\statevec{u}_t = \left(\frac{\partial s}{\partial\statevec{u}}\right)^{\!T}\statevec{u}_t  = s_t(\statevec{u}),
\label{eq:wsContraction}
\end{equation}
and
\begin{equation}
\statevec{w}^T\, \vecNablaX\cdot\bigstatevec{f} = \vecNablaX\cdot\spacevec{f}^{\halfComma\ent}.
\label{eq:spatial_ent_contract}
\end{equation}
The contraction allows us to convert a system of advection equations (in this instance the compressible Euler equations)
\begin{equation}
\statevec u_t +  \vecNablaX\cdot\bigstatevec{f} = 0
\end{equation}
to a scalar advection equation for the entropy simply by multiplying by the entropy variables $\statevec w$,
\begin{equation}
\statevec w^T\left(\statevec u_t +  \vecNablaX\cdot\bigstatevec{f}\right) = s_t +\vecNablaX\cdot\spacevec{f}^{\halfComma\ent}  = 0.
\end{equation}
 (Cf. how multiplying by the solution state $\statevec u$ in the linear analysis above converts the system to a scalar equation for the mathematical energy.)

Furthermore, the viscous flux can be rewritten in terms of the gradient of the entropy variables
\begin{equation}
\label{eq:visc_flux_alt}
\bigstatevec{f}^v\left({\statevec u,{\vecNablaX}\statevec u}\right) = \bigmatrix{B}^\ent\,\vecNablaX\statevec{w},
\end{equation}
where $\bigmatrix{B}^\ent$ satisfies
\begin{equation}
\label{eq:spd_viscous_terms}
\mmatrix{B}^\ent_{ij} = (\mmatrix{B}^\ent_{ji})^T,\qquad \sum\limits_{i=1}^d\sum\limits_{j=1}^d\left(\frac{\partial \statevec w}{\partial x_i}\right)^T\,\mmatrix{B}^\ent_{ij}\,\left(\frac{\partial \statevec w}{\partial x_j}\right)\geqslant 0,\quad \forall \statevec w,
\end{equation}
if $p>0$ and $\mu>0$ \citep{carpenter_esdg,tadmor2006,Dutt:1988}. Formally, the entropy variables \eqref{eq:entVariables} can be used to rewrite the compressible Navier-Stokes equations into a symmetric and nonnegative form as shown by \citet{Dutt:1988}, again analogous to the linear analysis symmetrization procedure.

Using the contraction properties of the entropy variables, we can construct a bound on the mathematical entropy of the form
\begin{equation}
\overline{s}(T)\leqslant \overline{s}(0)
\label{eq:FinalEntropyBound_{DAK}}
\end{equation}
where 
\begin{equation}
\overline{s} = \iprod{s(\statevec{u}),1}=\int\limits_{\Omega} s(\statevec{u}) \,\mathrm{dV}
\label{eq:EntropyBound_{DAK}}
\end{equation}
is the total entropy,
provided that suitable boundary conditions are applied.

We find the bound \eqref{eq:FinalEntropyBound_{DAK}} much as we did when we found the energy bound for the linear system. We multiply the first equation by the entropy variables and the second equation by the viscous flux, and integrate over the domain to get the weak forms
\begin{equation}
\begin{split}
\iprod{\statevec{w}(\statevec u), {\statevec u_t}} + \iprod{\statevec{w}(\statevec u),\vecNablaX\cdot\bigstatevec{f}} &= \overRe\,\iprod{\statevec{w}(\statevec u),\vecNablaX\cdot\bigstatevec{f}^v} ,\\
\iprod{\bigstatevec{q}, \bigstatevec{f}^v}&=\iprod{\vecNablaX\statevec{w},\bigstatevec{f}^v}.
\end{split}
\label{eq:weakformNS1}
\end{equation}
Next we use the properties of the entropy pair \eqref{eq:wsContraction} and \eqref{eq:spatial_ent_contract} to contract the left side of the first equation of \eqref{eq:weakformNS1}, and use multidimensional integration-by-parts on the right hand side to get
\begin{equation}
\begin{split}
\iprod{s_t(\statevec{u}),1} + \iprod{\vecNablaX\cdot\spacevec{f}^{\halfComma\ent},1} &= \overRe\,\int\limits_{\partial\Omega}\statevec{w}^T(\statevec u)\,\left(\bigstatevec{f}^v\cdot\spacevec{n}\right) \dS  - \overRe\,\iprod{\vecNablaX\statevec{w}(\statevec u),\bigstatevec{f}^v} ,\\
\iprod{\bigstatevec{q}, \bigstatevec{f}^v}&=\iprod{\vecNablaX\statevec{w},\bigstatevec{f}^v}.
\end{split}
\label{eq:weakformNS2}
\end{equation}
Inserting the second equation of \eqref{eq:weakformNS2} into the first and applying the identity \eqref{eq:visc_flux_alt} gives
\begin{equation}
\begin{split}
\iprod{s_t(\statevec{u}),1} + \iprod{\vecNablaX\cdot\spacevec{f}^{\halfComma\ent},1} &= \overRe\,\int\limits_{\partial\Omega}\statevec{w}^T(\statevec u)\,\left(\bigstatevec{f}^v\cdot\spacevec{n}\right) \dS  - \overRe\,\iprod{\bigstatevec{q}, \bigmatrix{B}^\ent\,\bigstatevec{q}}.
\end{split}
\end{equation}
Finally, we use the property \eqref{eq:spd_viscous_terms} as well as multidimensional integration-by-parts to integrate the flux divergence on the left side to get the estimate
\begin{equation}
\label{eq:continuous_NSE_entropy_estimate}
\begin{split}
\frac{d}{dt}\overline{s}  \leqslant \int\limits_{\partial\Omega}\left( -\spacevec{f}^{\halfComma\ent} \cdot\spacevec{n} + \overRe\,\statevec{w}^T(\statevec u)\left(\bigstatevec{f}^v\cdot\spacevec{n}\right)\right)\dS.
\end{split}
\end{equation}
This entropy estimate is precisely that given previously in \eqref{eq:Entropy_stab_general_statment} except the form of the physical boundary terms is now explicitly given for the compressible Navier-Stokes equations. Boundary conditions then need to be specified so that the bound on the entropy depends only on the boundary data. We will assume here that boundary data are given so that the right hand side is non-positive so that the entropy will not increase in time. For more thorough discussion on boundary conditions for the Navier-Stokes equations, see, e.g. \citep{Dutt:1988,dalcin2019conservative,hindenlang2019stability}. Integrating \eqref{eq:continuous_NSE_entropy_estimate} in time then gives the desired result, \eqref{eq:FinalEntropyBound_{DAK}}.

\CCLsection{Construction of Curvilinear Spectral Elements}\label{sec:section_SEM_GG}
The general goal for the nodal DG method is to use the Lagrange polynomial basis with Gauss-Lobatto interpolation nodes to approximate the solutions as high-order polynomial interpolants, \eqref{eq:Interpolant2D} or \eqref{eq:Interpolant3D}. This appears to restrict the approximation to the simple quadrilateral $E_{2}$ or hexahedron $E_{3}$. To overcome this severe limitation, we use a process to extend the methods to completely general geometries.
That process consists of three steps:
\begin{enumerate}
\item The domain $\Omega$ is subdivided into quadrilateral or hexahedral elements, $e^{k}$, $k=1,2,\ldots, K$.
\item A mapping is created from the computational space coordinate $\spacevec{\xi}$ on the reference element, $E_{2}$ or $E_{3}$, onto the physical space coordinate $\spacevec{x}$ for each element $e^{k}$. 
\item The equations are re-written in terms of the computational space coordinate on the reference element.
\end{enumerate}
It is on the reference element with the mapped equations that the DG approximation is then created using the spectral approximation tools derived in Sec. \ref{sec:spectral_calc}. The result will be a DG spectral element method to approximate the solution of conservation laws in three-dimensional geometries, e.g. \eqref{eq:NSEqnsInBlockForm_{DAK}}, that has as many of the properties of the continuous equations as possible, e.g. \eqref{eq:ContinuousBound} or \eqref{eq:FinalEntropyBound_{DAK}}. 

A significant advantage of approximating the equations in computational (or reference) space is that they can be derived independently of the element shape and depend only on the transformation defined in Step 2. A specific advantage is that high-order spectral approximations for the reference domain $E_{3}$ have been previously described in Sec.~\ref{sec:InterpMultipleSpace_DAK}, with spectral accuracy coming from the Gauss-Lobatto quadrature and the Lagrange polynomial basis ansatz.

\CCLsubsection{Subdividing the Domain: Spectral Element Mesh Generation}\label{sec:meshing}

\begin{figure}[!ht]
 \begin{center}
    \begin{minipage}[b]{\textwidth}\centering
       \includegraphics[scale = 0.1465]{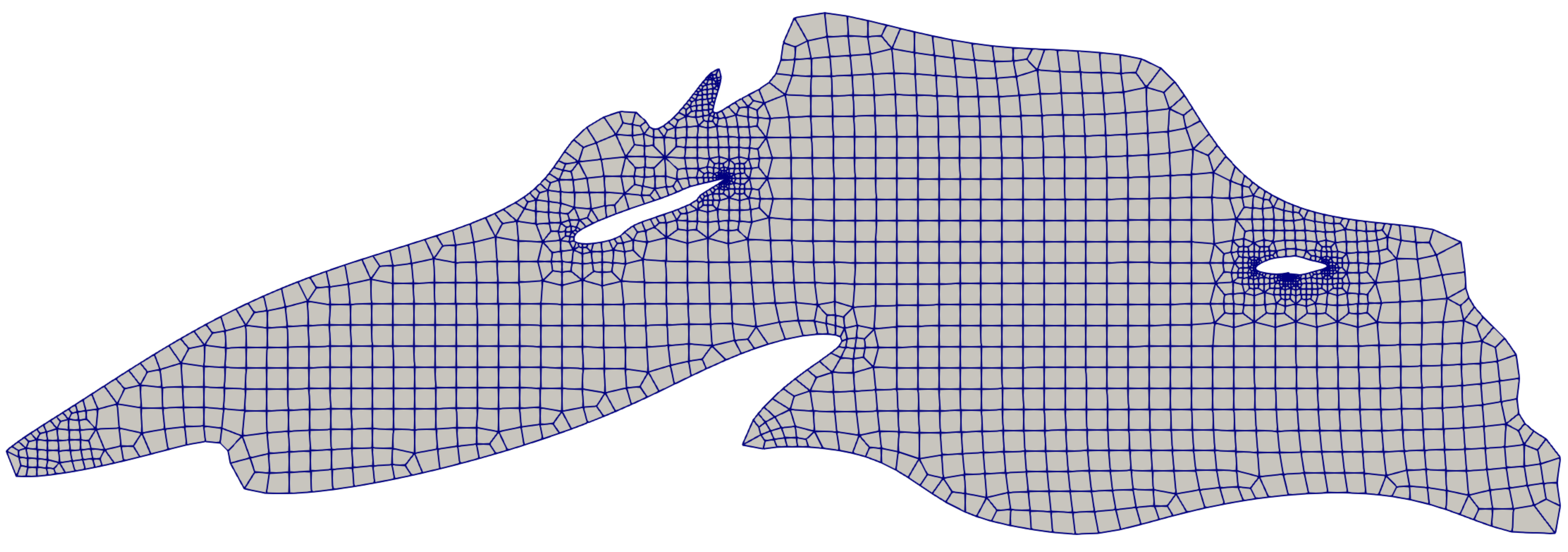}\\
       {(a)} Quadrilateral mesh of Lake Superior.
     \end{minipage}
     \vspace{0.25cm}
 \\[-3ex]
     \begin{minipage}[b]{\textwidth}\centering
       \includegraphics[scale=0.144]{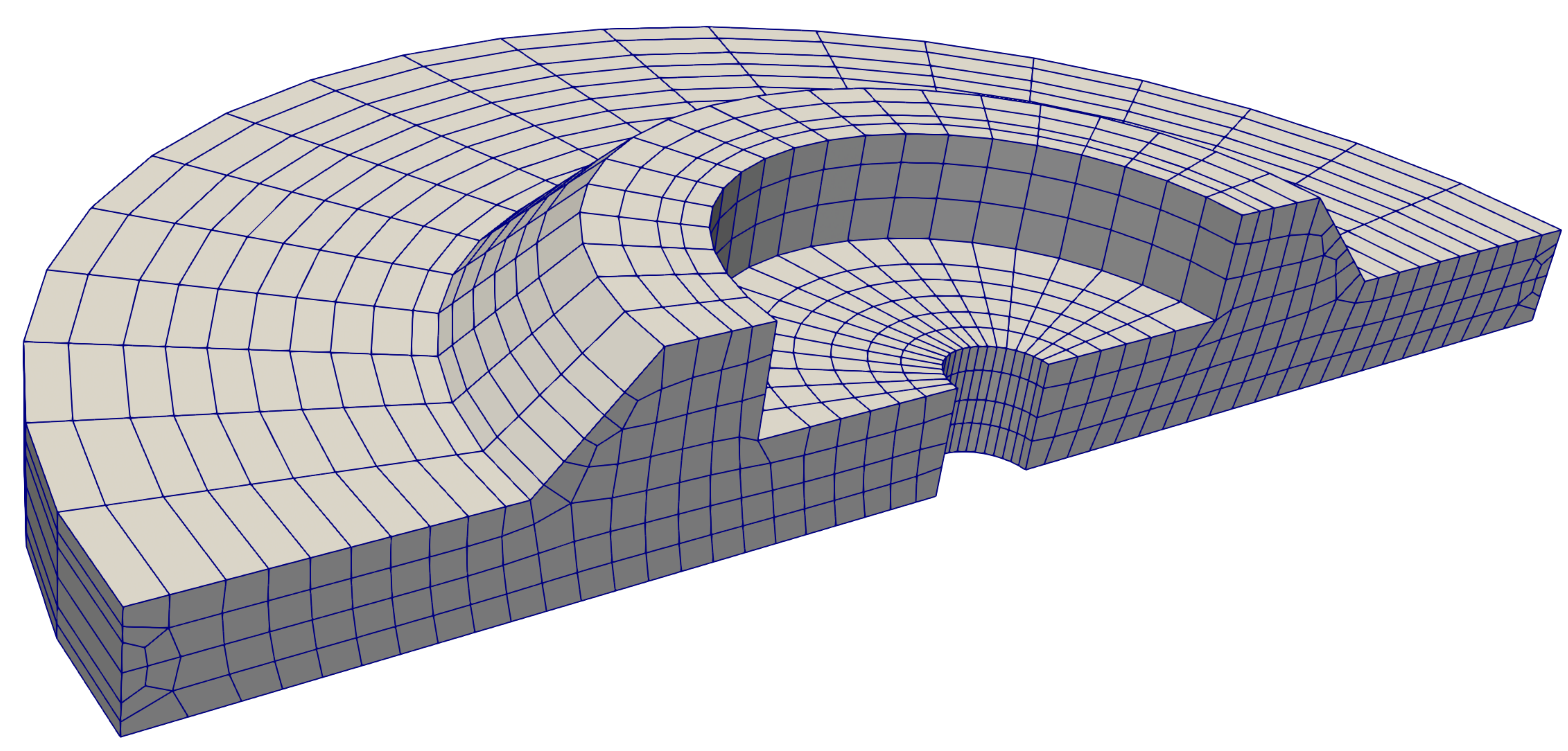}\\[-0.25ex]
       {(b)} Hexahedral mesh.
     \end{minipage}
  \caption{Example meshes in (a) two and (b) three spatial dimensions.}
  \label{fig:mesh_examples}
 \end{center}
\end{figure}

The first step in the approximation is to subdivide a domain $\Omega$ into a mesh of non-overlapping elements $e^{k}$, $k=1,2,\ldots, K$. We restrict the discussion to quadrilateral or hexahedral elements because the forthcoming DG approximation will be built from a tensor product ansatz as in Sec.~\ref{sec:InterpMultipleSpace_DAK}. Two examples of such meshes are given in Figure \ref{fig:mesh_examples}. The generation of such meshes, especially with curved boundary information, is outside the scope of this chapter. As mentioned in the Prologue, high-order mesh generation is a difficult task with many open issues. For instance it is necessary that the elements are valid and non-inverted, with non-negative mapping Jacobians. This, however, is non-trivial for boundary layer meshes with high aspect ratios near boundaries with high curvature. We refer the interested reader to \citet{geuzaine2009gmsh,Hindenlang:2014gl} and the references therein for more details.

\CCLsubsection{Mapping Elements from the Reference Element}\label{sec:element_mapping}
 
 More complex elements than $E_{2}$ or $E_{3}$ that have curved boundaries can be accommodated by a transformation from the reference elements onto the physical elements $e^{k}$. We adopt the naming convention where the domain $E_{d}$ is called the \emph{reference element} or \emph{computational domain} and the element onto which computational domain is mapped is the \emph{physical domain}.

From the mesh in the previous section, the physical domain has been divided into a set of (possibly curved) elements $\{e^k\}_{k=1}^K$. As before, the physical domain coordinates will be denoted by $\spacevec{x} = (x,y,z)^T$ within an element $e^k$. Analogously, the computational domain coordinates are defined as
\begin{equation}
\spacevec{\xi} = (\xi\,,\,\eta\,,\,\zeta)^T = \left(\xi^1\,,\,\xi^2\,,\,\xi^3\right)^T,
\end{equation}
in the reference element $E_3$. Points in the reference element are mapped to each of the element, $e^k$, in physical space, with a polynomial mapping 
\begin{equation}\label{eq:generalMapping}
\spacevec{x} = \spacevec{X}^{k}(\spacevec{\xi}).
\end{equation}
In the following we will ignore the index, $k$, and realize that all expressions relate to any given element, $e^{k}$.

To maintain generality, \eqref{eq:generalMapping} is an algebraic transformation that maps the boundaries of the reference element to the boundaries of the physical element and interior to interior. In this section, we demonstrate how to create a three-dimensional transformation from $E_3$ to a curved hexahedron. The two dimensional mappings for quadrilateral elements are described in the book by \citet{Kopriva:2009nx}.

The most common approach to generate the mapping $\spacevec{X}(\spacevec{\xi})$ is to use transfinite interpolation introduced by \citet{Gordon&Hall1973a}. The idea is to interpolate between (possibly curved) boundaries with a polynomial to guarantee a smooth transformation between the computational and physical domains. The simplest transfinite interpolation, and the one almost always used in practice, is the linear blending formula, which uses a linear interpolation between boundaries. 

\begin{figure}[!ht]
\begin{center}
\includegraphics[width=\textwidth]{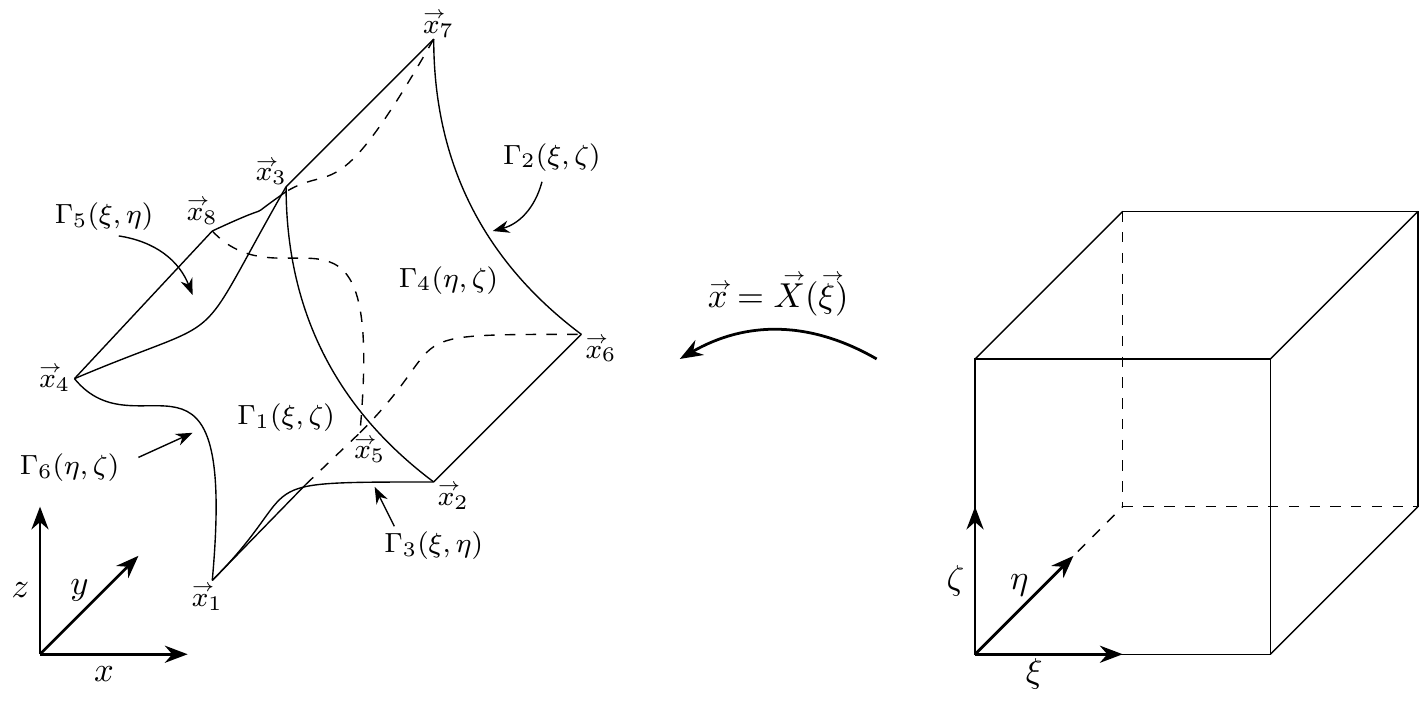}
\caption{Example mapping from computational to physical coordinates.}
\label{fig:3D_mapping}
\end{center}
\end{figure}

In three space dimensions the physical domain is bounded by six curved faces $\Gamma_i$, $i = 1,\ldots,6$, as depicted in Figure \ref{fig:3D_mapping}.
Although it may be possible to define the boundary curves through analytic functions, we show later that constraints like free-stream preservation require that the curves are polynomials in their arguments. As a result, the faces are approximated by polynomials of degree $N$, written in the Lagrange basis.  For example the third boundary face $\Gamma_3$ is approximated as
\begin{equation}
\Gamma_3\approx \IN{\Gamma_3} = \sum_{i,j=0}^N\left(\Gamma_3\right)_{ij}\ell_i(\xi)\ell_j(\eta).
 \end{equation}
Approximating the boundary to the same polynomial order, $N$, as the solution is called \emph{isoparametric}.

The transformation is derived by linear interpolation between opposing faces. As such, one first creates a linear interpolation between two faces, say $ \Gamma_3$ and $ \Gamma_5$
\begin{equation}
\spacevec{X}_{35}(\xi,\eta) = \frac{1}{2}\{(1-\zeta)\Gamma_3(\xi,\eta) + (1+\zeta)\Gamma_5(\xi,\eta)\}.
\end{equation}
Similarly, linear interpolation is constructed for the other four faces as
\begin{equation}
\begin{aligned}
\spacevec{X}_{12}(\xi,\zeta) = \frac{1}{2}\{(1-\eta)\Gamma_1(\xi,\zeta) + (1+\eta)\Gamma_2(\xi,\zeta)\}, \\
\spacevec{X}_{64}(\eta,\zeta) = \frac{1}{2}\{(1-\xi)\Gamma_6(\eta,\zeta) + (1+\xi)\Gamma_4(\eta,\zeta)\}.
\end{aligned}
\end{equation}
The final mapping will be a combination of the six face interpolants and three linear interpolations between them, starting with the sum
\begin{equation}\label{3DMappingCombo}
\begin{aligned}
\spacevec\Sigma(\xi,\eta,\zeta) &= \frac{1}{2}\left\{(1-\xi)\Gamma_6(\eta,\zeta) + (1+\xi)\Gamma_4(\eta,\zeta) + (1-\eta)\Gamma_1(\xi,\zeta)\right. \\
                                  &\;\;\;\; + \left.(1+\eta)\Gamma_2(\xi,\zeta) +(1-\zeta)\Gamma_3(\xi,\eta) + (1+\zeta)\Gamma_5(\xi,\eta)\right\}.
\end{aligned}
\end{equation}
Unfortunately, the combination \eqref{3DMappingCombo} no longer always matches at the faces
\begin{align}
\spacevec\Sigma(-1,\eta,\zeta) = \Gamma_6(\eta,\zeta) \;+\; &\frac{1}{2}\left\{(1-\eta)\Gamma_1(-1,\zeta)+ (1+\eta)\Gamma_2(-1,\zeta)\right.\nonumber\\
                                  &\;\; +\left.(1-\zeta)\Gamma_3(-1,\eta) + (1+\zeta)\Gamma_5(-1,\eta)\right\},\label{extraTerms3DMapping1}\\[0.1cm]
\spacevec\Sigma(1,\eta,\zeta) = \Gamma_4(\eta,\zeta) \;+\; &\frac{1}{2}\left\{(1-\eta)\Gamma_1(1,\zeta)+ (1+\eta)\Gamma_2(1,\zeta)\right.\nonumber\\
                                  &\;\; +\left.(1-\zeta)\Gamma_3(1,\eta) + (1+\zeta)\Gamma_5(1,\eta)\right\}, \label{extraTerms3DMapping2}\\[0.1cm]
\spacevec\Sigma(\xi,-1,\zeta) = \Gamma_1(\xi,\zeta) \;+\; &\frac{1}{2}\left\{(1-\xi)\Gamma_6(-1,\zeta) + (1+\xi)\Gamma_4(-1,\zeta)\right.\nonumber\\
                                  &\;\; +\left.(1-\zeta)\Gamma_3(\xi,-1) + (1+\zeta)\Gamma_5(\xi,-1)\right\}, \label{extraTerms3DMapping3}\\[0.1cm]
\spacevec\Sigma(\xi,1,\zeta) = \Gamma_2(\xi,\zeta) \;+\; &\frac{1}{2}\left\{(1-\xi)\Gamma_6(1,\zeta) + (1+\xi)\Gamma_4(1,\zeta)\right.\nonumber\\
                                  &\;\; +\left.(1-\zeta)\Gamma_3(\xi,1) + (1+\zeta)\Gamma_5(\xi,1)\right\}, \label{extraTerms3DMapping4}\\[0.1cm]
\spacevec\Sigma(\xi,\eta,-1) = \Gamma_3(\xi,\eta) \;+\; &\frac{1}{2}\left\{(1-\eta)\Gamma_1(\xi,-1)+ (1+\eta)\Gamma_2(\xi,-1)\right.\nonumber\\
                                  &\;\; +\left.((1-\xi)\Gamma_6(\eta,-1) + (1+\xi)\Gamma_4(\eta,-1)\right\},  \label{extraTerms3DMapping5}\\[0.1cm]
\spacevec\Sigma(\xi,\eta,1) = \Gamma_5(\xi,\eta) \;+\; &\frac{1}{2}\left\{(1-\eta)\Gamma_1(\xi,1)+ (1+\eta)\Gamma_2(\xi,1)\right.\nonumber\\
                                  &\;\; +\left.(1-\xi)\Gamma_6(\eta,1) + (1+\xi)\Gamma_4(\eta,1)\right\}.  \label{extraTerms3DMapping6}
\end{align}
To match the faces, correction terms must be subtracted in the $\xi$, $\eta$, and $\zeta$ directions to cancel the additional terms that appear in the braces of \eqref{extraTerms3DMapping1} -- \eqref{extraTerms3DMapping6}. These linear corrections are
\begin{equation}\label{xiEdgeCorrection}
\resizebox{\textwidth}{!}{$
\begin{aligned}
\spacevec{C}^{\xi} = &\left(\frac{1-\xi}{4}\right)\big[(1-\eta)\Gamma_1(-1,\zeta)+ (1+\eta)\Gamma_2(-1,\zeta)+(1-\zeta)\Gamma_3(-1,\eta)+ (1+\zeta)\Gamma_5(-1,\eta)\big] \\[0.05cm]
+&\left(\frac{1+\xi}{4}\right)\big[(1-\eta)\Gamma_1(1,\zeta)+ (1+\eta)\Gamma_2(1,\zeta)+(1-\zeta)\Gamma_3(1,\eta) + (1+\zeta)\Gamma_5(1,\eta)\big],
\end{aligned}
$}
\end{equation}
\begin{equation}\label{etaEdgeCorrection}
\resizebox{\textwidth}{!}{$
\begin{aligned}
\spacevec{C}^{\eta} = &\left(\frac{1-\eta}{4}\right)\big[(1-\xi)\Gamma_6(-1,\zeta) + (1+\xi)\Gamma_4(-1,\zeta)+(1-\zeta)\Gamma_3(\xi,-1) + (1+\zeta)\Gamma_5(\xi,-1)\big]\\[0.05cm]
 + &\left(\frac{1+\eta}{4}\right)\big[(1-\xi)\Gamma_6(1,\zeta) + (1+\xi)\Gamma_4(1,\zeta)+(1-\zeta)\Gamma_3(\xi,1) + (1+\zeta)\Gamma_5(\xi,1)\big],
\end{aligned}
$}
\end{equation}
and
\begin{equation}\label{zetaEdgeCorrection}
\resizebox{\textwidth}{!}{$
\begin{aligned}
\spacevec{C}^{\zeta} = &\left(\frac{1-\zeta}{4}\right)\big[(1-\eta)\Gamma_1(\xi,-1)+ (1+\eta)\Gamma_2(\xi,-1)+(1-\xi)\Gamma_6(\eta,-1) + (1+\xi)\Gamma_4(\eta,-1)\big]\\
 +& \left(\frac{1+\zeta}{4}\right)\big[(1-\eta)\Gamma_1(\xi,1) + (1+\eta)\Gamma_2(\xi,1)+(1-\xi)\Gamma_6(\eta,1) + (1+\xi)\Gamma_4(\eta,1)\big]. \\
\end{aligned}
$}
\end{equation}
However, subtracting the correction terms \eqref{xiEdgeCorrection}, \eqref{etaEdgeCorrection}, and \eqref{zetaEdgeCorrection} from \eqref{3DMappingCombo} removes the interior contribution twice. Thus, to complete the correction to \eqref{3DMappingCombo}, one adds back the transfinite map of the reference cube to a straight-sided hexahedral element,
\begin{equation}\label{StandardHexMap}
\begin{aligned}
\spacevec{X}_{H}(\spacevec\xi) =  &\frac{1}{8}\left\{\spacevec{x}_1(1-\xi)(1-\eta)(1-\zeta)+\spacevec{x}_2(1+\xi)(1-\eta)(1-\zeta)\right.\\
            & + \spacevec{x}_3(1+\xi)(1+\eta)(1-\zeta)+\spacevec{x}_4(1-\xi)(1+\eta)(1-\zeta) \\
            & + \spacevec{x}_5(1-\xi)(1-\eta)(1+\zeta)+\spacevec{x}_6(1+\xi)(1-\eta)(1+\zeta) \\
            & + \left.\spacevec{x}_7(1+\xi)(1+\eta)(1+\zeta)+\spacevec{x}_8(1-\xi)(1+\eta)(1+\zeta)\right\},
\end{aligned}
\end{equation}
where $\spacevec{x}_i$, $i=1,\ldots,8$ are the locations of the corners of the hexahedron. 

The final transfinite interpolation with linear blending for a curved-sided hexahedron is therefore
\begin{equation}\label{3DTransfinite}
\spacevec{X}(\spacevec\xi) = \spacevec\Sigma(\spacevec\xi) -\frac{1}{2}\left[\spacevec{C}^{\xi}+\spacevec{C}^{\eta}+\spacevec{C}^{\zeta}\right] + \spacevec{X}_H(\spacevec\xi),
\end{equation}
where the correction terms  \eqref{xiEdgeCorrection}, \eqref{etaEdgeCorrection}, and \eqref{zetaEdgeCorrection} are further divided by two, otherwise they would contribute at each of the twelve edges twice.

\CCLsubsection{Transforming Equations from Physical to Reference Domains}\label{sec:equation_mapping}

%

The mapping \eqref{3DTransfinite} provides a mechanism to connect differential operators in the computational domain to the physical domain. Under the mapping, the equations themselves, e.g. \eqref{eq:NSEqnsInBlockForm_{DAK}}, are transformed as well, essentially being an exercise of the chain rule. Specifically, the differential operators of the divergence, gradient and curl change form due to the mapping. 

Rather than simply apply the chain rule, we summarize a general approach that uses ideas from differential geometry to transform equations between the reference and physical coordinate systems. This approach better exposes properties of the transformations that should be satisfied by the approximation. Full discussions of these general derivations can be found in \citet{Farrashkhalvat:2003pd,Knupp:1993bh} and in \citet{hesthaven2008,Kopriva:2009nx} for the particular context of spectral methods.

The differential transformations are described in terms of two sets of independent coordinate basis vectors. The first is the \textit{covariant basis}
\begin{equation}
\spacevec{a}_{\halfComma i}=\pderivative{\spacevec{X}}{\xi^{i}}\quad i=1,2,3,
\label{eq:covariant_vectors}
\end{equation}
whose components lie tangent to the transformation of a coordinate line in the computational space.
Conveniently, the covariant basis vectors can be computed directly from the mapping between the reference element and physical space, $\spacevec X$, \eqref{3DTransfinite}. 

The second basis
is the \textit{contravariant basis}, whose components are normal to the transformation of the coordinate lines,
\begin{equation}
\spacevec{a}^{\halfComma i}=\vecNablaX\xi^{i},\; i=1,2,3.
\label{eq:contravariantBasisDef}
\end{equation}
The contravariant basis vectors, for instance, point in the direction of the normal at a physical boundary. These two bases are not necessarily orthogonal, and will not be, unless the transformation is conformal. At first glance, it appears that the inverse mapping $\spacevec{\xi} = \spacevec{X}^{-1}(\spacevec{x})$ is needed to compute the contravariant basis vectors $\spacevec{a}^{\halfComma i},\; i=1,2,3.$ But this is not the case, once we have a way to represent the gradient in reference space.

A differential surface element can be written terms of the reference space coordinates by way of the cross product
\begin{equation}\label{eq:co_surf_elem}
\text{dS}^i = \spacevec{a}_j\text{d}\xi^j \times \spacevec{a}_k\text{d}\xi^k = (\spacevec{a}_j\times\spacevec{a}_k)\text{d}\xi^j\text{d}\xi^k, \quad (i,j,k)\;\text{ cyclic} 
\end{equation}
from which a volume element can be generated by extending the the surface element \eqref{eq:co_surf_elem} in the normal direction
\begin{equation}\label{eq:co_vol_elem}
\text{dV} = \spacevec{a}_i\cdot(\spacevec{a}_j\times\spacevec{a}_k)\text{d}\spacevec{\xi} = J\text{d}\spacevec{\xi}, \quad (i,j,k)\;\text{ cyclic}.
\end{equation}
Writing the volume element this way exposes the Jacobian of the transformation in terms of the covariant basis vectors,
\begin{equation}
J = \spacevec{a}_1\cdot(\spacevec{a}_2\times\spacevec{a}_3).
\end{equation}

Using the usual pillbox approach, the divergence is derived from the surface and volume differentials as
\begin{equation}\label{eq:divForm_map}
\vecNablaX\cdot\spacevec{f} = \frac{1}{J}\sum_{i=1}^3\pderivative{}{\xi^i}\left((\spacevec{a}_j\times\spacevec{a}_k)\cdot\spacevec{f}\,\right).
\end{equation}

From the divergence it is possible to find an important identity satisfied by the covariant basis vectors. Under the assumption that the flux vector $\spacevec f$ is an arbitrary constant state, i.e. $\spacevec{f} = \spacevec{c}$, \eqref{eq:divForm_map} simplifies to
\begin{equation}\label{eq:metricIDs_weirdForm}
0 = \sum_{i=1}^3\pderivative{}{\xi^i}(\spacevec{a}_j\times\spacevec{a}_k).
\end{equation}
The statement \eqref{eq:metricIDs_weirdForm} is one form of the \textit{metric identities}. 
From \eqref{eq:metricIDs_weirdForm} it is possible to rewrite the divergence \eqref{eq:divForm_map} into an equivalent form
\begin{equation}\label{eq:divForm_map_nonCons}
\vecNablaX\cdot\spacevec{f} = \frac{1}{J}\sum_{i=1}^3(\spacevec{a}_j\times\spacevec{a}_k)\cdot\pderivative{\spacevec{f}}{\xi^i}.
\end{equation}

From the alternative form of the divergence \eqref{eq:divForm_map_nonCons} it is straightforward to see that the \textit{gradient} of some scalar function $g$ in reference coordinates is
\begin{equation}\label{eq:gradForm_map_nonCons}
\vecNablaX g = \frac{1}{J}\sum_{i=1}^3(\spacevec{a}_j\times\spacevec{a}_k)\pderivative{g}{\xi^i}.
\end{equation}

If we replace $g$ by $g = \xi^i$, $i = 1,2,3$ in the gradient, \eqref{eq:gradForm_map_nonCons}, we relate the contravariant vectors \eqref{eq:contravariantBasisDef} to the covariant
\begin{equation}
\vecNablaX \xi^i = \frac{1}{J}\sum_{m=1}^3(\spacevec{a}_j\times\spacevec{a}_k)\pderivative{\xi^i}{\xi^m}.
\end{equation}
But $\pderivative{\xi^i}{\xi^m} = \delta_{im}$ so the sum simplifies to a definition of the volume weighted contravariant vectors in terms of the covariant,
\begin{equation}\label{eq:contra_vec_map}
J\vecNablaX \xi^i =  J\spacevec{a}^{\halfComma i} =\spacevec{a}_j\times\spacevec{a}_k,\quad(i,j,k)\,\text{ cyclic.}
\end{equation}
Therefore, the contravariant basis can be computed from the covariant basis, which in turn can be computed directly from the transformation of the reference element to a physical element. 

Now it is possible to write the metric identities \eqref{eq:metricIDs_weirdForm} compactly in terms of the contravariant vectors,
\begin{equation}\label{eq:metricIDs_better}
\sum_{i=1}^3 \pderivative{J\spacevec{a}^{\halfComma i}}{\xi^i} = \vecNablaXi \cdot \left(J\spacevec{a}^{\halfComma i} \right)= 0.
\end{equation}

Since a contravariant vector points in the direction of the normal along an element face, it is easy to construct a normal. The (outward) pointing normal vectors in the physical coordinate in terms of the reference coordinates are
\begin{equation}\label{eq:normals}
\spacevec{n}_i  = \frac{J\spacevec{a}^{\halfComma i}}{|J\spacevec{a}^{\halfComma i}|}= \frac{|J|}{J}\frac{\spacevec{a}_j\times\spacevec{a}_k}{|\spacevec{a}_j\times\spacevec{a}_k|}.
\end{equation}

The transformation allows the normal in physical space to be related to the normal in reference space.
The reference space normal (but not normalized) vectors are directly written in the contravariant basis \eqref{eq:contra_vec_map}
\begin{equation}\label{eq:reference_normals}
\hat{n}^i = J\spacevec{a}^{\halfComma i},\quad i = 1,2,3.
\end{equation}
Going back to \eqref{eq:co_surf_elem} and \eqref{eq:contra_vec_map}, we can write the size of the surface differential elements as $\hat{s}$, 
\begin{equation}
\label{eq:surface_metrics}
\begin{aligned}
     \xi =\pm 1  \quad&:\quad            \hat{s}(\eta,\zeta) =\left|J\spacevec{a}^{1}(\pm1, \eta,\zeta)\right|,
\\  \eta =\pm 1  \quad&:\quad            \hat{s}(\xi,\zeta)  =\left|J\spacevec{a}^{2}( \xi, \pm1,\zeta)\right|,
\\ \zeta =\pm 1  \quad&:\quad            \hat{s}(\xi,\eta)   =\left|J\spacevec{a}^{3}( \xi, \eta, \pm1)\right| \,.
\end{aligned}
\end{equation}
The surface elements \eqref{eq:surface_metrics} are continuous across the element interface since the contravariant vectors \eqref{eq:contra_vec_map} 
and the covariant vectors \eqref{eq:covariant_vectors} are defined to be tangent to the shared face. 

Now we can relate the two normal vector representations, either in physical space \eqref{eq:normals} or reference space \eqref{eq:reference_normals}, through
\begin{equation}\label{eq:normal_vec_equiv}
\hat{n}^i = \frac{J\spacevec{a}^{\halfComma i}}{\hat{s}}\,\hat{s} = \spacevec{n}_i\,\hat{s},\quad i = 1,2,3,
\end{equation}
for the appropriate surface element $\hat{s}$ corresponding to a particular face. Since the $\hat s$ are continuous between elements sharing a face,
the normal vector only changes sign.

\CCLsubsubsection{Summary}

To extend approximations defined on a square or cube to a general quadrilateral or hexahedron, we re-write differential operators in physical coordinates in terms of reference space coordinates through the contravariant basis vectors \eqref{eq:contra_vec_map}. We summarize the common differential operators of divergence, gradient, and curl in Table \ref{tab:Mappings}.
{\begin{table}[!ht]
		\centering
		\caption{Differential operators in physical and computational coordinates.}
		\label{tab:Mappings}
\begingroup
\renewcommand{\arraystretch}{2.25} 
\aboverulesep=0ex
\belowrulesep=0ex
\begin{tabu} to \textwidth {@{}X[l]|X[r]@{}}
				\toprule
				Physical element & Reference element\\
				\midrule
				\midrule
				$\vecNablaX\cdot\spacevec{f}$ & $\displaystyle\frac{1}{J}\sum_{i=1}^3\pderivative{}{\xi^i}\left(J\spacevec{a}^{\halfComma i}\cdot\spacevec{f}\,\right)$ \\
				$\vecNablaX g$ & $\displaystyle\frac{1}{J}\sum_{i=1}^3 J{\spacevec {a}^{\halfComma i}} \pderivative{g}{\xi^i}$ \\
				$\vecNablaX\times\spacevec{f}$ & $\displaystyle\frac{1}{J}\sum_{i=1}^3\pderivative{}{\xi^i}\left(J\spacevec{a}^{\halfComma i}\times\spacevec{f}\,\right)$ \\
				\bottomrule
       	\end{tabu}
		\endgroup
	\end{table}}
	
The divergence operator can be written compactly by defining the \emph{volume weighted contravariant flux} vector $\contraspacevec{f}$, whose components are $\contravec{f}^{i}={{J}{\spacevec{a}^{\halfComma i}} \cdot \spacevec{ f}}$. In terms of the contravariant flux, the divergence looks similar in both physical and reference spaces, 
\begin{equation}
\vecNablaXi\cdot  \spacevec{f}
= \frac{1}{{J}}\sum\limits_{i = 1}^3 {\frac{\partial }{{\partial {\xi ^i}}}\left( {{J}{\spacevec{a}^{\halfComma i}} \cdot \spacevec{ f}} \right)}
= \frac{1}{{J}}\sum\limits_{i = 1}^3 {\frac{\partial \tilde f^{i}}{{\partial {\xi^i}}}}
= \frac{1}{{J}}\vecNablaXi \cdot \contraspacevec{ f}.
\label{eq:CompSpaceDivergence_DAK}
\end{equation}

\CCLsection{Building a Modern Discontinuous Galerkin Spectral Element Approximation}\label{sec:curvi_DG}
%
%




In this section, we use the derivations of Sec.~\ref{sec:equation_mapping} and apply a mapping from physical to reference space to the compressible Navier-Stokes equations written in mixed form \eqref{eq:navierstokes_mixed} and derive a DG spectral element approximation for that system.

To extend the transformation of the gradient and divergence operators of Table \ref{tab:Mappings}  to a system of partial differential equations, we define a matrix of the metric terms,
\begin{equation}\label{eq:metric_matrix}
\bigmatrix{M} = \begin{bmatrix}
Ja_1^1 \mmatrix{I} & Ja_1^2 \mmatrix{I} & Ja_1^3 \mmatrix{I}\\[0.05cm]
Ja_2^1 \mmatrix{I} & Ja_2^2 \mmatrix{I} & Ja_2^3 \mmatrix{I}\\[0.05cm]
Ja_3^1 \mmatrix{I} & Ja_3^2 \mmatrix{I} & Ja_3^3 \mmatrix{I}\\[0.05cm]
\end{bmatrix}
\end{equation}
with a $5\times 5$ identity matrix $\mmatrix{I}$, to match the size of the Navier-Stokes state variables. With \eqref{eq:metric_matrix}, the transformation of the gradient of a state vector is
\begin{equation}\label{eq:state_gradient}
\vecNablaX\statevec u
=
\begin{bmatrix}
  \statevec{u}_x \\ 
  \statevec{u}_y \\ 
  \statevec{u}_z 
\end{bmatrix}
=
 \frac{1}{J}\bigmatrix{M}
 \begin{bmatrix}
  \statevec{u}_{\xi} \\ 
  \statevec{u}_{\eta} \\ 
  \statevec{u}_{\zeta} 
\end{bmatrix}
= 
\frac{1}{J} \bigmatrix{M}\vecNablaXi\statevec{u}
\end{equation}
and the transformation of the divergence is
\begin{equation}\label{eq:state_divergence}
\vecNablaX  \cdot \bigstatevec{f} = \frac{1}{J}\vecNablaXi\cdot\left({\bigmatrix{M}^T}\bigstatevec{f}\right).
\end{equation} 
Moreover, the matrix \eqref{eq:metric_matrix} allows us to define contravariant block vectors
\begin{equation}\label{eq:block_contra}
\bigcontravec{\!f} = \bigmatrix{M}^T\bigstatevec{f}.
\end{equation}

Applying the differential operator transformations \eqref{eq:state_gradient} and \eqref{eq:state_divergence} as well as the contravariant block vector notation \eqref{eq:block_contra}, we get the transformed compressible Navier-Stokes equations
\begin{equation}\label{eq:mapped_NS}
\begin{aligned}
J\statevec{u}_t + \vecNablaXi\cdot\bigcontravec{\!f}(\statevec{u})&= \overRe \vecNablaXi \cdot \bigcontravec{\!f}^v\left( \statevec{u},\bigstatevec{q} \right),  \\[0.05cm]
J\bigstatevec{q} &= \bigmatrix{M}\vecNablaXi\statevec{w}\,.
\end{aligned}
\end{equation}
Note, to build an approximation that accounts for the entropy, we have taken the auxiliary variable $\bigstatevec{q}$ to be the gradient of the entropy variables \eqref{eq:entVariables}, to match the continuous equations \eqref{eq:weakformNS1}.

We first apply the polynomial ansatz from the DG toolbox described in Sec.~\ref{sec:InterpMultipleSpace_DAK} to approximate the solution, fluxes, entropy variables, etc. as interpolants written in the Lagrange basis \eqref{eq:Interpolant3D}. These quantities are denoted with capital letters, e.g. $\statevec{u}\approx\statevec{U}$. 

We then generate weak forms of the equations as in Sec.~\ref{sec:energy_entropy_continuous}: We multiply the transformed equations \eqref{eq:mapped_NS} by test functions $\testfuncOne$ and $\testfuncTwo$, also polynomials, and integrate over the reference element $E_{d}$. Any integrals in the weak formulation are approximated with the Gauss-Lobatto quadrature and the quadrature points are collocated with the interpolation points, as discussed in Sec.~\ref{sec:disc_IP_and_SBP}. Finally, we apply multidimensional summation-by-parts \eqref{eq:DiscreteGreens_DAK} to move derivatives off the fluxes and onto the test functions, generating boundary terms. The result is the set of two weak forms,
\begin{equation}\label{eq:DG_weak_form}
\begin{aligned}
\iprodN{\IN{J}\statevec{U}_t,\testfuncOne} &+ \isurfEN\testfuncOne^T\left\{\statevec{F}_n - \statevec{F}^v_n\right\}\hat{s}\dS - \iprodN{\bigcontravec{F},\vecNablaXi\testfuncOne} = -\overRe\iprodN{\bigcontravec{F}^v,\vecNablaXi\testfuncOne},\\[0.05cm]
\iprodN{\IN{J}\bigstatevec{Q},\testfuncTwo} &= \isurfEN \statevec{W}^T\left(\testfuncTwo\cdot\spacevec{n}\right)\hat{s}\dS - \iprodN{\statevec{W},\vecNablaXi\cdot\left(\bigmatrix{M}^T\testfuncTwo\right)}.
\end{aligned}
\end{equation}
Here, we use a compact notation for the normal fluxes, i.e., the normal flux in physical space $\statevec F_n=\left(\bigstatevec F \cdot \spacevec{n}\right)$.

A result of the DG polynomial ansatz is that solution values at element interfaces are discontinuous, and thus, the surface fluxes are not uniquely defined. This presents a problem to uniquely determine the normal fluxes, $\statevec F_n$. To resolve this, the elements are coupled through the boundaries as in a finite volume scheme with appropriate \textit{numerical flux functions} denoted by $\statevec F_n^{*}$, $\statevec F_n^{v,*}$ and $\statevec W^{*}$. 

The numerical fluxes are functions of two states, one to the left and one to the right of the interface, e.g. $\statevec F_n^{*}\left( \statevec U_{L}, \statevec U_{R}\right)$. They must also be \emph{consistent}, i.e. $\statevec F_n^{*}\left( \statevec U, \statevec U\right) = \bigstatevec f(\statevec U )\cdot\spacevec n$ so that if there is no jump, the exact flux in the normal direction is recovered. Other conditions, we will see, are still needed to ensure stability of the numerical scheme.

With the discontinuities at element interfaces resolved, we can perform another application of the multidimensional summation-by-parts \eqref{eq:DiscreteGreens_DAK} on the first equation in \eqref{eq:DG_weak_form} to move derivatives from the test functions back to the transformed flux vectors
\begin{equation}\label{eq:DG_strong_form}
\begin{aligned}
\iprodN{\IN{J}\statevec{U}_t,\testfuncOne} &+ \iprodN{\vecNablaXi\cdot\IN{\bigcontravec{F}},\testfuncOne} + \isurfEN\testfuncOne^T\left\{\statevec{F}^*_n - \statevec{F}_n\right\}\hat{s}\dS \\
&=\overRe \iprodN{\vecNablaXi\cdot\IN{\bigcontravec{F}^v},\testfuncOne} + \overRe\isurfEN\testfuncOne^T\left\{\statevec{F}^{v,*}_n - \statevec{F}^v_n\right\}\hat{s}\dS,\\[0.05cm]
\iprodN{\IN{J}\bigstatevec{Q},\testfuncTwo} &= \isurfEN \statevec{W}^{*,T}\left(\testfuncTwo\cdot\spacevec{n}\right)\hat{s}\dS - \iprodN{\statevec{W},\vecNablaXi\cdot\left(\bigmatrix{M}^T\testfuncTwo\right)}.
\end{aligned}
\end{equation}
The first equation of \eqref{eq:DG_strong_form} is called the \emph{strong form} of the DG approximation. It will be used later to create an entropy stable method. Note that the surface contributions within each element for the first equation resemble a penalty method in this form, in that it is proportional to the difference between the numerical flux and the flux computed from the interior. Note also that at this point, the multidimensional summation-by-parts,  \eqref{eq:DiscreteGreens_DAK}, says that \eqref{eq:DG_weak_form} and \eqref{eq:DG_strong_form} are algebraically equivalent. 

Approximations like \eqref{eq:DG_weak_form} and \eqref{eq:DG_strong_form} have been used in practice for many years, and they ``usually'' work. Sometimes, however, they are known to be unstable in that the computations blow up with unbounded energy or entropy. Since we have already shown in \eqref{eq:EnergyNormTimeDerivativeContinuous_{DAK}} and \eqref{eq:continuous_NSE_entropy_estimate} that the energy of the linear equations and the entropy of the nonlinear equations is bounded by the boundary contributions, we should expect the numerical schemes to share these properties: i.e. they should be \emph{stable}. 

The problem is that \textit{even for linear fluxes, the approximations \eqref{eq:DG_weak_form} and \eqref{eq:DG_strong_form} are not necessarily stable}. Based on what we have seen so far, we should require the numerical schemes to mimic the properties of the continuous equations. To that end, it should not be surprising that we should start with the same split form of the equation, \eqref{eq:SplitLinearSystem0}, that was used to show boundedness of the continuous solution.

Starting with \eqref{eq:SplitLinearSystem0}, we construct an alternative, \emph{split form} DG approximation, by approximating the divergence of the flux with an approximation of the average of the conservative and nonconservative forms,
\begin{equation}
\vecNablaXi\cdot\IN{\bigcontravec{F}} \approx\oneHalf \mathbb I^{N}\left\{ \vecNablaXi\cdot\IN{\bigcontravec{F}} +
\left(\bigcontravec A\right)^{T}{\vecNablaXi }\statevec U 
+ \left(\vecNablaXi\cdot\bigcontravec A\right)\statevec U \right\},
\label{eq:SplitFormApproxLin_{DAK}}
\end{equation}
where $\bigcontravec  A = \IN{\bigmatrix{M}^{T}\IN{\bigstatevec A}}$ and $\bigcontravec{F} = \IN{\bigcontravec  A\statevec U}$. Since we have already assumed that $\vecNabla\cdot\bigstatevec A=0$ to ensure that any energy growth is the system is due solely to boundary conditions, we will also assume in the following that $\vecNabla\cdot\bigcontravec  A=0$. Alternatively, we can simply drop that term from the approximation, since it will be a spectrally accurate approximation to zero. Doing so will lead to an approximation that is conservative only to within spectral accuracy, but that is less critical for linear systems of equations than for nonlinear.

With the split form approximation to the divergence, the DG approximation of advection terms of \eqref{eq:DG_strong_form} becomes
\begin{equation}
\begin{split}
 \oneHalf\iprodN{\vecNablaXi\cdot\IN{\bigcontravec{F}},\testfuncOne} 
+ \oneHalf\iprodN{\left({\bigcontravec A}\right)^{T}{\vecNablaXi }\statevec U ,\testfuncOne} 
+ \isurfEN\testfuncOne^T\left\{\statevec{F}^*_n - \statevec{F}_n\right\}\hat{s}\dS.
\end{split}
\label{eq:FirstDGLinearSplitForm}
\end{equation}

We can write \eqref{eq:FirstDGLinearSplitForm} in any one of many algebraically equivalent forms, and then use whichever is convenient for a given purpose.
For instance, we can apply the multidimensional summation-by-parts, \eqref{eq:DiscreteGreens_DAK}, to the first term and move the coefficient matrix in the second 
to get the algebraically equivalent form
\begin{equation}
\begin{split}
- \oneHalf\iprodN{\bigcontravec{F},\vecNablaXi\testfuncOne} 
+ \oneHalf\iprodN{{\vecNablaXi }\statevec U ,\bigcontravec{F}^{(T)}\left(\testfuncOne\right)} 
+ \isurfEN\testfuncOne^T\left\{\statevec{F}^*_n - \oneHalf\statevec{F}_n\right\}\hat{s}\dS,
\end{split}
\label{eq:LinearDirectlyStableForm}
\end{equation}
where
\begin{equation}
 \bigcontravec{F}^{\,(T)}\left(\testfuncOne \right) = \IN{\bigcontravec  A^{T}\testfuncOne} =\mathbb{I}^{N}\left[ {\begin{array}{*{20}{c}}
  {{\tilde{\mmatrix A}^{T}_1}\testfuncOne } \\ 
  {{\tilde{\mmatrix A}^{T}_2}\testfuncOne } \\ 
  {{\tilde{\mmatrix A}^{T}_3}\testfuncOne } 
\end{array}} \right],\,
\end{equation}
is the test function flux composed with the transpose of the coefficient matrices.

Continuing on, we can apply the multidimensional summation-by-parts rule to the second term of \eqref{eq:LinearDirectlyStableForm} to get another algebraically equivalent form
\begin{equation}
\begin{split}
- \oneHalf\iprodN{\bigcontravec{F},\vecNablaXi\testfuncOne} 
- \oneHalf\iprodN{\statevec U ,{\vecNablaXi }\cdot\bigcontravec{F}^{(T)}\left(\testfuncOne\right)} 
+ \isurfEN\testfuncOne^T\statevec{F}^*_n \hat{s}\dS.
\end{split}
\label{eq:LinearConservativeForm}
\end{equation}
In this form, all derivatives are on the test functions.

Finally, we can add and subtract the second term in \eqref{eq:FirstDGLinearSplitForm} and combine the difference to get
\begin{equation}
\begin{split}
 \iprodN{\vecNablaXi\cdot\IN{\bigcontravec{F}},\testfuncOne} 
 &+ \isurfEN\testfuncOne^T\left\{\statevec{F}^*_n - \statevec{F}_n\right\}\hat{s}\dS
\\&+ \oneHalf\iprodN{\IN{{\bigcontravec A}}^{T}{\vecNablaXi }\statevec U -\vecNablaXi\cdot\IN{\bigcontravec{F}},\testfuncOne}.
\end{split}
\label{eq:ThirdDGLinearSplitForm}
\end{equation}
We summarize the equivalent forms for the approximation of the divergence of the flux in Table \ref{tab:EquivLinearForms}. 

{\begin{table}[ht]
		\centering
		\caption{Equivalent DG approximations to the advective flux divergence.}
		\label{tab:EquivLinearForms}
\begingroup
\setlength{\tabcolsep}{6pt} 
\renewcommand{\arraystretch}{2.6} 
\aboverulesep=0ex
\belowrulesep=0ex
\begin{adjustbox}{max width=\textwidth}
			\begin{tabular}{@{}l|r@{}}
				\toprule
				Form & Approximation\\
				\midrule
				\midrule
				Strong [S] & $\displaystyle\oneHalf\iprodN{\vecNablaXi\cdot\IN{\bigcontravec{F}},\testfuncOne} 
+ \oneHalf\iprodN{\IN{{\bigcontravec A}}^{T}{\vecNablaXi }\statevec U ,\testfuncOne} 
+ \isurfEN\testfuncOne^T\left\{\statevec{F}^*_n - \statevec{F}_n\right\}\hat{s}\dS
$ \\
				Weak [W]& $ - \displaystyle\oneHalf\iprodN{\bigcontravec{F},\vecNablaXi\testfuncOne} 
- \oneHalf\iprodN{\statevec U ,{\vecNablaXi }\cdot\bigcontravec{F}^{(T)}\left(\testfuncOne\right)} 
+ \isurfEN\testfuncOne^T\statevec{F}^*_n \hat{s}\dS$\\
				Directly Stable [DS] & $-\displaystyle \oneHalf\iprodN{\bigcontravec{F},\vecNablaXi\testfuncOne} 
+ \oneHalf\iprodN{{\vecNablaXi }\statevec U ,\bigcontravec{F}^{(T)}\left(\testfuncOne\right)} 
+ \isurfEN\testfuncOne^T\left\{\statevec{F}^*_n - \oneHalf\statevec{F}_n\right\}\hat{s}\dS
$ \\
				Strong + Correction [SC]& 
				$\displaystyle \iprodN{\vecNablaXi\cdot\IN{\bigcontravec{F}},\testfuncOne} 
 + \isurfEN\testfuncOne^T\left\{\statevec{F}^*_n - \statevec{F}_n\right\}\hat{s}\dS
+ \oneHalf\iprodN{\IN{{\bigcontravec A}}^{T}{\vecNablaXi }\statevec U -\vecNablaXi\cdot\IN{\bigcontravec{F}},\testfuncOne}
$\\
				\bottomrule
		\end{tabular}
		\end{adjustbox}
		\endgroup
	\end{table}}

Since the approximations in Table \ref{tab:EquivLinearForms} are algebraically
equivalent, we can choose which one to use depending on what property of the equations we wish to study. Additionally, their equivalence can be exploited in practice as they 
can be reduced to the same implementation in code.

For instance, 
the form of the advective volume quadratures in the  [DS] form \eqref{eq:LinearDirectlyStableForm} look like those in the continuous form \eqref{eq:SplitLinearSystem}, which was used to show energy boundedness. On the other hand, if we set $\testfuncOne = 1$ for each component in the  [W] form \eqref{eq:LinearConservativeForm}, the first two terms vanish, leaving only the surface quadrature. The result implies that the approximation
is conservative, since the integral (the quadrature is exact when the test function is a polynomial of degree zero) over the volume of the divergence is equal to the integral of the flux over the surface.

Finally, written in the  [SC] form, we see that the split form \eqref{eq:FirstDGLinearSplitForm} is the original conservative form, \eqref{eq:DG_strong_form}, plus a \emph{correction} that is the discrete projection of  
\begin{equation}
\textrm{CN} = \IN{{\bigcontravec A}}^{T}{\vecNablaXi }\statevec U -\vecNablaXi\cdot\IN{\bigcontravec{F}}.
\end{equation}
The quantity $\textrm{CN}$ is the amount by which the product rule fails to hold due to aliasing when taking the divergence of the linear flux, $\bigcontravec{F} = \IN{\bigcontravec  A\statevec U}$. In other words, the split form approximation in \eqref{eq:SplitFormApproxLin_{DAK}} serves to cancel the product rule (aliasing) error in the divergence approximation.

Finally, since the problem is linear, we use the state rather than the entropy variables in the second equation of \eqref{eq:DG_strong_form} for the diffusion approximation
\begin{equation}
\iprodN{\IN{J}\bigstatevec{Q},\testfuncTwo} = \isurfEN \statevec{U}^{*,T}\left(\testfuncTwo\cdot\spacevec{n}\right)\hat{s}\dS - \iprodN{\statevec{U},\vecNablaXi\cdot\left(\bigmatrix{M}^T\testfuncTwo\right)},
\end{equation}
to match the second equation of \eqref{eq:SplitLinearSystem0}.

\CCLsubsection{Role of the Split Form Approximation}\label{sec:linear_split_form}
It is desirable that the approximation match as many properties of the original PDE as possible. Important properties include boundedness of the solution (stability),
conservation, free-stream preservation, phase and dissipation properties. In this section we will show that the split form approximation is stable, and if the metric terms are computed so that they satisfy the metric identities discretely, the approximation is free-stream (or constant state) preserving.
\CCLsubsubsection{Stability}
In this section we show that the discontinuous Galerkin approximation to the linear system of equations is stable if the split form approximation of the divergence is used. In the process, we will see precisely why the straight forward divergence approximation is not stable, but will often run stably.

Using the [DS] form of the advective terms from Table \ref{tab:EquivLinearForms}, the split form discontinuous Galerkin approximation of the linear Navier-Stokes equations is\begin{equation}\label{eq:DG_Split_formForStability}
\resizebox{\textwidth}{!}{$\displaystyle{
\begin{aligned}
\iprodN{\IN{J}\statevec{U}_t,\testfuncOne} &
- \oneHalf\iprodN{\bigcontravec{F},\vecNablaXi\testfuncOne} 
+ \oneHalf\iprodN{{\vecNablaXi }\statevec U ,\bigcontravec{F}^{(T)}\left(\testfuncOne\right)} 
+ \isurfEN\testfuncOne^T\left\{\statevec{F}^*_n - \oneHalf\statevec{F}_n\right\}\hat{s}\dS \\
&= \overRe\iprodN{\vecNablaXi\cdot\IN{\bigcontravec{F}^v},\testfuncOne} + \overRe\isurfEN\testfuncOne^T\left\{\statevec{F}^{v,*}_n - \statevec{F}^v_n\right\}\hat{s}\dS,
\\[0.05cm]
\iprodN{\IN{J}\bigstatevec{Q},\testfuncTwo} &= \isurfEN \statevec{U}^{*,T}\left(\testfuncTwo\cdot\spacevec{n}\right)\hat{s}\dS - \iprodN{\statevec{U},\vecNablaXi\cdot\left(\bigmatrix{M}^T\testfuncTwo\right)}.
\end{aligned}
}$}
\end{equation}

To assess stability, we follow the steps as to show energy boundedness in Sec. \ref{sec:energy_entropy_continuous}.  This time we first we set $\testfuncTwo =  \IN{\left(\bigmatrix S^{-1}\right)^{T}\bigmatrix S^{-1}\bigmatrix B\bigstatevec Q}$ in the second equation of \eqref{eq:DG_Split_formForStability}. 
Using the fact that $\left(\bigmatrix S^{-1}\right)^{T}$ commutes with $\bigmatrix M^{T}$,
\begin{equation}
\resizebox{\textwidth}{!}{$\displaystyle{
\begin{aligned}
\iprodN{  \IN{J}\bigmatrix S^{-1}\bigstatevec Q,\bigmatrix S^{-1}\bigmatrix B\bigstatevec Q}
&= 
\iprodN{ \IN{J}\bigstatevec Q^{s},\bigmatrix B^{s}\bigstatevec Q^{s}}
\\&= \isurfEN {\left\{{{\statevec U}^{s,*}-\statevec U^{s}}\right\}^{T}\statevec{F}^{v,s}_n \hat{s}\dS} 
 +\iprodN{\vecNablaXi{\statevec U^{s}, \bigcontravec F^{v,s}}},
\end{aligned}
}$}
\end{equation}
where $\statevec{F}^{v,s}_n = \bigstatevec{F}^{v,s} \cdot \spacevec{n}$ and $\bigcontravec F^{v,s}=\mathbb I^{N}\left(\bigmatrix M^{T}\bigmatrix B^{s}\bigstatevec Q^{s}\right)$. Therefore,
\begin{equation}
\iprodN{\vecNablaXi\statevec U^{s}, \bigcontravec F^{v,s}}
= \iprodN{ \IN{J}\bigstatevec Q^{s},\bigmatrix B^{s}\bigstatevec Q^{s}}
- \isurfEN {\left\{{{\statevec U}^{s,*}-\statevec U^{s}}\right\}^{T} \statevec{F}^{v,s}_n \hat{s}\dS }, 
\end{equation}
where 
\begin{equation}
\iprodN{ \IN{J}\bigstatevec Q^{s},\bigmatrix B^{s}\bigstatevec Q^{s}}\geqslant  0.
\label{eq:LinearDissibationBound}
\end{equation}

Next, 
we set $\testfuncOne = \left(\bigmatrix S^{-1}\right)^{T}\bigmatrix S^{-1}\statevec U = \left(\bigmatrix S^{-1}\right)^{T}\statevec U^{s}$ in the first equation of \eqref{eq:DG_Split_formForStability}.  The time derivative term becomes
\begin{equation}
\begin{split}
\iprodN{\IN{J}\statevec{U}_t,\left(\bigmatrix S^{-1}\right)^{T}\bigmatrix S^{-1}\statevec U}
&=\iprodN{ \IN{J}\statevec U^{s}_{t},\statevec U^{s}}
\\&=\oneHalf\frac{d}{dt}\sum_{ijk=0}^{N} \IN{J}_{ijk}\left|\statevec U^{s}\right|^{2}w_{ijk}.
\end{split}
\label{eq:TimeDerivInnerProduct}
\end{equation}
For \eqref{eq:TimeDerivInnerProduct} to represent a norm, equivalent to the continuous energy norm \eqref{eq:continuousTimeDepNorm}, and for \eqref{eq:LinearDissibationBound} to hold, it is necessary that
$\IN{J}_{ijk}>0$ for all $N$. This fact should be remembered in the grid generation process to ensure that the energy is always positive. If this is true, then we can write
\begin{equation}
\iprodN{ \IN{J}\statevec U^{s}_{t},\statevec U^{s}} \equiv \oneHalf\frac{d}{dt}\inorm{\statevec U^{s}}^{2}_{J,N}.
\end{equation}

The advective volume terms in \eqref{eq:DG_Split_formForStability} cancel when we substitute for the test function, as they did for the continuous terms, \eqref{eq:ExactFluxCancellation}, leaving us with
\begin{equation}
\begin{split}
\frac{1}{2}\frac{d}{dt} \inorm{\statevec U^{s}}^{2}_{J,N}  =  &-  \isurfEN {\left(\statevec U^{s}\right)^T\left( {\left\{\statevec{F}^{s,*}_n - \frac{1}{2}\statevec F^{s}_n\right\}  - \overRe\statevec F^{v,s,*}_n} \right) \hat{s}\dS}
\\&+\overRe\int\limits_{\partial E , N} {\left\{{{\statevec U}^{s,*}-\statevec U^{s}}\right\}^{T} \statevec F^{v,s}_n\hat{s}\dS}\\&- \overRe\iprodN{ J\bigstatevec Q^{s},\bigmatrix B^{s}\bigstatevec Q^{s}}\,.
\end{split}
\end{equation}
Separating the advective and viscous boundary terms, the elemental contribution to the total energy is
\begin{equation}
\begin{aligned}
\frac{1}{2}\frac{d}{dt} \inorm{\statevec U^{s}}^{2}_{J,N}=&- 
\isurfEN {{\left(\statevec U^{s}\right)^T\left\{{\statevec F}^{s,*}_n - \frac{1}{2}\statevec F^{s}_n\right\}} \hat{s}\dS}\\&
+  \overRe\int\limits_{\partial E,N } {\left\{\left(\statevec U^{s}\right)^T\statevec{F}_{n}^{v,s,*} +\left(\statevec U^{s,*}\right)^{T}\statevec F^{v,s}_n -\left(\statevec U^{s}\right)^T\statevec{F}^{v,s}_n\right\} \hat{s}\dS}
\\&- \overRe\iprodN{ J\bigstatevec Q^{s},\bigmatrix B^{s}\bigstatevec Q^{s}}.
\end{aligned}
\end{equation}

The total energy is found by summing over all of the elements. At the element faces, there will be jumps in the solution states and fluxes. To represent those jumps, we introduce the \emph{jump operator}: For a quantity $V$ defined on the left, $L$, and right, $R$, side of an interface with respect to the outward normal, 
\begin{equation}
\jump{V}\equiv V_{R}-V_{L}
\label{eq:jump_operator}
\end{equation}
 is the jump operator.
 
Summing over all elements,
\begin{equation}
\resizebox{\textwidth}{!}{$\displaystyle{
\begin{aligned}
\frac{1}{2}\frac{d}{dt} \sum\limits_{k=1}^{K}\inorm{\statevec U^{s,k}}^{2}_{J,N}=& 
 \sum\limits_{\interiorfaces} \int\limits_{N } {\left\{\jump{\statevec U^{s}}^T{\statevec F}_{n}^{s,*}-\frac{1}{2}\jump{\left(\statevec U^{s}\right)^{T}\bigstatevec F^{s}}\cdot\spacevec{n} \right\} \hat{s}\dS}\\&
- \overRe \sum\limits_{\interiorfaces} \int\limits_{N } {\left\{\jump{\statevec U^{s}}^T\statevec F_{n}^{v,s,*}+\left(\statevec U^{s,*}\right)^{T}\jump{\bigstatevec F^{v,s} }\cdot\spacevec{n}-\jump{\left(\statevec{U}^{s}\right)^T\bigstatevec F^{v,s}}\cdot\spacevec{n}\right\} \hat{s}\dS}
\\&-\overRe \sum\limits_{k=1}^{K}\iprodN{ J\bigstatevec Q^{s,k},\bigmatrix B^{s}\bigstatevec Q^{s,k}}
\\& + \PBT\,,
\end{aligned}
}$}
\label{eq:LinearEnergySumBound}
\end{equation}
where $\statevec U^{s,k}$ is the (symmetric) solution vector on element $k$. The quantity PBT represents the physical boundary terms, which we assume are dissipative, i.e. $\PBT\leqslant 0$.

Sufficient conditions for stability are those for which the right hand side of \eqref{eq:LinearEnergySumBound} is always non-positive. Since the third term is always non-positive because $\bigmatrix B^{s}>0$, sufficent conditions are that at each node on the interior element faces, the numerical values $\statevec F_{n}^{v,s,*}$, $\statevec{F}_{n}^{s,*}$ and $\statevec U^{s,*}$ are chosen so that 
\begin{equation}
\jump{\statevec U^{s}}^T{\statevec F}_{n}^{s,*}-\frac{1}{2}\jump{\left(\statevec U^{s}\right)^{T}\bigstatevec F^{s}}\cdot\spacevec{n} \leqslant 0,
\label{eq:AdvectiveConditionLinear}
\end{equation}
and
\begin{equation}
\jump{\statevec U^{s}}^T\statevec F_{n}^{v,s,*}+\left(\statevec U^{s,*}\right)^{T}\jump{\bigstatevec F^{v,s} }\cdot\spacevec{n}-\jump{\left(\statevec{U}^{s}\right)^T\bigstatevec F^{v,s}}\cdot\spacevec{n} \geqslant  0.
\label{eq:ViscousLinearInterfaceCondition}
\end{equation}
With such conditions satisfied, the norm of the approximate solution is bounded by the initial conditions,
\begin{equation}
\frac{d}{dt} \inorm{\statevec U^{s}}^{2}_{J,N}\equiv \frac{d}{dt} \sum\limits_{k=1}^{K}\inorm{\statevec U^{s,k}}^{2}_{J,N}\leqslant 0\quad\Rightarrow\quad \inorm{\statevec U^{s}(T)}_{J,N}\leqslant \inorm{\statevec U_{0}^{s}}_{J,N}.
\label{eq:discreteLinearStability}
\end{equation}

So what remains is to find suitable numerical fluxes such that \eqref{eq:AdvectiveConditionLinear} and \eqref{eq:ViscousLinearInterfaceCondition} hold.
Since the advective part of the equations is hyperbolic, the advective flux can be split according the wave directions relative to the normal direction
\begin{equation}
\bigstatevec F^{s}\cdot \spacevec{n} = \left(\left(\bigmatrix M^{T}\bigstatevec {A}^{\!s}\right)\cdot\hat n\right)\statevec{U}^{s}
=\left(\bigcontravec{A}^{\!s}\cdot \hat n\right)\statevec{U}^{s}\;\equiv\; \tilde {\mmatrix A}_n^{\!s}\statevec{U}^{s} =  \left(\tilde {\mmatrix A}_n^{\!s,+}+\tilde {\mmatrix A}_n^{\!s,-}\right)\statevec{U}^{s}\,,
\end{equation}
where
\begin{equation}
\tilde{\mmatrix A}_n^{\!s,\pm} = \frac{1}{2}\left(\tilde{\mmatrix A}_n^{\!s}\pm \left| \tilde{\mmatrix A}_n^{\!s}\right| \right)\,.
\label{eq:MatrixEVSplit}
\end{equation}
From that splitting, we can write the numerical advective flux choosing the left and right states according to the wave direction given by the sign of the eigenvalues as
\begin{equation}
\statevec F_{n}^{s,*}\left(\statevec U^{s}_{L}, \statevec U^{s}_{R}\right) =\tilde {\mmatrix A}_n^{\!s,+}\statevec{U}^{s}_{L}+\tilde {\mmatrix A}_n^{\!s,-}\statevec{U}^{s}_{R}.
\label{eq:upwinding}
\end{equation}
We substitute \eqref{eq:MatrixEVSplit} into \eqref{eq:upwinding} and rearrange to get a numerical flux
\begin{equation}
\begin{split}
\statevec F_{n}^{s,*}\left(\statevec U^{s}_{L}, \statevec U^{s}_{R}\right)&=\frac{\tilde{\mmatrix A}_n^{\!s}\statevec U^{s}_{L} + \tilde{\mmatrix A}_n^{\!s}\statevec U^{s}_{R}}{2}+\frac{\sigma}{2} \left| \tilde{\mmatrix A}_n^{\!s}\right|\left(\statevec U^{s}_{L} - \statevec U^{s}_{R}\right) \\[0.1cm]
&= \tilde{\mmatrix A}_n^{\!s}\avg{\statevec U^{s}} - \frac{\sigma}{2}\left| \tilde{\mmatrix A}_n^{\!s}\right|\jump{\statevec U^{s}}\,,
\end{split}
\label{eq:LinearNumericalFlux}
\end{equation}
where $\avg{\statevec U} = \oneHalf\left(\statevec U_{L}+ \statevec U_{R}\right)$ is the arithmetic mean.
For additional flexibility, we have added the parameter $\sigma$ so that the fully upwind numerical flux corresponds to $\sigma = 1$, whereas $\sigma = 0$ gives the central flux. 

With either the upwind or central numerical flux \eqref{eq:LinearNumericalFlux}, the contribution of the advective fluxes at the faces is dissipative. For any two state vectors,
\begin{equation}
\begin{split}
\jump{\statevec a^{T}\statevec b} &= \jump{\sum\limits_{m=1}^{5}a_{m}b_{m}} = \sum\limits_{m=1}^{5}\jump{a_{m}b_{m}} \\&= \sum\limits_{m=1}^{5}\left( \avg{a_{m}}\jump{b_{m}}+\jump{a_{m}}\avg{b_{m}}\right) = \avg{\statevec a}^{T}\jump{\statevec b} + \jump{\statevec a}^{T}\avg{\statevec b}\,.
\end{split}
\label{eq:KnightsFormerlyKnownAsLemma1}
\end{equation}
Therefore
\begin{equation}
\begin{split}
\jump{\left(\statevec U^{s}\right)^T\bigstatevec{F}^{s}} \cdot \spacevec{n} &=  \avg{\statevec U^{s}}^T\jump{\bigstatevec F^{s}} \cdot \spacevec{n}  + \jump{\statevec U^{s}}^T\avg{\bigstatevec F^{s}}\cdot \spacevec{n}\\[0.1cm]
&=\avg{\statevec U^{s}}^{T}\left(\bigstatevec{A}^{\!s} \cdot \spacevec{n} \right)\jump{\statevec U^{s}} + \jump{\statevec U^{s}}^{T}\left(\bigstatevec{A}^{\!s} \cdot \spacevec{n} \right)\avg{\statevec U^{s}} \\[0.1cm]
&=\avg{\statevec U^{s}}^{T}\tilde{\mmatrix A}_n^{\!s}\jump{\statevec U^{s}} + \jump{\statevec U^{s}}^{T}\tilde{\mmatrix A}_n^{\!s}\avg{\statevec U^{s}} \\[0.1cm]
&= 2\jump{\statevec U^{s}}^{T}\tilde{\mmatrix A}_n^{\!s}\avg{\statevec U^{s}}\,,
\end{split}
\end{equation}
so
\begin{equation}
\jump{\statevec U^{s}}^T\statevec F_{n}^{s,*} - \frac{1}{2}\jump{\left(\statevec U^{s}\right)^T\bigstatevec{F}^{s}} \cdot \spacevec{n} = -\frac{\sigma}{2}\jump{\statevec U^{s}}^{T}\left| \tilde{\mmatrix A}_n^{s}\right|\jump{\statevec U^{s}}\leqslant 0\,,
\label{eq:LinearInterfaceDissipation}
\end{equation}
which satisfies condition \eqref{eq:AdvectiveConditionLinear} for either the central numerical flux or an upwind flux,
and the contribution of the advective interface terms to the energy in \eqref{eq:LinearEnergySumBound} is nonpositive. 

We are now left to satisfy \eqref{eq:ViscousLinearInterfaceCondition} for the viscous terms. The simplest choice is to match the equality, which can be done with the \mbox{Bassi-Rebay-1} (or BR1 for short) numerical flux from \citet{Bassi&Rebay:1997:B&F97}, which computes the interface values as simple arithmetic means
\begin{equation}\label{eq:BR1_with_U}
 \begin{gathered}
  {\statevec{U}^{*}} = \frac{{{\statevec{U}_L} + {\statevec{U}_R}}}{2} = \avg{ \statevec{U}}  \hfill \\
  \statevec F_{n}^{v,s,*} = \left(\frac{\bigstatevec{F}_{L}^{v,s} + \bigstatevec{F}_{R}^{v,s}}{2}\right) \cdot \spacevec{n} = \avg{ \bigstatevec{F}^{v,s}}\cdot \spacevec{n}\,.  \hfill \\ 
 \end{gathered}
\end{equation}
When we make the substitution of the BR1 fluxes into the left side of \eqref{eq:ViscousLinearInterfaceCondition} it becomes (factoring out the normal direction)
\begin{equation}
\label{eq:linear_BR_step}
\left(\jump{\statevec U^{s}}^T \avg{\bigstatevec{F}^{v,s}}+\avg{ \statevec{U}}^{T}\jump{\bigstatevec F^{v,s} } - \jump{\left(\statevec{U}^{s}\right)^T\bigstatevec F^{v,s}}\right) \cdot\spacevec{n}.
\end{equation}
Then replacing the jump in the product using the identity \eqref{eq:KnightsFormerlyKnownAsLemma1}, the approximation satisfies \eqref{eq:ViscousLinearInterfaceCondition} because
\begin{equation}
\jump{\statevec U^{s}}^T\avg{\bigstatevec F^{v,s}} + \avg{\statevec U^{s}}^{T}\jump{\bigstatevec F^{v,s}}- \avg{\statevec U^{s}}^{T}\jump{\bigstatevec F^{v,s}} - \jump{\statevec U^{s}}^{T}\avg{\bigstatevec{F}^{v,s}} = \spacevec{0}\,,
\label{eqBR1DissipationLinear}
\end{equation}
so \eqref{eq:linear_BR_step} vanishes exactly. Therefore, the split form approximation \eqref{eq:DG_Split_formForStability} is stable in the sense of \eqref{eq:discreteLinearStability}.

We are now in the position to also see why the standard scheme, using only the divergence of the flux polynomial, can work, but is not guaranteed to be stable. In Table \ref{tab:EquivLinearForms} [SC] shows that the split form approximation of the advective terms is the standard scheme plus a correction term. Alternatively, the standard approximation is the split form minus that correction. 

If we subtract the correction term from the split form approximation and insert the results \eqref{eq:LinearInterfaceDissipation}, \eqref{eqBR1DissipationLinear} and \eqref{eq:LinearDissibationBound}, then the standard approximation satisfies
\begin{equation}
\begin{split}
\frac{1}{2}\frac{d}{dt} \inorm{\statevec U^{s}}^{2}_{J,N} = &-\overRe\inorm{\bigstatevec Q^{s}}_{\bigmatrix B^{s},N}^{2} - 
\frac{\sigma}{2} \sum\limits_{\interiorfaces} \int\limits_{N } { \jump{\statevec U^{s}}^{T}\left| \tilde{\mmatrix A}_n^{s}\right|\jump{\statevec U^{s}}\dS} \\&+
 \frac{1}{2}\sum\limits_{k=1}^{K}\left|{\iprodN {\left\{ {\IN{{\bigcontravec A}}^{T}  \vecNablaXi \statevec U^{s,k} -\vecNablaXi  \cdot \bigcontravec F^{s}(\statevec U^{s,k})} \right\},\statevec U^{s,k}} }\right|
 \\ &+ \PBT\,,
\end{split}
\label{eq:DGSEMBoundWithAliasing}
\end{equation}
which is \eqref{eq:LinearEnergySumBound} plus the correction term contribution.
The additional volume term is due to the failure of the product rule to hold for polynomial interpolants due to aliasing. 

Equation \eqref{eq:DGSEMBoundWithAliasing} shows that the physical diffusion and/or the dissipation associated with the numerical flux could counterbalance the product rule error and make the right hand side nonpositive. 
For well resolved solutions, the product rule error will be spectrally small, making it likely that the physical and interface dissipations are sufficiently large for stabilization. For under resolved solutions, the aliasing errors may be too large for the approximate solution to stay bounded. 
For large Reynolds numbers, the physical dissipation may be too small. The artificial dissipation due to the numerical fluxes might be sufficiently large, depending on flux solver is chosen. (For example, a Lax-Friedrichs numerical flux will be more dissipative than the exact upwind one.) 
Finally, changing from the BR1 to another viscous coupling procedure, coupled with a more dissipative numerical flux \emph{might} be enough to counteract the aliasing term. But, ultimately, the key to a stable discontinuous Galerkin spectral element method (DGSEM) is the stable approximation of the advective terms, as given in the split form approximation.

\CCLsubsection{The Importance of the Metric Identities}\label{sec:free_stream_AW}

One simple property of fluid flows and the solutions of the associated linearized equations with constant coefficient matrices is that a constant solution stays constant. This property is usually known as \emph{free-stream preservation} for fluid flows and we will use that term here. It is desirable that free-stream preservation holds for the approximate solution for if it doesn't, waves can spontaneously appear and propagate in an initially constant state even without applied external forces, see, e.g. \cite{Kopriva:2006er}. 

We now show that the split form spatial approximation of the constant coefficient linearized Euler equations is free-stream preserving provided that the approximations of metric terms satisfy a form of the metric identities, \eqref{eq:metricIDs_better}. Using the form [W] in Table \ref{tab:EquivLinearForms} for the advection terms, the DGSEM approximation of the advection equation is
\begin{equation}
\begin{split}
\iprodN{\IN{J}\statevec{U}_t,\testfuncOne} &+ \displaystyle\oneHalf\iprodN{\vecNablaXi\cdot\IN{\bigcontravec{F}},\testfuncOne} 
+ \oneHalf\iprodN{\left({\bigcontravec A}\right)^{T}{\vecNablaXi }\statevec U ,\testfuncOne} 
\\&+ \isurfEN\testfuncOne^T\left\{\statevec{F}^*_n - \statevec{F}_n\right\}\hat{s}\dS = 0
\end{split}
\label{eq:DGSEMForFSP}
\end{equation}
on each element. We can ignore the contribution of the diffusion terms since they are automatically zero when the gradients are zero.

When $\statevec U = \statevec C$ is constant over all elements, its gradient vanishes and the surface term in \eqref{eq:DGSEMForFSP} vanishes by consistency of the numerical flux. Therefore, if we write out the contravariant flux, 
\begin{equation}
\iprodN{\IN{J}\statevec{U}_t,\testfuncOne} = -\oneHalf\iprodN{\vecNablaXi\cdot\IN{\bigmatrix M \bigstatevec A\statevec C},\testfuncOne}
\end{equation}
Since $\bigstatevec A\statevec C$ is a constant, and $\testfuncOne$ is arbitrary, the right hand side of \eqref{eq:DGSEMForFSP} vanishes if and only if for each block of $\bigmatrix M$,
\begin{equation}
\vecNablaXi\cdot\IN{ J\spacevec a^{\halfComma i}}=0 \quad i = 1,2,3
\label{eq:discreteMetricIdentities}
\end{equation}
If we compare \eqref{eq:discreteMetricIdentities} with the metric identities \eqref{eq:metricIDs_better}, we see that the interpolant of the volume weighted contravariant basis vectors must vanish for the approximation to be free-stream preserving.
Since differentiation and interpolation do not commute, it is not immediately true that if the metric terms analytically satisfy the metric identities then their interpolants do as well. 

It is relatively straightforward to satisfy the metric identities in two spatial dimensions if the boundaries of the elements are polynomials. For such domains,
\begin{equation}
\begin{aligned}
J\spacevec a^{1} &=Y_{\xi}\hat x - X_{\eta}\hat y, \\
J\spacevec a^{2} &= -Y_{\xi}\hat x + X_{\xi}\hat y.
\end{aligned}
\end{equation}
Therefore, if the mapping $\spacevec X\in\PN{N}$, which it is if the boundary curves are isoparametric (polynomials of degree $N$) or less, $\IN{ J\spacevec a^{\halfComma i}}=J\spacevec a^{\halfComma i}$, and so \eqref{eq:discreteMetricIdentities} holds.

It is more complicated to satisfy the metric identities for general hexahedral elements in three spatial dimensions. Direct approximation of the cross product form of the metric terms \eqref{eq:contra_vec_map} will not satisfy the discrete metric identities except in special cases because
\begin{equation}
\sum\limits_{i = 1}^3 {\frac{\partial }{{\partial {\xi ^i}}}\left( \IN{J{\spacevec{a}^{\halfComma i}}} \right)}  = \sum\limits_{i = 1}^3 {\frac{\partial }{{\partial {\xi ^i}}}\left( \IN{\frac{{\partial \spacevec X}}{{\partial {\xi ^j}}} \times \frac{{\partial \spacevec X}}{{\partial {\xi ^k}}}} \right)} .
\end{equation}
Even if $\spacevec X \in\PN{N}$, the cross product is a polynomial of degree $2N$. Thus, aliasing errors will not allow the outer differentiation to commute with the interpolation operator to allow the terms to cancel. 

Special cases for which the cross product form \eqref{eq:contra_vec_map} can be used, then, are those where $\IN{\frac{{\partial \spacevec X}}{{\partial {\xi ^j}}} \times \frac{{\partial \spacevec X}}{{\partial {\xi ^k}}}} = \frac{{\partial \spacevec X}}{{\partial {\xi ^j}}} \times \frac{{\partial \spacevec X}}{{\partial {\xi ^k}}}$.
Such special cases include
\begin{itemize}
\item $\spacevec X\in\mathbb{P}^{N/2}$. If the faces of the hexahedral elements are approximated by half the order of the solution, then the product is a polynomial of degree $N$ and the interpolation is exact.
\item The element faces are planar and $N\geqslant 2$. A special case of item 1, the cross product form can be used if the faces are flat.
\end{itemize}

To avoid such limitations, a general formulation of the metric terms is necessary that satisfies the metric identities. This is achieved by writing the contravariant vector components in a curl form \citet{Kopriva:2006er}, for instance
\begin{equation}
\label{eq:curlMetrics}
Ja_n^i = - \hat{x}_i \cdot \vecNablaXi\times \left( \IN{{X_l}{\vecNablaXi }{X_m}} \right),\quad i = 1,2,3,\; n = 1,2,3,\;\;(n,m,l)\;\text{ cyclic}.
\end{equation}
Computed this way, $\IN{J\spacevec{a}^{\halfComma i}} = J\spacevec{a}^{\halfComma i}$ and the divergence of the curl is explicitly zero without the need to commute interpolation and differentiation.
Written out in full \eqref{eq:curlMetrics} reads
\begin{equation}
\resizebox{\textwidth}{!}{$
\begin{gathered}
  J{{\spacevec a}^1} = \left[ {{{\left( \IN{{Y_\eta }Z}\right)}_{\!\zeta} } - {{\left( \IN{{Y_\zeta }Z}\right)}_{\!\eta} }} \right]\hat x + \left[ {{{\left( \IN{{Z_\eta }X}\right)}_{\!\zeta} } - {{\left(\IN{{Z_\zeta }X} \right)}_{\!\eta} }} \right]\hat y + \left[ {{{\left( \IN{{X_\eta }Y}\right)}_{\!\zeta} } - {{\left( \IN{{X_\zeta }Y}\right)}_{\!\eta} }} \right]\hat z, \hfill \\
  J{{\spacevec a}^2} = \left[ {{{\left( \IN{{Y_\zeta }Z}\right)}_{\!\xi} } - {{\left(\IN{{Y_\xi }Z}\right)}_{\!\zeta} }} \right]\hat x + \left[ {{{\left( \IN{{Z_\zeta }X} \right)}_{\!\xi} } - {{\left( \IN{{Z_\xi }X} \right)}_{\!\zeta} }} \right]\hat y + \left[ {{{\left( \IN{{X_\zeta }Y} \right)}_{\!\xi} } - {{\left( \IN{{X_\xi }Y}\right)}_{\!\zeta} }} \right]\hat z, \hfill \\
  J{{\spacevec a}^3} = \left[ {{{\left( \IN{{Y_\xi }Z} \right)}_{\!\eta} } - {{\left( \IN{{Y_\eta }Z} \right)}_{\!\xi} }} \right]\hat x + \left[ {{{\left(\IN{{Z_\xi }X} \right)}_{\!\eta} } - {{\left( \IN{{Z_\eta }X} \right)}_{\!\xi} }} \right]\hat y + \left[ {{{\left( \IN{{X_\xi }Y} \right)}_{\!\eta} } - {{\left( \IN{{X_\eta }Y} \right)}_{\!\xi} }} \right]\hat z. \hfill \\ 
\end{gathered}
$}
\end{equation}

\CCLsubsection{The Concept of Flux Differencing and Two-Point Fluxes}\label{sec:yet_another_split_AW}

Although it may not appear so, one feature of the split form approximation is that it can be implemented by a simple modification of the volume integral of a standard DGSEM approximation, thereby taking a code that usually works and transforming it into a code that is provably stable. 

To get the implementation form, we use the form $[S]$ and take a tensor product of the Lagrangian basis functions to be the test functions, i.e. $\testfuncOne = {\ell _i}{\ell _j}{\ell _k}\,\statevec e_p$, where $p$ indicates a component of the state vector and $\statevec e_p$ the corresponding unit vector. Since for any state-vector polynomial $\statevec V\in\PN{N}$,
\begin{equation}{\iprodN{\statevec V,{\ell _i}{\ell _j}{\ell _k}\,\statevec e_p} } = \sum\limits_{n,m,l = 0}^N {{ V^p_{nml}}{\ell _i}( {{\xi _n}} ){\ell _j}( {{\eta _m}} ){\ell _k}( {{\zeta _l}} ){w_{nml}}}  = {V^p_{ijk}}{w_{ijk}}.
\end{equation}
Choosing the test functions in this way for all state components gives for the the volume term with the temporal derivative in $[S]$ 
\begin{equation}\iprodN{ {J{\statevec U_t},\testfuncOne } } \to {J_{ijk}}{{\dot {\statevec U}}_{ijk}}{w_{ijk}}.\end{equation}


Similarly,
\begin{equation}
{\iprodN{ {\vecNablaXi  \cdot \bigcontravec{F}\left( \statevec U \right),\testfuncOne }}} \to {w_{ijk}}\left\{ {\sum\limits_{n = 0}^N {{\contravec {\statevec F}^{1}_{njk}}{\dmat_{in}}}  + \sum\limits_{n = 0}^N {{\contravec {\statevec F}^{2}_{ink}}{\dmat_{jn}}}  + \sum\limits_{n = 0}^N {{\contravec {\statevec F}^{3}_{ijn}}{\dmat_{kn}}} } \right\},
\label{eq:DivFTerm}
\end{equation}
and
\begin{equation}
\resizebox{\textwidth}{!}{$\displaystyle{
{\iprodN{ {\IN{{\bigcontravec A}}^{T}\vecNablaXi \statevec U,\testfuncOne } }} \to {w_{ijk}}\left\{ {\tilde {\mmatrix A}_{ijk}^{1}\sum\limits_{n = 0}^N {{\statevec U_{njk}}{\dmat_{in}}}  + \tilde {\mmatrix A}_{ijk}^{2}\sum\limits_{n = 0}^N {{\statevec U_{ink}}{\dmat_{jn}}}  + \tilde {\mmatrix A}_{ijk}^{3}\sum\limits_{n = 0}^N {{\statevec U_{ijn}}{\dmat_{kn}}} } \right\},
\label{eq:ADotGradUTerm}
}$}
\end{equation}
where the $\tilde {\mmatrix A}^{i}=J\spacevec a^{i}\cdot\spacevec {\mmatrix A}$ are the contravariant coefficient matrix components.

If we add the vanishing terms of the divergence of the coefficient matrices,
\begin{equation}
\resizebox{\textwidth}{!}{$\displaystyle{
{\left( {\vecNablaXi  \cdot \IN {\bigcontravec A}\statevec U,\testfuncOne } \right)_N} \to {w_{ijk}}\left\{ {\sum\limits_{n = 0}^N {\tilde {\mmatrix A}_{njk}^{1}{\dmat_{in}}}  + \sum\limits_{n = 0}^N {\tilde {\mmatrix A}_{ink}^{2}{\dmat_{jn}}}  + \sum\limits_{n = 0}^N {\tilde {\mmatrix A}_{ijn}^{3}{\dmat_{kn}}} } \right\}{\statevec U_{ijk}},
\label{eq:DivATerm}
}$}
\end{equation}
we can gather the three terms, \eqref{eq:DivFTerm}, \eqref{eq:ADotGradUTerm} and \eqref{eq:DivATerm},  to see that
\begin{equation}
\begin{split}
{\iprodN{{{\vecNablaXi  } \cdot \tilde {\statevec{F}}\left( {\mathbf{U}} \right),\testfuncOne }}} 
&+ {\iprodN{ {\IN {\bigcontravec A}  \cdot {\vecNablaXi  }{\mathbf{U}},\testfuncOne } }} 
+ {\iprodN{ {{\vecNablaXi  } \cdot \IN{\bigcontravec A}{\mathbf{U}},\testfuncOne })}}
\to 
\\&
w_{ijk}\sum\limits_{n = 0}^N {\left\{ {{{\tilde {\statevec F}}^{1}_{njk}} + \tilde {\mmatrix A}_{ijk}^{1}{\statevec U_{njk}} + \tilde {\mmatrix A}_{njk}^{1}{\statevec U_{ijk}}} \right\}{\dmat_{in}}}  
\\&+ w_{ijk}\sum\limits_{n = 0}^N {\left\{ {{{\tilde {\statevec F}}^{2}_{ink}} + \tilde {\mmatrix A}_{ijk}^{2}{\statevec U_{ink}} + \tilde {\mmatrix A}_{ink}^{2}{\statevec U_{ijk}}} \right\}{\dmat_{jn}}} 
 \\&
+ w_{ijk}\sum\limits_{n = 0}^N {\left\{ {{{\tilde {\statevec F}}^{3}_{ijn}} + \tilde {\mmatrix A}_{ijk}^{3}{\statevec U_{ijn}} + \tilde {\mmatrix A}_{ijn}^{3}{\statevec U_{ijk}}} \right\}{\dmat_{kn}}}.
\end{split}
\label{eq:TwoPointFluxStart}
 \end{equation}
 
 The quantities in the braces in \eqref{eq:TwoPointFluxStart} can be interpreted as two point fluxes. For example, the quantity in braces in the first sum depends on the points $ijk$ and $njk,\;n=0,\ldots,N$. So let us define the two-point fluxes
 \begin{equation}
 \begin{aligned}
 \overline{\statevec F}^{1}_{(n,i)jk} &= {{{\tilde {\statevec F}}^{1}_{njk}} + \tilde {\mmatrix A}_{ijk}^{1}{\statevec U_{njk}} + \tilde {\mmatrix A}_{njk}^{1}{\statevec U_{ijk}}},\\[0.1cm]
 \overline{\statevec F}^{2}_{i(n,j)k} &= {{{\tilde {\statevec F}}^{2}_{ink}} + \tilde {\mmatrix A}_{ijk}^{2}{\statevec U_{ink}} + \tilde {\mmatrix A}_{ink}^{2}{\statevec U_{ijk}}},\\[0.1cm]
 \overline{\statevec F}^{3}_{ij(n,k)} &= {{{\tilde {\statevec F}}^{3}_{ijn}} + \tilde {\mmatrix A}_{ijk}^{3}{\statevec U_{ijn}} + \tilde {\mmatrix A}_{ijn}^{3}{\statevec U_{ijk}}}.
\end{aligned}
 \label{eq:modifiedLinearFlux}
 \end{equation}
 
 With two-point fluxes the advective part of the split form approximation \eqref{eq:DG_Split_formForStability} looks like the DSGEM implementation presented by \citet{Kopriva:2009nx} except that the fluxes in the derivative sums have been replaced
  \begin{equation}
\resizebox{\textwidth}{!}{$
\begin{aligned}
\dot{\statevec U}_{ijk} + \frac{1}{J_{ijk}}&
\left\{ \left[ {\left\{\contravec{\statevec{F}}^{*}_{Njk}-{{{{\contravec {\statevec F}}_{Njk}} \cdot \hat \xi}}\right\}\frac{\delta_{iN}}{{{w_i}}} 
- \left\{\contravec{\statevec{F}}^{*}_{0jk}-{{{{\contravec {\statevec F}}_{0jk}} \cdot \hat \xi}}\right\}\frac{\delta_{i0}}{{{w_i}}} 
+ \oneHalf\sum\limits_{n = 0}^N {{\overline{\statevec F}^{1}_{(n,i)jk}}{{ \dmat}_{in}}} } \right] \right.
 \\&+ \left[ {\left\{\contravec{\statevec{F}}^{*}_{iNk}-{{{{\contravec {\statevec F}}_{iNk}} \cdot \hat \eta}}\right\}\frac{\delta_{jN}}{{{w_j}}} 
- \left\{\contravec{\statevec{F}}^{*}_{i0k}-{{{{\contravec {\statevec F}}_{i0k}} \cdot \hat \eta}}\right\}\frac{\delta_{j0}}{{{w_j}}} 
+ \oneHalf\sum\limits_{n = 0}^N {{\overline{\statevec F}^{2}_{i(n,j)k}}{{ \dmat}_{jn}}} } \right] 
\\& + \left.\left[ {\left\{\contravec{\statevec{F}}^{*}_{ijN}-{{{{\contravec {\statevec F}}_{ijN}} \cdot \hat \zeta}}\right\}\frac{\delta_{kN}}{{{w_k}}} 
- \left\{\contravec{\statevec{F}}^{*}_{ij0}-{{{{\contravec {\statevec F}}_{ij0}} \cdot \hat \zeta}}\right\}\frac{\delta_{0k}}{{{w_k}}} 
+ \oneHalf\sum\limits_{n = 0}^N {{\overline{\statevec F}^{3}_{ij(n,k)}}{{ \dmat}_{kn}}} } \right] \right\}\\&=0.
 \end{aligned}
$}
 \label{eq:PointwiseWithModifiedFlux}
 \end{equation}

We can go further and re-write each of these two point fluxes in terms of two point averages. To do so, we add the derivative of a constant, which following \eqref{eq:sbp_property1_GG} is zero,
 \begin{equation}
 0 = \tilde {\mmatrix A}_{ijk}^{1}{\statevec U_{ijk}}\sum\limits_{n = 0}^N {{\dmat_{in}}} = \sum\limits_{n = 0}^N {\tilde {\mmatrix A}_{ijk}^{1}{\statevec U_{ijk}}{\dmat_{in}}}.
 \label{eq:DerivOfConstIsZero}
 \end{equation}
 Adding \eqref{eq:DerivOfConstIsZero} to \eqref{eq:modifiedLinearFlux}, we get a new two point flux, the first component of which we define by
 \begin{equation}
 \begin{split}
4{{\contravec{\statevec F}}^{\#,1}_{(n,i)jk}} 
 &= \tilde {\mmatrix A}_{njk}^{1}{\statevec U_{njk}} + \tilde {\mmatrix A}_{ijk}^{1}{\statevec U_{ijk}} + \tilde {\mmatrix A}_{ijk}^{1}{\statevec U_{njk}} + \tilde {\mmatrix A}_{njk}^{1}{\statevec U_{ijk}} ,
 \end{split}
 \label{eq:FirstLinearFSharpDef}
 \end{equation}
 and similarly for the other contravariant fluxes.
The right hand side of \eqref{eq:FirstLinearFSharpDef}  can be split into a product of two factors,
\begin{equation}
{{\contravec{\statevec F}}^{\#,1}_{(n,i)jk}}  = \frac{{\left( {\tilde {\mmatrix A}_{njk}^{1} + \tilde {\mmatrix A}_{ijk}^{1}} \right)}}{2}\frac{{\left( {{\statevec U_{ijk}} + {\statevec U_{njk}}} \right)}}{2} ,
\end{equation}
so that
\begin{equation}
 {{\contravec{\statevec F}}^{\#,1}_{(n,i)jk}}= {\avg{\tilde {\mmatrix A}^{1}}_{(n,i)jk}}{\avg{\statevec U}_{(n,i)jk}}.
 \label{eq:split_form_1}
\end{equation}
Since the linear flux $\bigcontravec f = \bigcontravec A \statevec U$, we see that the two-point flux whose divergence is equal to the split form approximation to the divergence of the flux can be expressed as the product of two averages.

With the definition of the two-point flux, we can re-write the sums in \eqref{eq:PointwiseWithModifiedFlux} with $\bigcontravec F^{\#}$, for example
\begin{equation}
\oneHalf\sum\limits_{n = 0}^N {{\bar{\statevec F}^{1}_{(n,i)jk}}{{ \dmat}_{in}}} = \sum\limits_{n = 0}^N 2{{\contravec{\statevec F}^{\#,1}_{(n,i)jk}}{{ \dmat}_{in}}},
\label{eq:split_form_2}
\end{equation}
to give the approximation at each point (and what one would code)
  \begin{equation}
\resizebox{\textwidth}{!}{$
\begin{aligned}
\dot{\statevec U}_{ijk} + \frac{1}{J_{ijk}}&
\left\{ \left[ {\left\{\contravec{\statevec{F}}^{*}_{Njk}-{{{{\contravec {\statevec F}}_{Njk}} \cdot \hat \xi}}\right\}\frac{\delta_{iN}}{{{w_i}}} 
- \left\{\contravec{\statevec{F}}^{*}_{0jk}-{{{{\contravec {\statevec F}}_{0jk}} \cdot \hat \xi}}\right\}\frac{\delta_{i0}}{{{w_i}}} 
+ \sum\limits_{n = 0}^N {2{\contravec{\statevec F}^{\#,1}_{(n,i)jk}}{{ \dmat}_{in}}} } \right] \right.
 \\&+ \left[ {\left\{\contravec{\statevec{F}}^{*}_{iNk}-{{{{\contravec {\statevec F}}_{iNk}} \cdot \hat \eta}}\right\}\frac{\delta_{jN}}{{{w_j}}} 
- \left\{\contravec{\statevec{F}}^{*}_{i0k}-{{{{\contravec {\statevec F}}_{i0k}} \cdot \hat \eta}}\right\}\frac{\delta_{j0}}{{{w_j}}} 
+ \sum\limits_{n = 0}^N {2{\contravec{\statevec F}^{\#,2}_{i(n,j)k}}{{ \dmat}_{jn}}} } \right] 
\\& + \left.\left[ {\left\{\contravec{\statevec{F}}^{*}_{ijN}-{{{{\contravec {\statevec F}}_{ijN}} \cdot \hat \zeta}}\right\}\frac{\delta_{kN}}{{{w_k}}} 
- \left\{\contravec{\statevec{F}}^{*}_{ij0}-{{{{\contravec {\statevec F}}_{ij0}} \cdot \hat \zeta}}\right\}\frac{\delta_{0k}}{{{w_k}}} 
+ \sum\limits_{n = 0}^N {2{\contravec{\statevec F}^{\#,3}_{ij(n,k)}}{{ \dmat}_{kn}}} } \right] \right\}\\&=0.
 \end{aligned}
 $}
 \label{eq:PointwiseLinearWithSharpFlux}
 \end{equation}

As a shorthand, we write the divergence operator implied by the three summations in \eqref{eq:PointwiseLinearWithSharpFlux} as
\begin{equation}
\begin{split}
\spacevec{\mathbb{D}} \cdot (\bigcontravec{F} )^{\#}(\xi,\eta,\zeta)\equiv2\sum_{n=0}^N 
&\quad \ell'_n(\xi)   \contrastatevec{F}^{\#,1}(\xi,\eta,\zeta;\xi_n,\eta,\zeta)\\[-1ex]       
&+     \ell'_n(\eta)  \contrastatevec{F}^{\#,2}(\xi,\eta,\zeta;\xi,\eta_n,\zeta)\\[1ex]    
&+     \ell'_n(\zeta) \contrastatevec{F}^{\#,3}(\xi,\eta,\zeta;\xi,\eta,\zeta_n),                                                                                                                                              
\end{split}
\label{eq:DFLinSharp}
\end{equation}
which allows us to add one more equivalent form for the divergence approximation to Table \ref{tab:EquivLinearForms}, namely
\begin{equation}
\textrm{Two-Point [T]}: \quad \displaystyle\iprodN{\spacevec{\mathbb{D}} \cdot(\bigcontravec{F} )^{\#},\testfuncOne} 
+ \isurfEN\testfuncOne^T\left\{\statevec{F}^*_n - \statevec{F}_n\right\}\hat{s}\dS.
\end{equation}
In summary, the split form approximation to the divergence, which includes three terms for $\vecNablaXi\cdot\bigcontravec F$, $\left(\vecNablaXi\cdot\bigcontravec A\right)\statevec U$ and $\bigcontravec A^{T}\vecNablaXi\statevec U$, can be re-written, represented, and coded in terms of a single two-point flux, $\bigcontravec F^{\#}$. Since it is algebraically equivalent, the approximation written this way is stable, and free-stream preserving if the metric terms are divergence free.


\CCLsection{The Final Assembly: A Robust DGSEM}\label{sec:section_splitformDG_GG}


This section serves as the culmination of this book chapter, where we present the details needed to construct an entropy stable DG approximation for the compressible Navier-Stokes equations. In the end, the solution of the numerical method will possess a discrete entropy bound which directly mimics that from the continuous analysis \eqref{eq:Entropy_stab_general_statment}. Great care was taken in the previous sections to discuss, contextualise and analyse the components of the DG approximation piecemeal, e.g. high-order accuracy, aliasing, free-stream preservation, split forms, etc., so that we are now fully equipped with a powerful spectral DG toolbox to address this final task.

Systematically, the split form DG approximation for nonlinear PDEs is outlined as follows:
\begin{enumerate}
\item The formulation for the linear problem in Sec. \ref{sec:yet_another_split_AW} is generalised. Herein it is highlighted that many components are similar, but the approximation of the advective terms undergo a fundamental change in structure.
\item Advective terms are given a primary focus because their proper treatment is critical to demonstrate entropy stability.
\item The stage is then set to present the discrete stability statement.
\end{enumerate}


We create the split form DG approximation from the strong form DG formulation \eqref{eq:DG_strong_form}, where the auxiliary variable for the approximation to the solution gradient, $\bigstatevec{Q}$, is written in terms of the entropy variables $\statevec{W}$. Just as in the previous section, we must address how to treat the contributions of the nonlinear advective and viscous fluxes in the volume of an element (quantities containing a divergence operator) as well as along its surface (quantities denoted with a star).

We approximate the viscous flux contributions first, because their treatment in the high-order split DG approximation for nonlinear problems is straightforward and utilises well-developed, \emph{standard} components of the DG toolbox, e.g. \citep{Hindenlang201286}, shown in Table \ref{tab:CalcComparison_DAK}. The divergence of the viscous fluxes in the volume are approximated by
\begin{equation}\label{eq:viscVolume}
\resizebox{\textwidth}{!}{$\displaystyle{
\begin{aligned}
\vecNablaXi \cdot \IN{\bigcontravec F^v}\approx \DDs\bigcontravec{F}^{v} = \sum_{m=0}^N \dmat_{im}\left(\contrastatevec{F}^{v}_1\right)_{\!mjk} + \sum_{m=0}^N \dmat_{jm}\left(\contrastatevec{F}^{v}_2\right)_{\!imk} + \sum_{m=0}^N \dmat_{km}\left(\contrastatevec{F}^{v}_3\right)_{\!ijm},
\end{aligned}
}$}
\end{equation}
where the metric terms are \emph{included} in the transformed viscous fluxes $\contrastatevec{F}^v_{l}$, $l=1,2,3$. 

Analogous to the linear approximation, we approximate the surface contribution of the viscous fluxes with the BR1 numerical flux
\begin{equation}\label{eq:BR1Fluxes_with_W}
\statevec F_n^{v,*} = \avg{\bigstatevec{F}^v}\cdot\spacevec{n},\quad \statevec{W}^*=\avg{\statevec{W}},
\end{equation}
where, again, the compact notation is used for the arithmetic mean. The only difference in this treatment of the viscous fluxes is the use of the discrete entropy variables and gradients in the auxiliary variable, in contrast to \eqref{eq:BR1_with_U} where the solution quantity $\statevec{U}$ was used. The BR1 terms are neutrally stable for the split form DG approximation of the nonlinear problem, as we show later in Sec. \ref{sec:es_proof_AW}.

Our focus now turns to the advective components in the approximation, which require greater care to produce an entropy stable, split form DG method. From a physical perspective, it makes sense that the advective terms tend to be more troublesome compared to the ``nice'' viscous terms. The split formulation fundamentally changes the structure of the flux divergence of the advective flux components.

Proper treatment of the volume contribution of the advective flux divergence needed to produce an entropy stable approximation is built from works in the finite difference community \citep{fisher2012,fisher2013,skew_sbp2,lefloch_rhode_2000} and the DG community \citep{carpenter_esdg,Gassner:2016ye,Gassner_BR1}. With notation introduced by \citet{Gassner_BR1}, we define the split form DG divergence approximation as
\begin{equation}
\label{eq:entropy-cons_volint}
\begin{split}
\vecNablaXi \cdot \IN{\bigcontravec F}\approx \DD\bigcontravec F^{\#} = 
&\quad 2\sum_{m=0}^N \dmat_{im}\left(\bigstatevec{F}^{\#}(\statevec{U}_{ijk}, \statevec{U}_{mjk})\cdot\avg{J\spacevec{a}^{\,1}}_{(i,m)jk}\right)\\ 
&+     2\sum_{m=0}^N \dmat_{jm}\left(\bigstatevec{F}^{\#}(\statevec{U}_{ijk}, \statevec{U}_{imk})\cdot\avg{J\spacevec{a}^{\,2}}_{i(j,m)k}\right)\\ 
&+     2\sum_{m=0}^N \dmat_{km}\left(\bigstatevec{F}^{\#}(\statevec{U}_{ijk}, \statevec{U}_{ijm})\cdot\avg{J\spacevec{a}^{\,3}}_{ij(k,m)}\right)\,
\end{split}
\end{equation}
for each Gauss-Lobatto node $i,j,k$ of an element. As with the approximation of the linear equations, \eqref{eq:DFLinSharp}, we have introduced $\bigstatevec{F}^{\#}$, an additional \emph{two-point volume flux} that is symmetric with respect to its arguments and consistent. We write the arithmetic mean in each spatial direction, defined as in \eqref{eq:split_form_1} compactly, e.g., in the $\xi-$direction as
\begin{equation}
\avg{\cdot}_{(i,m)jk} = \frac{1}{2}\left(\left(\cdot\right)_{ijk}+\left(\cdot\right)_{mjk}\right).
\end{equation}

Note, that the split formulation \eqref{eq:entropy-cons_volint} is analogous to that from the linear analysis \eqref{eq:DFLinSharp}; however, the treatment of the metric terms differs. The mapping terms have been ``peeled off'' from the physical flux components. Separating the arithmetic mean of the metric terms from the physical fluxes corresponds to a \textit{dealiasing} of the metric terms, c.f , e.g. \cite{Kopriva:2014yq,kopriva2019free}, (which are variable functions themselves when elements have curved sides) and affects stability. Furthermore, it remains crucial for the discrete entropy estimate that the approximation retains free-stream preservation and that the discrete divergence of the metric terms vanish, as discussed in Sec. \ref{sec:free_stream_AW}.


Substituting the divergence discretizations \eqref{eq:viscVolume} and \eqref{eq:entropy-cons_volint} as well as the BR1 coupling of the viscous fluxes \eqref{eq:BR1Fluxes_with_W} into the strong form DG approximation \eqref{eq:DG_strong_form} gives the \emph{split form DGSEM}:
\begin{equation}\label{eq:DG_split_form}
\begin{aligned}
\iprodN{\IN{J}\statevec{U}_t,\testfuncOne} &+ \iprodN{\DD\bigcontravec F^{\#},\testfuncOne} + \isurfEN\testfuncOne^T\left\{\statevec{F}^*_n - \statevec{F}_n\right\}\hat{s}\dS \\
&= \overRe\iprodN{\DDs\bigcontravec{F}^{v},\testfuncOne} + \overRe \isurfEN\testfuncOne^T\left\{\statevec{F}^{v,*}_n - \statevec{F}^v_n\right\}\hat{s}\dS,\\[0.05cm]
\iprodN{\IN{J}\bigstatevec{Q},\testfuncTwo} &= \isurfEN \statevec{W}^{*,T}\left(\testfuncTwo\cdot\spacevec{n}\right)\hat{s}\dS - \iprodN{\statevec{W},\vecNablaXi\cdot\left(\bigmatrix{M}^T\testfuncTwo\right)}.
\end{aligned}
\end{equation}

As presented, the split formulation \eqref{eq:DG_split_form} is incomplete because the two-point volume fluxes, $\bigstatevec{F}^{\#}$, are not yet defined, and the surface coupling of the advective fluxes through $\statevec{F}^*_n$ remains open. To partially close the question of the surface contributions we connect the choice of the volume flux to the choice of the surface numerical flux through
\begin{equation}\label{eq:vol_connect_surf}
\statevec{F}^*_n = \bigstatevec{F}^{\#}(\statevec{U}_L,\statevec{U}_R)\cdot\vec{n} - \frac{\lambda_{\max}}{2}\jump{\statevec{W}},
\end{equation}
where $\lambda_{\mathrm{max}}$ is an estimate of the fastest wave speed at the point in question. 
We use the local Lax-Friedrichs (LLF) numerical flux function as a blueprint in \eqref{eq:vol_connect_surf} to add numerical dissipation at the element surfaces. This is motivated by the simplicity of the LLF flux and the fact that LLF leads to an entropy stable formulation, as will be shown in Sec. \ref{sec:es_proof_AW}. There are, however, more complex and \textit{selective} dissipation terms (analogous to a \cite{Roe:1981:JCP81} flux) available in the literature, e.g. in \cite{barth99,winters2017uniquely}.

So, the final form of the split form DG approximation now completely hinges on the selection of the two-point numerical volume fluxes $\bigstatevec{F}^{\#}$. From Sec. \ref{sec:yet_another_split_AW}, we know that for linear problems the split formulation is algebraically equivalent to a DG approximation of the advective terms. This equivalence remains true for \eqref{eq:entropy-cons_volint}, as well as being a high-order accurate approximation in the volume \citep{fisher2013,Gassner:2016ye,ranocha2018comparison}. But what is the ``action'' of a particular choice of the numerical volume fluxes? 

\CCLsubsection{The Choice of the Two-Point Flux}\label{sec:sharp_fluxes}

In essence, the split form DG divergence \eqref{eq:entropy-cons_volint} is an \emph{abstraction} or extension of the standard DGSEM approximation. It encompasses the ``classical'' DG divergence operator, but offers an impressive ability to recover discrete approximations of alternative forms of the governing equations. To describe how such high-order discretizations of alternative forms of the advective terms are achieved requires an examination of how components of the numerical volume fluxes are constructed. By assumption, the two-point volume fluxes $\bigstatevec{F}^{\#}$ are symmetric with respect to their arguments. So, it is natural that the components of the numerical volume fluxes will be built from some average state (arithmetic or otherwise).

The abstraction of the split form DGSEM \eqref{eq:DG_split_form} provides a powerful framework that offers a convenient \emph{construct} to generalize, through particular selections of the two-point volume fluxes, a well-trodden technique from the finite difference community for developing split forms of the original governing equations, e.g. \citep{ducros2000,Pirozzoli20107180,kennedy2008,sjogreen2017skew}, which was put it into the nodal DG context by \citet{gassner_skew_burgers,Gassner:2016ye}. We showed in \eqref{eq:split_form_1} and \eqref{eq:split_form_2} that selecting the product of two arithmetic averages is equivalent to a discrete approximation of the split form of a quadratic product, i.e., the average of the conservative form and the advective form of the equations. 


The use of two point fluxes is even more general because it offers a direct translation of various split forms from the continuous level onto the discrete level depending on what is averaged. For example, one can approximate the \textit{cubic split form} for the $x$-momentum flux divergence in the $x-$direction proposed by \citet{kennedy2008} as
\begin{equation}
\resizebox{\textwidth}{!}{$
\begin{aligned}
\frac{1}{4}&\left[(\rho v_1 v_2)_x+\rho (v_1v_2)_x + v_1(\rho v_2)_x + v_2(\rho v_1)_x + v_1v_2(\rho)_x + \rho v_2(v_1)_x + \rho v_1(v_2)_x\right]\\
&\qquad\qquad \approx 2\sum_{m=0}^N\dmat_{im}\avg{\rho}_{(i,m)jk}\avg{v_1}_{(i,m)jk}\avg{v_2}_{(i,m)jk}.
\end{aligned}
$}
\end{equation} 
\citet{Gassner:2016ye} provides a dictionary that one can use in order to immediately construct a discrete split formulation from a proposed continuous splitting.  The nodal split form DGSEM \textit{inherits} the underlying split form of the equations with a high-order spatial accuracy and remains conservative. This is a somewhat surprising result because the split form is created by averaging particular combinations of the conservative and non-conservative forms of the PDEs. See  \citet{gassner_skew_burgers,Gassner:2016ye} for details.

It appears, though, that translation of a given split form to a two-point flux form is only possible if the splitting on the continuous level is explicitly known. This is a problem when we want to get an entropy conserving (or decreasing) approximation. \citet{Tadmor1984428} showed that there always \emph{exists} a split form (also referred to as a skew-symmetric form) that preserves the mathematical entropy of a PDE system for smooth solutions. Unfortunately, the explicit form of such an entropy conservative splitting is \textit{unknown} for many physically relevant and interesting systems of conservation laws like the compressible Euler equations. Therefore, we need an alternative approach to develop numerical approximations that conserve (or dissipate) the mathematical entropy.

Working with finite volume methods, \citet{Tadmor1987_2} developed a condition to guarantee that the numerical flux function is entropy conservative. He eschews any knowledge of the split formulation and focuses, instead, on the contraction of the flux derivative into entropy space \eqref{eq:spatial_ent_contract}, which we restate here in one space dimension for convenience:
\begin{equation}
\label{eq:1d_ent_contract}
\statevec{w}^T\statevec{f}_x = f^S_x.
\end{equation}

The contraction \eqref{eq:1d_ent_contract} relies on the chain rule, whose discrete recovery is extraordinarily difficult, or often impossible, in practice \citep{tadmor2016perfect}. To circumvent this obstacle, one applies the product rule to the entropy contraction \eqref{eq:1d_ent_contract} to re-write it as
\begin{equation}
\label{eq:rewrittenCompatibility}
 \statevec{w}_{x}^T\statevec{f} = (\statevec{w}^T\statevec{f})_x - f^S_{x} = (\statevec{w}^T\statevec{f} - f^S)_x.
\end{equation}
Tadmor analysed this equivalent compatibility condition \eqref{eq:rewrittenCompatibility} to determine a numerical surface flux function for finite volume schemes that is discretely entropy conservative. The  finite volume method takes the unknowns in each element to be mean values that are naturally discontinuous across element interfaces, see, e.g., \citet{leveque2002} for complete details. 

The numerical flux Tadmor derived for finite volume approximations carries over to DG approximations, since, as mentioned in Sec. \ref{sec:curvi_DG}, the idea to resolve discontinuities with numerical surface fluxes is also used in their construction. We describe the result by considering the contraction \eqref{eq:rewrittenCompatibility} at an arbitrary surface point. The flux depends on the discrete values in the current element, denoted with $L$, and the direct neighbour of that element, denoted with $R$. Approximating the derivatives in \eqref{eq:rewrittenCompatibility} with first order differences gives Tadmor's \textit{entropy conservation} condition on the numerical surface flux
\begin{equation}
\left(\frac{\statevec{w}_R - \statevec{w}_L}{\Delta x}\right)^T\statevec{f}^{\ec}(\statevec{u}_L,\statevec{u}_R) = \frac{\left(\statevec{w}_R^T\statevec{f}_R - {f}_R^S\right) - \left(\statevec{w}_L^T\statevec{f}_L - {f}_L^S\right)}{\Delta x}
\end{equation}
where $\Delta x$ is the size of a grid cell. Multiplying through by $\Delta x$ and utilising the jump notation \eqref{eq:jump_operator} the entropy conservation condition on the numerical surface flux is written compactly
\begin{equation}
\label{eq:compact_ec_fv_cond}
\jump{\statevec{w}}^T\statevec{f}^{\ec}(\statevec{u}_L,\statevec{u}_R) = \jump{\statevec{w}^T\statevec{f} - f^S}.
\end{equation}
It is important to reiterate that the entropy conservative flux, $\statevec{f}^{\ec}$, is symmetric in its arguments and consistent to the physical flux in the sense that for identical arguments one recovers the physical flux, i.e., $\statevec{f}^{\ec}(\statevec{u},\statevec{u})=\statevec{f}(\statevec{u})$.

Two interesting aspects of Tadmor's work are: (1) That constructing an entropy conservative surface flux from \eqref{eq:compact_ec_fv_cond} produces a consistent, low-order finite volume approximation \textit{without} the need to solve a Riemann problem and (2) For systems of nonlinear hyperbolic conservation laws \eqref{eq:compact_ec_fv_cond} is a single algebraic condition for a vector of unknown flux quantities. Therefore, there exist many ``solutions'' for $\statevec{f}^{\ec}$ that yield a numerical surface flux that is entropy conservative by satisfying \eqref{eq:compact_ec_fv_cond}. Care must be taken so that entropy conservative numerical flux function remains physically consistent. One such numerical flux, originally proposed by \cite{Tadmor1987_2}, is defined as a phase integral. Though theoretically useful, this phase integral form is computationally impractical. Thus, over the past 20 years \emph{affordable} versions of the entropy conservative finite volume surface flux have been developed for a variety of nonlinear hyperbolic systems, c.f e.g. \cite{Chandrashekar2012,fjordholm2011,Winters:2016cq}.

We provide here a brief summary of one particular affordable numerical surface flux for the $x-$direction of the Euler equations with the ideal gas law, since it is relevant to the development of an entropy stable DG approximation for compressible flows. Complete details are provided by \cite{Chandrashekar2012}. The crucial idea behind finding a numerically tractable version of an entropy conservative surface flux is the evaluation of its components at various means states between $\statevec{u}_L$ and $\statevec{u}_R$. These mean state expressions can take on incredibly complex forms that depend on the arithmetic mean, the product of arithmetic means or more uncommon quantities like the logarithmic mean \citep{carlson1972logarithmic}. We have already introduced notation for the arithmetic mean, e.g. \eqref{eq:BR1_with_U}. The logarithmic mean of two quantities $a_L$ and $a_R$ takes the form
\begin{equation}
a^{\ln} = \frac{a_L - a_R}{\ln(a_L) - \ln(a_R)}.
\end{equation}
Note that care must be taken for the logarithmic mean to remain numerically stable when the states are close, $a_L\approx a_R$, as discussed by \citet{IsmailRoe2009}. Also, we introduce an auxiliary variable proportional to the inverse temperature
\begin{equation}
\beta = \frac{p}{2\rho},
\end{equation}
which simplifies the form of the entropy variables \eqref{eq:entVariables} to
\begin{equation}
\statevec{w} = \left[\frac{\gamma - \varsigma}{\gamma-1} - \beta v^2\,,\,2\beta v_1\,,\,2\beta v_2\,,\,2\beta v_3\,,\,-2\beta\right]^T.
\end{equation}
Then, Tadmor's entropy conservation condition \eqref{eq:compact_ec_fv_cond} and many algebraic manipulations determine an analytical expression of an entropy conservative numerical flux for the compressible Euler equations
\begin{equation}
\label{eq:CH_flux}
\statevec{f}^{\ec}(\statevec{u}_{L},\statevec{u}_{R}) = 
\begin{bmatrix}
 \rho^{\ln}\avg{v_1} \\[0.1cm]
{\rho^{\ln}}\avg{v_1}^2 + \widehat{p} \\[0.1cm]
{\rho^{\ln}}\avg{v_1}\avg{v_2} \\[0.1cm]
{\rho^{\ln}}\avg{v_1}\avg{v_3} \\[0.1cm]
{\rho^{\ln}}\avg{v_1}\widehat{H}
\end{bmatrix},
\end{equation}
where particular average states for the pressure and enthalpy are needed
\begin{equation}
\resizebox{\textwidth}{!}{$
\begin{aligned}
\widehat{p} &= \frac{\avg{\rho}}{2\avg{\beta}},\\[0.1cm] 
\widehat{H} &= \frac{1}{2\beta^{\ln}(\gamma-1)} + \frac{\widehat{p}}{\rho^{\ln}} + \avg{v_1}^2+ \avg{v_2}^2+ \avg{v_3}^2 -\oneHalf\left(\avg{v_1^2} + \avg{v_2^2} + \avg{v_3^2}\right).
\end{aligned}
$}
\end{equation}
The numerical surface flux \eqref{eq:CH_flux} is obviously symmetric with respect to its arguments and it is consistent to the physical flux given in \eqref{eq:NSEquationFluxes_DAK}.

This slight detour of the discussion to numerical fluxes for low-order finite volume methods actually serves as the backbone for the construction of an entropy conservative (or stable) split form DG approximation, for entropy conservative finite volume flux functions are precisely those to be used as the two-point volume flux functions $\bigstatevec{F}^{\#}$ in the  DG flux divergence approximation, \eqref{eq:entropy-cons_volint}. In this way, the two-point divergence approximation can recover the action of the entropy conservative split form \emph{without} an explicit expression of the original equations! The remarkable property of \eqref{eq:entropy-cons_volint} is that it extends the entropy conservative flux form from low-order finite volume approximations to high-order accuracy, as was first demonstrated by \cite{fisher2013}, \textbf{provided} the derivative matrix of the high-order method is a diagonal norm SBP operator, as introduced in Sec. \ref{sec:disc_IP_and_SBP}. This result unlocks the true power of the two-point DG approximation because, in a sense, the entropy analysis of the high-order numerical scheme reduces to the (somewhat simpler) finite volume problem, \eqref{eq:compact_ec_fv_cond}.

We can finally state a complete version of the split form DG approximation that is entropy conservative (or stable) for nonlinear problems. To do so, we take the volume flux functions to be the entropy conservative surface flux functions from the finite volume approximation, e.g. \eqref{eq:CH_flux}, in each Cartesian direction,
\begin{equation}
\bigstatevec{F}^{\#} = \bigstatevec{F}^{\ec} = \left(\statevec{F}^{\ec}_1\,,\,\statevec{F}^{\ec}_2\,,\,\statevec{F}^{\ec}_3\right)^T.
\end{equation}
Note that the method automatically operates on curvilinear geometries because the mapping terms have been separated from the physical flux components in \eqref{eq:entropy-cons_volint}. The final split form DG approximation takes the form
\begin{equation}\label{eq:DG_split_form_with_EC}
\begin{aligned}
\iprodN{\IN{J}\statevec{U}_t,\testfuncOne} &+ \iprodN{\DD\bigcontravec F^{\ec},\testfuncOne} + \isurfEN\testfuncOne^T\left\{\statevec{F}^*_n - \statevec{F}_n\right\}\hat{s}\dS \\
&= \overRe \iprodN{\DDs\bigcontravec{F}^{v},\testfuncOne} + \overRe \isurfEN\testfuncOne^T\left\{\statevec{F}^{v,*}_n - \statevec{F}^v_n\right\}\hat{s}\dS\\[0.05cm]
\iprodN{\IN{J}\bigstatevec{Q},\testfuncTwo} &= \isurfEN \statevec{W}^{*,T}\left(\testfuncTwo\cdot\spacevec{n}\right)\hat{s}\dS - \iprodN{\statevec{W},\vecNablaXi\cdot\left(\bigmatrix{M}^T\testfuncTwo\right)},
\end{aligned}
\end{equation}
where the surface contributions of the advective fluxes have the form \eqref{eq:vol_connect_surf}.


\CCLsubsection{The Boundedness of the Discrete Entropy}\label{sec:es_proof_AW}


The stage is now set to demonstrate semi-discrete entropy stability of the split form DG method \eqref{eq:DG_split_form_with_EC} for the compressible Navier-Stokes equations. 

Fundamentally, the goal of the discrete entropy analysis is to mimic the continuous analysis performed in Sec. \ref{sec:energy_entropy_continuous}. The key to the continuous entropy analysis was two-fold and required: (i) Integration-by-parts \eqref{eq:IntegrationByParts1D2_DAK} and (ii) Proper contraction of the physical fluxes to become the entropy fluxes \eqref{eq:spatial_ent_contract} (essentially the chain rule). We have on hand discrete, high-order equivalents of both necessary components: (i) A derivative matrix $\dmat$ with the SBP property \eqref{eq:SBP1D2_DAK} and (ii) Two-point flux functions that satisfy Tadmor's entropy conservation condition \eqref{eq:compact_ec_fv_cond}, which are lifted to high-order with the split form DG divergence \eqref{eq:entropy-cons_volint}. 

We begin from the final split form DG approximation for nonlinear problems \eqref{eq:DG_split_form_with_EC} and mimic the continuous entropy analysis and closely as possible to get a discrete bound on the entropy. The test function in the first equation is replaced with the polynomial interpolant of the entropy variables $\testfuncOne \gets \statevec{W}$ and the test function in the second equation is replaced with the viscous fluxes, $\testfuncTwo \gets \bigstatevec{F}^v$ to obtain
\begin{equation}\label{eq:DG_split_form_with_EC_proof_1}
\begin{aligned}
\iprodN{\IN{J}\statevec{U}_t,\statevec{W}} &+ \iprodN{\DD\bigcontravec F^{\ec},\statevec{W}} + \isurfEN\statevec{W}^T\left\{\statevec{F}^*_n - \statevec{F}_n\right\}\hat{s}\dS \\
&= \overRe \iprodN{\DDs\bigcontravec{F}^{v},\statevec{W}} + \overRe \isurfEN\statevec{W}^T\left\{\statevec{F}^{v,*}_n - \statevec{F}^v_n\right\}\hat{s}\dS,\\[0.05cm]
\iprodN{\IN{J}\bigstatevec{Q},\bigstatevec{F}^v} &= \isurfEN \statevec{W}^{*,T}\left(\bigstatevec{F}^v\cdot\spacevec{n}\right)\hat{s}\dS - \iprodN{\statevec{W},\vecNablaXi\cdot\left(\bigmatrix{M}^T\bigstatevec{F}^v\right)}.
\end{aligned}
\end{equation}
It is possible to condense the expressions in the second equation of \eqref{eq:DG_split_form_with_EC_proof_1} from previously introduced notation. That is, $\bigcontravec{F}^v = \IN{\bigmatrix{M}^T\bigstatevec{F}^v}$, $\bigstatevec{F}^v\cdot\spacevec{n} = \statevec{F}^v_n$ and the standard DG divergence applied to the viscous flux terms \eqref{eq:viscVolume} so that
\begin{equation}\label{eq:DG_split_form_with_EC_proof_2}
\begin{aligned}
\iprodN{\IN{J}\statevec{U}_t,\statevec{W}} &+ \iprodN{\DD\bigcontravec F^{\ec},\statevec{W}} + \isurfEN\statevec{W}^T\left\{\statevec{F}^*_n - \statevec{F}_n\right\}\hat{s}\dS \\
&= \overRe \iprodN{\DDs\bigcontravec{F}^{v},\statevec{W}} + \overRe \isurfEN\statevec{W}^T\left\{\statevec{F}^{v,*}_n - \statevec{F}^v_n\right\}\hat{s}\dS,\\[0.05cm]
\iprodN{\IN{J}\bigstatevec{Q},\bigstatevec{F}^v} &= \isurfEN \statevec{W}^{*,T}\,\statevec{F}^v_n\,\hat{s}\dS - \iprodN{\statevec{W},\DDs\bigcontravec{F}^v}.
\end{aligned}
\end{equation}

In the continuous entropy analysis, we showed in Sec. \ref{sec:section_NSE_GG} that the volume has no contribution to the entropy estimate because the contraction of the physical flux divergence into entropy space becomes the entropy flux on the boundary. 
The split form DG flux divergence \eqref{eq:entropy-cons_volint} that uses the entropy conservative finite volume fluxes precisely mimics this structure discretely. Therefore, it is possible to replace the volume integral (quadrature) of the advective flux in \eqref{eq:DG_split_form_with_EC_proof_2} by a surface integral (quadrature) as demonstrated by \cite{Gassner_BR1}
\begin{equation}
\label{eq:move_to_boundary}
\iprodN{\DD\bigcontravec F^{\ec},\statevec{W}} = \isurfEN\left(\spacevec{F}^{\ent}\cdot\spacevec{n}\right)\hat{s}\dS = \isurfEN{F}^{\ent}_n\,\hat{s}\dS.
\end{equation}
At its core, the proof of this property relies on the SBP property and discrete metric identities \eqref{eq:discreteMetricIdentities} as well as Tadmor's discrete entropy conservation condition. From \eqref{eq:move_to_boundary}, the split form DG approximation \eqref{eq:DG_split_form_with_EC_proof_2} becomes
\begin{equation}\label{eq:DG_split_form_with_EC_proof_3}
\begin{aligned}
\iprodN{\IN{J}\statevec{U}_t,\statevec{W}} &+ \isurfEN \left({F}^{\ent}_n + \statevec{W}^T\left\{\statevec{F}^*_n - \statevec{F}_n\right\}\right)\hat{s}\dS \\
&= \overRe \iprodN{\DDs\bigcontravec{F}^{v},\statevec{W}} + \overRe \isurfEN\statevec{W}^T\left\{\statevec{F}^{v,*}_n - \statevec{F}^v_n\right\}\hat{s}\dS,\\[0.05cm]
\iprodN{\IN{J}\bigstatevec{Q},\bigstatevec{F}^v} &= \isurfEN \statevec{W}^{*,T}\,\statevec{F}^v_n\,\hat{s}\dS - \iprodN{\statevec{W},\DDs\bigcontravec{F}^v}.
\end{aligned}
\end{equation}


For the compressible Euler equations there are infinitely many convex entropy functions, $s(\statevec{u})$, that symmetrize the equations, as shown by \cite{harten1983}; however, \cite{Dutt:1988} demonstrated that only the entropy function \eqref{eq:entropy_function} simultaneously symmetrizes the advective and viscous components of the compressible Navier-Stokes equations.

With this built-in symmetrization in mind, we next examine the first term of the second equation of \eqref{eq:DG_split_form_with_EC_proof_3}. It is possible to cast the viscous fluxes into an alternative form \eqref{eq:visc_flux_alt} as the gradients of the entropy variables
\begin{equation}
\bigstatevec{F}^v = \bigmatrix{B}^\ent\bigmatrix{M}^T\vecNablaXi\statevec{W} = \bigmatrix{B}^{\ent}\bigstatevec{Q}.
\end{equation}
The viscous flux matrices $\bigmatrix{B}^\ent$ are symmetric positive definite \eqref{eq:spd_viscous_terms} and leads to the manipulation
\begin{equation}
\label{eq:posivitve_visc_contribution}
\iprodN{\IN{J}\bigstatevec{Q},\bigstatevec{F}^v} = \iprodN{\IN{J}\bigstatevec{Q},\bigmatrix{B}^{\ent}\bigstatevec{Q}} \geqslant \min_{E,N}(\IN{J})\iprodN{\bigstatevec{Q},\bigmatrix{B}^{\ent}\bigstatevec{Q}} \geqslant 0,
\end{equation}
provided the interpolant of the element Jacobian is non-negative at the Gauss-Lobatto nodes. Again, see  \cite{Gassner_BR1} for details.

Finally, we substitute the second equation of \eqref{eq:DG_split_form_with_EC_proof_3} into the first and apply the estimate \eqref{eq:posivitve_visc_contribution}. This yields an inequality where the volume contribution of the time derivative term is dictated only through the surface contributions of an element
\begin{equation}\label{eq:DG_split_form_with_EC_proof_4}
\begin{aligned}
\iprodN{\IN{J}\statevec{U}_t,\statevec{W}} &+ \isurfEN\left({F}^{\ent}_n + \statevec{W}^T\left\{\statevec{F}^*_n - \statevec{F}_n\right\}\right)\hat{s}\dS \\
&\leqslant \overRe \isurfEN\left(\statevec{W}^{*,T}\,\statevec{F}^v_n + \statevec{W}^T\left\{\statevec{F}^{v,*}_n - \statevec{F}^v_n\right\}\right)\hat{s}\dS.
\end{aligned}
\end{equation}

We take a moment to interpret the crucial steps that have just occurred to arrive at the expression \eqref{eq:DG_split_form_with_EC_proof_4}. The combination of the discrete entropy analysis and the SBP property allowed us to move the advective and viscous flux contributions out of the volume, where we have no control on its behavior, and onto the element boundary, where we \emph{do} have control through the influence of element neighbours and/or boundary conditions by way of the numerical fluxes. 

This movement of all flux influences onto each element's boundaries is a critical intermediate step to mimic the continuous entropy analysis. Now that we have shown how each element contributes to its \textit{local} entropy we are prepared to examine how the discrete entropy will evolve in time globally over the entire domain.

Under the assumption that the chain rule with respect to differentiation in time holds (semi-discrete analysis), the remaining volume term in \eqref{eq:DG_split_form_with_EC_proof_4} is the time rate of change of the entropy in an element. From the contraction property of the entropy variables \eqref{eq:wsContraction} we see that
\begin{equation}
\begin{split}
\iprodN{\IN{J}\statevec{U}_t,\statevec{W}} &= \sum_{i,j,k,}^N \omega_{ijk}J_{ijk}\statevec{W}_{ijk}^T\frac{d\statevec{U}_{ijk}}{dt}\\ 
&= \sum_{i,j,k}^N\omega_{ijk}J_{ijk}\frac{dS_{ijk}}{dt} = \iprodN{\IN{J}S_t,1}.
\end{split}
\end{equation}
Moreover, we get the total discrete entropy by summing over all elements in the mesh
\begin{equation}
\frac{d}{dt}\overline{S} = \sum_{k=1}^K\iprodN{\IN{J^k}S^k_t,1}.
\end{equation}

Just as in the linear analysis, summing over all elements generates jump terms in the fluxes (advective, viscous and entropy) as well as the entropy variables, whereas the numerical surface flux functions are unique. The physical normal vector $\spacevec{n}$ is defined uniquely at surfaces to point outward from the current element and into its neighbour so that $\spacevec{n} = \spacevec{n}_L = -\spacevec{n}_R$. With all this in mind, we find that the total discrete entropy satisfies the inequality
\begin{equation}\label{eq:totat_ent_est_1}
\begin{aligned}
\frac{d}{dt}\overline{S} &\leqslant \sum\limits_{\interiorfaces} \int\limits_{N } \left\{\jump{\spacevec{F}^{\ent}}\cdot\spacevec{n} + \jump{\statevec{W}}^T\statevec{F}^*_n - \jump{\statevec{W}^T\bigstatevec{F}}\cdot\spacevec{n}\right\}\hat{s}\dS \\
&- \overRe  \sum\limits_{\interiorfaces} \int\limits_{N } \left\{\statevec{W}^{*,T}\jump{\bigstatevec{F}^v}\cdot\spacevec{n} + \jump{\statevec{W}}^T\statevec{F}^{v,*}_n - \jump{\statevec{W}^T\bigstatevec{F}^v}\cdot\spacevec{n}\right\}\hat{s}\dS \\ &+ \PBT,
\end{aligned}
\end{equation}
where $\PBT$ are the physical boundary terms with proper outward pointing normal orientation
\begin{equation}
\label{eq:phys_bndy_terms}
\begin{aligned}
\PBT=&\sum\limits_{\boundaryfaces} \int\limits_{N} - \spacevec{F}^\ent\cdot\spacevec{n} + \overRe \statevec{W}^T\left(\bigstatevec{F}^{v}\cdot\spacevec{n}\right)\dS\\
&+  \sum\limits_{\boundaryfaces} \int\limits_{N}  \statevec{W}^T\statevec{F}_{n}^{*}\dS +   \overRe\sum\limits_{\boundaryfaces} \int\limits_{N}  \statevec{W}^{*,T}\statevec{F}_{v,n} + \statevec{W}^T\statevec{F}^{v,*}_{n} \dS.
\end{aligned}
\end{equation}
Notice that the discrete physical boundary contributions precisely mimic those present in the continuous estimate \eqref{eq:continuous_NSE_entropy_estimate} except for additional dissipation due to the surface fluxes $\statevec{F}_{n}^{*}$, $\statevec{W}^{*}$ and $\statevec{F}^{v,*}_{n}$ evaluated at the boundaries.

We first investigate the contribution from the advective flux terms at each quadrature point on the interior element faces. The advective numerical surface flux was selected to take the form \eqref{eq:vol_connect_surf}. So, the first part of the total discrete entropy estimate \eqref{eq:totat_ent_est_1} will be
\begin{equation}
\label{eq:final_adv_surface}
\resizebox{\textwidth}{!}{$\displaystyle{
\begin{aligned}
\jump{\spacevec{F}^{\ent}}\cdot\spacevec{n} + \jump{\statevec{W}}^T\statevec{F}^*_n - \jump{\statevec{W}^T\bigstatevec{F}}\cdot\spacevec{n} &= \jump{\spacevec{F}^{\ent}}\cdot\spacevec{n} + \jump{\statevec{W}}^T\left(\bigstatevec{F}^{\ec}\cdot\spacevec{n} - \frac{\lambda_{\max}}{2}\jump{\statevec{W}}\right) - \jump{\statevec{W}^T\bigstatevec{F}}\cdot\spacevec{n}\\[0.1cm]
&= \left(\jump{\spacevec{F}^{\ent}} + \jump{\statevec{W}}^T\bigstatevec{F}^{\ec} - \jump{\statevec{W}^T\bigstatevec{F}}\right)\cdot\spacevec{n} - \frac{\lambda_{\max}}{2}\jump{\statevec{W}}^T\jump{\statevec{W}}\\[0.1cm]
&=0 - \frac{\lambda_{\max}}{2}\jump{\statevec{W}}^T\jump{\statevec{W}} \leqslant 0,
\end{aligned}
}$}
\end{equation}
where the terms involving $\bigstatevec{F}^{\ec}$ vanish by construction from the entropy conservation condition \eqref{eq:compact_ec_fv_cond} in each Cartesian direction. Also, we note that dissipation must be introduced in an appropriate fashion to ensure the correct sign. In this instance, we took the LLF-type dissipation in terms of the jump in the entropy variables that leads to a guaranteed negative contribution.

Next, we address how the viscous flux terms contribute at the interior element faces. The BR1 discretization \eqref{eq:BR1Fluxes_with_W} was selected for the numerical surface viscous fluxes so that the second part on the right hand side of \eqref{eq:totat_ent_est_1} becomes
\begin{equation}
\label{eq:final_visc_surface_1}
\resizebox{\textwidth}{!}{$\displaystyle{
\statevec{W}^{*,T}\jump{\bigstatevec{F}^v}\cdot\spacevec{n} + \jump{\statevec{W}}^T\statevec{F}^{v,*}_n - \jump{\statevec{W}^T\bigstatevec{F}^v}\cdot\spacevec{n} = \left(\avg{\statevec{W}}^{T}\jump{\bigstatevec{F}^v} + \jump{\statevec{W}}^T\avg{\bigstatevec{F}^v} - \jump{\statevec{W}^T\bigstatevec{F}^v}\right)\cdot\spacevec{n}
}$}.
\end{equation}
From identity \eqref{eq:KnightsFormerlyKnownAsLemma1},
\begin{equation}
\jump{\statevec{W}^T\bigstatevec{F}^v} = \avg{\statevec{W}}^T\jump{\bigstatevec{F}^v} + \jump{\statevec{W}}^T\avg{\bigstatevec{F}^v},
\end{equation}
we see that the viscous numerical fluxes \eqref{eq:final_visc_surface_1} at the interior faces vanish exactly
\begin{equation}
\label{eq:final_visc_surface_2}
\statevec{W}^{*,T}\jump{\bigstatevec{F}^v}\cdot\spacevec{n} + \jump{\statevec{W}}^T\statevec{F}^{v,*}_n - \jump{\statevec{W}^T\bigstatevec{F}^v}\cdot\spacevec{n} = 0,
\end{equation}
as they did for the linear approximation.
In this sense, the BR1 treatment of the viscous terms is \emph{neutrally stable} for the nonlinear compressible flow problem. 

From \eqref{eq:final_adv_surface}, \eqref{eq:final_visc_surface_2} and \eqref{eq:phys_bndy_terms}, the final discrete entropy evolution statement is
\begin{equation}\label{eq:totat_ent_est_2}
\begin{aligned}
\frac{d}{dt}\overline{S} &\leqslant \sum\limits_{\boundaryfaces} \int\limits_{N} - \spacevec{F}^\ent\cdot\spacevec{n} + \overRe \statevec{W}^T\left(\bigstatevec{F}^{v}\cdot\spacevec{n}\right)\dS\\
&\quad - \sum\limits_{\boundaryfaces} \int\limits_{N}  \statevec{W}^T\bigstatevec{F}^{\ec}\cdot\spacevec{n} - \sum\limits_{\allfaces} \int\limits_{N} \frac{\lambda_{\max}}{2}\jump{\statevec{W}}^T\jump{\statevec{W}}\dS\\
&\quad + \overRe\sum\limits_{\boundaryfaces} \int\limits_{N} \statevec{W}^{*,T}\statevec{F}_{v,n} + \statevec{W}^T\statevec{F}^{v,*}_{n} \dS.
\end{aligned}
\end{equation}
Notice that the dissipation in the advective fluxes has an influence on the entropy estimate at every surface (physical and interior). Furthermore, the choice of these auxiliary physical boundary terms must ensure that their effect is dissipative to guarantee entropy stability. From another point of view, the additional term gives constraints on the boundary fluxes from which to derive stable boundary conditions, as explored by \citet{dalcin2019conservative,hindenlang2019stability}. If we assume that boundary data is given so that the entropy will not increase in time, e.g. \textit{periodic boundary conditions}, then
\begin{equation}\label{eq:totat_ent_est_3}
\frac{d}{dt}\overline{S} \leqslant 0.
\end{equation}
Integrating over the time interval $t\in[0,T]$, we see that
\begin{equation}
\overline{S}(T) \leqslant \overline{S}(0),
\end{equation}
which is a discrete equivalent to the entropy bound given in the continuous analysis \eqref{eq:FinalEntropyBound_{DAK}}. 

\CCLsubsubsection{Summary}

The analysis in this section served to culminate this chapter and describe the components of the entropy stable DGSEM for the compressible Navier-Stokes equations. The approach was systematic to highlight the similarities and differences of the split form DG approximation compared to the ``classical'' DGSEM. The most crucial change was the abstraction of the volume contributions in the discrete DG divergence operator via the use of a \textit{two-point volume numerical flux}. Furthermore, it served to clarify fundamental changes to the numerical approximation when studying solution estimates for linear and nonlinear problems. Most notably were the use of (i) The gradient of the entropy variables as an auxiliary quantity in the viscous components and (ii) Entropy conservative finite volume flux functions in the two-point volume split formulation.


At its heart, the strengths of the discontinuous Galerkin family of methods are its high-order solution accuracy and low dissipation/dispersion errors, e.g., \cite{winters2018compareSplit}. To retain these beneficial properties and expand the DGSEM to be provably entropy stable for nonlinear problems, one borrows some of the strongest features of other numerical methods: 
\begin{itemize}
\item Geometric flexibility from finite element methods.
\item Integration-by-parts (summation-by-parts) from spectral methods.
\item Split formulations of non-linear terms from finite difference methods. 
\item Entropy analysis tools from finite volume methods.
\end{itemize}
All these components were merged to create the nodal split form DG framework \eqref{eq:DG_split_form_with_EC} that can approximate the solution of general, nonlinear advection-diffusion equations.

\CCLsection{Epilogue}\label{sec:epilogue_GG}

We have surveyed the core components of the split form DG framework, a modern nodal DG variant, herein for the linearized and nonlinear compressible Navier-Stokes equations. A key feature of the split form framework is that it provides a demonstrable improvement to the robustness of the high-order numerical approximation, e.g. \citep{Gassner:2016ye}. A further exploration and (partial) explanation of this beneficial property is provided by \citet{Gassner_frontiers} or \citet{winters2018compareSplit}.

The response of the broader high-order numerics community, DG or otherwise, to the split form framework has been immense. As such, it remains an active area of research as the framework is developed and expanded upon in different contexts. For the compressible Navier-Stokes equations, this includes examinations into the development of provably stable boundary conditions \citep{Parsani201588,dalcin2019conservative,hindenlang2019stability} as well as explorations using the split form framework as a ``baseline'' to which turbulence modelling capabilities are added \citep{flad2017,manzanero2020design,flad2020large}.

A principle aspect of the split form DG framework is its generality. In essence, the stability estimates developed in Secs. \ref{sec:linear_split_form} and \ref{sec:es_proof_AW} rely only on:
\begin{enumerate}
\item The SBP property of the derivative matrix $\dmat$.
\item The formulation of a two-point symmetric flux function to be used in the volume and at the surface.
\end{enumerate}
Because of this, the framework is readily extended to other high-order numerical methods that feature the SBP property, e.g. multi-block finite difference methods \citep{Hicken2016,Crean2018} or alternative DG approaches \citep{pazner2019analysis}. Additionally, the split forms have been extended to many other systems of PDEs including:
\begin{itemize}
\item The shallow water equations \citep{Gassner2016291,Niklas-Wintermeyer:2016kq}.
\item Euler equations with alternative equations of state \citep{winters2019entropy}.
\item Ideal \citep{Liu2018} and resistive \citep{Bohm2018} magnetohydrodynamic (MHD) equations.
\item Relativistic Euler and MHD equations \citep{biswas2019entropy,wu2019entropy}.
\item Two-phase flows \citep{Renac2019}.
\item The Cahn-Hilliard equations \citep{manzanero2020free}.
\item Incompressible Navier-Stokes (INS) \citep{manzanero2020INS}.
\item Coupled Cahn-Hilliard and INS \citep{manzanero2020INSCH}.
\end{itemize}

The split form technique described herein was designed for curvilinear unstructured hexahedral meshes. Recent extensions have increased the generality and flexibility of the framework to include meshes that contain simplex elements \citep{chen2017entropy,chan2018}, are non-conforming \citep{Friedrich2018}, or move \citep{krais2020split,schnucke2020entropy,Kopriva2016148}. Furthermore, it is possible to create similar entropy stability estimates for interpolation/quadrature node sets that do not include the boundary points \citep{chan2019efficient}.

The split form technique, through the introduction of the two-point volume fluxes, increases the computational cost of the DG method \textit{locally} on each element. The coupling between DG elements remains weak and the split form DG framework retains the attractive, highly parallelizable nature of the DGSEM \citep{wintermeyer2018entropy}. 

As noted in Sec. \ref{sec:energy_entropy_continuous}, an assumption made within the entropy stability estimate is positivity of particular solution quantities, e.g. the density. For practical simulations, additional shock capturing measure must be incorporated to maintain positivity. However, this must be done carefully to maintain both high-order accuracy and provable entropy stability \citep{HennemannShock}.
Moreover, the analysis in Sec. \ref{sec:es_proof_AW} was done is a semi-discrete sense. Special considerations must be made to develop a fully discrete estimate on the entropy \citep{Friedrich2019,ranocha2020relaxation}.

Overall, exciting developments continue in the realm of high-order DG methods, and split formulations in general, where numerical methods are designed to mimic important continuous stability estimates of PDE systems. Interestingly, the numerical approximations described in this chapter are rapidly approaching the bleeding edge of the current mathematical knowledge we have for the physical models themselves. As it is quite difficult, and perhaps unwise, to discretely mimic physical properties that we do not understand, the further development of modern high-order methods should be done in close collaboration with researchers from physics, computer science and mathematics.

In closing, we took writing this book chapter as an opportunity to provide the interested reader with many details on the mathematical derivations but also of the numerical algorithms with the aim to provide a starting point for an actual implementation. In addition to this book chapter, we refer to the open source code FLUXO\footnote{github.com/project-fluxo} that implements the 3D curvilinear split-form DG methodology with different two-point fluxes. It is our hope that the self-contained derivations and discussions in the previous sections have clarified the motivation and construction of the split form DG method.

%
%
%
%
%
%
%
%
\bibliographystyle{abbrvnat}
\bibliography{final_build}   

\end{document}